\documentclass[12pt]{amsart}\usepackage{a4wide,enumerate,color,graphicx}
\usepackage{amsmath}
\usepackage{scrextend} % for labeling environment
\usepackage{amssymb}
\usepackage{esint}
\usepackage{amsfonts}
\usepackage{amsthm}
\usepackage{enumerate}
\usepackage{mathrsfs}
\usepackage{url,hyperref}
\usepackage{mhchem}     % chemical equations with \ce{}

\newtheorem{lemma}{Lemma}[section]
\newtheorem{prop}[lemma]{Proposition}
\newtheorem{theo}[lemma]{Theorem}
\newtheorem{defin}[lemma]{Definition}
\newtheorem{rem}[lemma]{Remark}
\newtheorem{coro}[lemma]{Corollary}

\DeclareMathOperator{\curl }{curl}
\DeclareMathOperator{\divv }{div}

\DeclareMathOperator{\sign }{sign}

\DeclareMathAlphabet{\mathdutchcal}{U}{dutchcal}{m}{n}

\newcommand{\calb}{\ensuremath\mathdutchcal{b}}

\newcommand{\calm}{\ensuremath\mathdutchcal{m}}

\newcommand{\calR}{\ensuremath\mathdutchcal{R}}
\newcommand{\calQ}{\ensuremath\mathdutchcal{Q}}
\newcommand{\calB}{\ensuremath\mathdutchcal{B}}

\allowdisplaybreaks

\begin{document}

\title[Elliptic system with an energy equation]{A regularity result for quasilinear elliptic systems with an equation of energy--type}

	% Math. Subject Classif.
%\pacs[2008]{82.45.Gj 82.45.Mp 82.60.Lf}				% ggf. Physics Astronomy Classif.

\author[P.-E.~Druet]{Pierre-Etienne Druet}
\address{TU Darmstadt, Department of Mathematics, Schlossgartenstr.\ 7,  64289 Darmstadt, Germany}
\email{druet@mma-tu-darmstadt.de}

\date{\today}

\begin{abstract}
We study the continuity of weak solutions for quasilinear elliptic systems with source terms of critical growth arising from a transport--energy structure. The latter occurs frequently in connection with the first balance principles of physics for internal energy or entropy. Relying on the properties of diffusive measures, we prove a H\"older continuity result for these PDE systems in the case of vectorial $p$-growth and of low regularity features such as mixed boundary conditions on Lipschitz domains. This result is valid in all space dimensions. In the special case where the principal part is linear in the second order, we also study the higher integrability of the gradients and the uniqueness principle in the small. Finally two illustrations in the context of electrothermal coupling are considered, the stationary thermistor model and the Nernst-Planck equations of electrochemistry with variable temperature.
%In particular, the production of internal energy due to friction, compression and electric resistance is locally given in the form of binary products that yield terms of critical growth.
\end{abstract}

\maketitle

% \setcounter{tocdepth}{2}
% \tableofcontents
\section{Introduction}

In this paper, we are concerned with the regularity of weak solutions for systems of partial differential equations of quasilinear elliptic type including terms of critical growth. The mathematical structures concerned occur frequently in connection with the first conservation principles of physics. In particular, the production of internal energy due to friction, compression and electric resistance is locally given in the form of binary products that yield terms of critical growth.

Let $\Omega$ be a bounded domain of $\mathbb{R}^n$ with $n \in \mathbb{N}$. For unknowns $\rho = (\rho_1, \ldots,\rho_N)$  (with $N \in \mathbb{N}$) and $u$ defined on $\Omega$ we consider the following system of second-order quasilinear PDEs
\begin{align}
	& \label{systemcomp} - \divv A_i(x, \, \rho,  \, u, \, \nabla \rho, \, \nabla u) + B_i(x, \, \rho,  \, u, \,  \nabla \rho, \ \nabla u) = \delta_{N+1}^i \, \Pi(x, \, \rho,  \, u, \, \nabla \rho, \, \nabla u)\, , \quad  x \in \Omega %\\
	%& \label{heat} - \divv A_{N+1}(x, \, \rho, \, u, \, \nabla \rho,\,  \nabla u) + B_{N+1}(x, \, \rho, \, u, \, \nabla \rho, \,  \nabla u) = \Pi(x, \, \rho,  \, u, \, \nabla \rho, \, \nabla u)   \, . 
\end{align}
for $i = 1,\ldots,N+1$, with the Kronecker symbol $\delta$. A prototype case arises if the functions $\rho_1,\ldots,\rho_N$ are obeying conservation laws, while the variable $u$ stands for the internal energy density or the temperature. For this reason, we shall call the $N+$first equation in \eqref{systemcomp} the energy equation. 

For $i = 1,\ldots,N+1$ and $j = 1,\ldots,n$ we assume that the coefficient functions $f = A_{ij}$, $f= B_i$ and $f = \Pi$ are defined in $\Omega \times \mathcal{O} \times \mathbb{R}^{(N+1) \times n} =: \Sigma$. Here $\mathcal{O}\subseteq \mathbb{R}^{N+1}$ is an open set accounting for possible state-constraints on the solution vector, as for instance positivity conditions in the common case of $\rho_i$ representing a density or $u$ the temperature.  
While we shall assume that the coefficents $A$ and $B$ satisfy standard $p-$growth conditions, the particular right-hand side $\Pi$ in the energy equation represents a critical term, typically with a product structure: For instance, if we interpret $-A_i$ as flux of the $i^{\rm th}$ physical quantity, and $\nabla \rho_i$ as its driving-force, then the expression $\Pi = A \,:\, \nabla \rho$ would typically occur as a source of energy or entropy in the system, a point that will be explained in more details just hereafter.

On $\partial \Omega$, we shall consider conditions of 
Dirichlet, Neumann or mixed type: 
\begin{align}\label{BC}
	\rho_i = \rho_i^S\quad  \text{ for }\quad  i =1,\ldots,N \, \quad \text{ and } \quad u = u^S \quad \text{ on } \quad S \, ,
\end{align}
corresponding to the cases $S = \emptyset$, $S = \partial \Omega$ and $S \subset \partial \Omega$ being a relatively open subset with positive area. In the two latter cases, $\rho^S_1,\ldots, \rho^S_N, \, u^S$ are given functions on $S$. For simplicity, we shall perform the proofs for homogeneous natural conormal derivative (Neumann boundary conditions) 
\begin{align}\label{BC2}
\nu \cdot A_i(x, \, \rho,  \, u, \, \nabla \rho, \, \nabla u) = 0  \quad \text{ on } \quad \partial \Omega \setminus \overline{S} \, ,
\end{align}
if the latter set possesses positive area. The proofs will make obvious that the case of different $S_i$'s in each equation is not excluded, but we restrict for simplicity to one single Dirichlet part and the same type of boundary conditions in all equations. The domain $\Omega$ is at first allowed to be general in the class $\mathscr{C}^{0,1}$ (local Lipschitz graph, uniform Lipschitz constant). In the case of an interior subset $S \subset \partial \Omega$ its boundary $\partial S$ relatively to $\partial \Omega$ is moreover asked to be a $n-2$ dimensional manifold of class $\mathscr{C}^{0,1}$ (see in particular the assumptions (C1) below). \\

%For the case of pure Neumann or better say co-normal derivative conditions the results are also valid. However, it is necessary to add dissipative processes on the boundary or algebraic gauge conditions for the fields in order to control their average. Such conditions are strongly application-dependent so we prefer not entering the particularities of this case.

{\bf The binary product structure.} In order to motivate the study of the system \eqref{systemcomp}, let us consider the stationary balance of internal energy in a fluid moving with velocity $v$, given in local form as 
\begin{align}\label{U00}
\divv (u \, v + q) = \Pi \, .
\end{align}
On the left-hand side we find under the divergence a quantity called internal energy flux, with the density of internal energy $u$ and the heat flux $q$. The right-hand side $\Pi$ describes the production of internal energy. If we account for compression, friction and input of electromagnetic power, it assumes the form
\begin{align}\label{D00}
\Pi = (-p \, \mathbb{I} + \mathbb{S}) \, : \, \nabla v + J \cdot (E+v\times B)\, 
\end{align}
with the pressure $p$, the viscous stress tensor $\mathbb{S}$, the density $J$ of electrical current, the electric field $E$ and the magnetic induction $B$. 

In general, the equation \eqref{U00} occurs in the description of non-isothermal systems, where the thermodynamic state-variables are the temperature $T$ and the density $\varrho$ of the medium, and the fields $v$, $E$ and $B$ subject, in addition to \eqref{U00}, to the conservation principles for mass, charge, momentum and electromagnetic impulse\footnote{See for instance \cite{MR4452316} for an example in multicomponent fluid dynamics without magnetic field, or \cite{MR2523256} for a single-component problem in magnetohydrodynamics.}. The other quantities $u$, $p$ and $q$, $\mathbb{S}$, $J$ are related to the state-variables by additional constitutive theory. We may assume $u = \hat{u}(T,\varrho)$ and $p = \hat{p}(T,\varrho)$ with certain given functions, while in the simplest, linear case, common choices for $q$, $\mathbb{S}$ and for $J$ is
\begin{gather*}
q = -\kappa \nabla T \quad \text{ -- Fourier--law}\, , \quad 	 \mathbb{S} = 2\eta \, (\nabla v)_{\rm sym} + \lambda \, \divv v \, \mathbb{I} \quad \text{ -- Newton--law}, \\
  J = \sigma \, (E+v\times B) \quad \text{ -- Ohm--law} \, ,
\end{gather*}
with phenomenological coefficients $\kappa$ (heat conductivity), $\eta$ and $\lambda$ (viscosity), and $\sigma$ (electrical conductivity). Then, the equation \eqref{U00} would assume the form
\begin{align*}
\divv (u(T,\, \varrho) \, v -\kappa  \, \nabla T) = -p(T,\, \varrho) \, \divv v \, + \mathbb{S} \, : \, \nabla v + \sigma \, |E+v\times B|^2 \, .
\end{align*}
From the viewpoint of mathematical analysis, the binary products on the right-hand side of \eqref{U00} are difficult to be handled if the main variables $v$, $E$ and $B$ are themselves subject to conservation laws. In fact the quadratic case impeaches in general bounded or H\"older continuous solutions according to the standard monographs (see example (2.33) in chapter 1, paragraph 2 of \cite{MR0244627} for a quasilinear equation, example (0.2) in chapter 8 of \cite{MR0244627} for a quasilinear system, and examples 3.6 of chapter II and 2.1 of chapter V for quasilinear systems in \cite{MR0717034}).\\

{\bf Distributional solution techniques.} In order to handle systems like \eqref{systemcomp} one can exploit the integrability of $\Pi$, which is a physical requirement. This principle allows to resort to the theory of equations with right-hand side $L^1$ (respectively right-hand side measure) and to construct distributional solutions. On the base of \eqref{U00}, the energy variable $u$ or the temperature $T$ are then found in a Sobolev--space $W^{1,\sigma}$ where $\sigma$ is arbitrary in the range $[1, \, n/(n-1)[$, with the space dimension $n \geq 2$. Examples of mathematical investigations using these techniques concern for instance resistance and induction heating problems \cite{MR2599927,MR2837572,MR4348166}, viscous and Ohmic dissipation in electrically or heat conducting fluids \cite{MR2185227,MR1459011,MR2545657,MR2491627}, the production of turbulent kinetic energy in Reynolds averaged Navier--Stokes models \cite{MR1296252,MR1418142,MR2506064}, or the wide field of thermal elasticity \cite{MR2495071}.

The original breakthrough for quasilinear elliptic equations with right-hand side measure is due to \cite{MR0251373}. More generally, for quasilinear equations of the $p-$Laplace type, the theory of elliptic equations with right-hand side in $L^1$ shows that distributional solutions can be constructed in $W^{1,\sigma}$ for $1 \leq \sigma < n(p-1)/(n-1)$, which requires $ p > 2-1/n$ \cite{MR1025884,MR1097910}. For equations with $p-$growth, $p>1$ arbitrary and right-hand side in $L^1$ weak solutions can still be defined via the concepts of entropy solutions or renormalised solutions introduced/used/developed in \cite{MR1201197}, \cite{MR1354907}, \cite{MR1409661}, \cite{MR1760541}.\footnote{There have been some similar developments for quasilinear systems, see \cite{MR1484710}.} Interestingly, for strict monotone constellations and scalar equations, this rather weak notions of solution preserve the uniqueness principle. Notice that solutions are typically constructed as accumulation points of sequences of solutions to the same boundary-value problems with the $L^1-$data replaced by regularised approximations. For $p-$growth equations in the range of $p >2-1/n$, this procedure yields a distributional solution which is identical both with the entropy and the renormalised solution, as shown for instance in \cite{MR1452883}. At least on $C^1$ domains and continuous $x-$ dependence of the main coefficients $A$, we can expect that the distributional solution is itself unique, but this does not need to be true for all constellations of growth exponents, coefficients and domains. A counterexample based on discontinuous coefficients is discussed in the Appendix IV of \cite{MR1354907}.

In fact, for energy systems like \eqref{systemcomp}, the $L^1$ scenario can often be slightly improved. The typical binary product structure $\Pi = A\, : \, \nabla \rho$ is inherently a $\divv-\curl$ structure, which allows for the Hardy--space regularity (for the $L^1_{\log}-$regularity if $\Pi$ is nonnegative) of the right-hand side, \cite{MR1225511}; Or it even helps reaching a pure divergence-form of \eqref{U00} for which the threshold regularity $u \in W^{1,n/(n-1)}$ is attained, e.g.\ Corollaire 3.-- 1) in \cite{MR1225511}, or also \cite{MR2332419}. Moreover, the Gehring--Lemma (reverse H\"older--inequality, see a.o.\ Chapter V of \cite{MR0717034} or Chapter 1.\ of \cite{MR1917320}) or, for linear principal part, the Meyers estimate (\cite{MR0159110}, \cite{MR0990595}, \cite{drunau}) are techniques applying to most of stationary energy-transport systems, so that $u \in W^{1,n/(n-1)+\epsilon}$ can very often be expected. In dimension $n = 2$, uniqueness in the small and the H\"older--continuity are always guaranteed. However, for space-dimensions $n \geq 3$ similar results seem not to be available. That it does not (even in the small) guarantees the validity of the uniqueness principle is a rather severe drawback of this type of solution theory. 
\\

{\bf Continuous weak solutions.} In this paper we discuss sufficient conditions allowing to avoid the low-regularity setting, obtaining solutions which are H\"older continuous and in the quadratic case satisfy the uniqueness principle in the small for higher space-dimensions. From the technical viewpoint, the essential question is whether we can understand the right-hand $\Pi$ in \eqref{systemcomp}, \eqref{U00} as a so-called diffusive measure, \cite{MR1233190,MR1290667,MR1887015,MR2352517}. We show that this is true under certain additional structure conditions. More specifically, our regularity results rely on the following pillars:
\begin{trivlist}\setlength{\itemindent}{9pt}
	\item[$\bullet$] The $p-$growth and ellipticity of the main coefficients;
	\item[$\bullet$] Suitable growth conditions for the coefficients in the $N$ first rows of \eqref{systemcomp} yield H\"older--estimates on $\rho_i$, so that the measure
	\begin{align}\label{dmu}
		d\calm = \Pi(x,w, \,\nabla w) \, dx \, ,
	\end{align}
	satisfies a capacitary estimate on balls: For all $x^0 \in \overline{\Omega}$ and $0<r< {\rm diam}(\Omega)$
	\begin{align}\label{capacitary}
		|\calm|(B_r(x^0) \cap \Omega) \leq c \, r^{n-p + \lambda \, p} \quad \text{ with } \quad 0< \lambda < 1 \text{ and } p > 1 \, .
	\end{align}
	\item[$\bullet$] Using the potential theory of the papers \cite{MR0998128}, \cite{MR1233190},  \cite{MR1290667} or \cite{MR2352517} for scalar equations, the solution $u$ to the last equation on \eqref{systemcomp} possesses the natural Sobolev--regularity and is H\"older--continuous with a certain exponent $\lambda >0$.
\end{trivlist}
Let us note that with a different method but the same framework of diffusive measures, the H\"older continuity and fractional differentiability of solutions was proved for quasilinear elliptic systems with $p$-growth and critical terms in \cite{MR3158841} (see also \cite{MR1944757} for the existence). However this result obviously focusses on other types of applications than the ones here considered. The term with critical growth is required to have a sign and exert some strong coercivity (see the conditions (13) and (14) for the system (17) in \cite{MR3158841}), while we are interested in the case that $\Pi$ does not satisfy sign or monotonicity conditions.

Let us note that right-hand side measures in equations and system that satisfy the diffusivity condition \eqref{capacitary} with $\lambda \leq 0$ are also studied: see among others \cite{MR2352517,MR2719762}. For different structure conditions allowing to handle terms of critical growth in elliptic systems, we refer to \cite{MR3397331}.\\

{\bf Two examples.} In order to explain our idea, we give two examples -- in which there is, symptomatically, no motion:
\begin{trivlist}\setlength{\itemindent}{9pt}
\item[\bf (P1)] {\bf Thermistor problem.} We assume in \eqref{U00}, \eqref{D00} that $v = 0$ and that $q = - \kappa(T) \, \nabla T$ (Fourier--law), where $T$ is the absolute temperature and $\kappa$ a given function. Further we let $E = - \nabla \varphi$ with the electric potential $\varphi$. If $\Omega$ represents an Ohmic conductor, we have $J = \sigma(T) \, E$ with the coefficient $\sigma$ of electrical conductivity, here assumed temperature-dependent. Applying the approximation of electrotechnics for the Maxwell equations, the current density satisfies $\divv J = 0$ and we obtain for the unknown $(\varphi, \, T)$ a system of two PDEs:
\begin{align*}
- \divv (\sigma(T) \, \nabla \varphi) =  0 \quad \text{ and } \quad 
- \divv (\kappa(T) \, \nabla T) =  \sigma(T) \, |\nabla \varphi|^2 + h(x)\, ,
\end{align*}
where $h$ is an applied heat source/sink. This system is often called the (stationary) thermistor-system, see among other \cite{MR0987900,MR3702856,MR4348166}. For this model we shall prove the H\"older continuity of weak solutions and the uniqueness in the small, see Section \ref{Thermite}.
\\

\item[\bf (P2)] {\bf Nernst--Planck equations.} We adopt the same assumptions for $v$ and $q$ as in Problem 1. Assume that instead of an Ohmic conductor we have an electrolyte described by the classical Nernst--Planck equations. Assuming $N-1$ ionic species, the electrical current due to free charges is expressed as
\begin{align*}
 J  = - \sum_{i=1}^{N-1} q_0z_i \,  m_{i} \, \Big(\nabla \frac{\mu_i}{T} + \frac{q_0 z_i}{T} \, \nabla \varphi\Big) \, .
\end{align*}
Here $q_0$ is the Faraday constant and $z_i \in \mathbb{Z}$ the charge multiple carried by the ion, and $m_i$ are the thermodynamic diffusivities of the ions, which are phenomenological functions of the state variables. The chemical potentials (Gibbs variables) $\mu_i$ of the ions, and the electric potential $\varphi$ are to be determined for the conservation principle for the mass and from the Poisson equation. We obtain the PDE system
\begin{gather}
 -\divv \Big[ m_{i} \, \Big(\nabla \frac{\mu_i}{T} + \frac{q_0z_i}{T} \, \nabla \varphi\Big)\Big] = 0 \, \quad \text{ for } \quad i=1,\ldots,N-1 \, \quad \text{ and } \quad    - \divv ( \epsilon \, \nabla \varphi) =  0 \, ,\nonumber\\
\label{npp}- \divv (\kappa \, \nabla T) =   \sum_{i=1}^{N-1} q_0z_i \,   m_{i} \, \Big(\nabla \frac{\mu_i}{T} + \frac{q_0z_i}{T} \, \nabla \varphi\Big) \cdot \nabla \varphi \, .
\end{gather}
\end{trivlist}
This form of the non-isothermal Nernst-Planck equations can be derived from the modelling of the electrolyte as a multicomponent fluid in a thermodynamic consistent way, see \cite{dreyer2013overcoming}. In particular, denoting by $j_i = -m_{i} \, (\nabla (\mu_i/T) + (q_0z_i/T) \, \nabla \varphi)$ the molar mass flux of the $i^{\rm th}$ ion species, the system admits a balance of entropy in the form
\begin{align*}
	\divv\Big(\frac{q - \sum_{i=1}^{N-1} j_i \, \mu_i}{T}\Big) = \sum_{i=1}^{N-1} \frac{|j_i|^2}{m_i} + q \cdot \nabla\frac{1}{T} \, .
\end{align*}
In section \ref{NPPSEC}, we establish for this problem the existence of continuous solutions and a uniqueness result in the small under general state-dependence of the phenomenological coefficients. Note that other known non-isothermal extensions of the Nernst--Planck model, \cite{MR4108365} or drift-diffusion models, \cite{MR1888855,MR2765805,MR3092288} lead to similar structures and problems in analysis. \\
%The fundamental feature that we exploit is that in both Problem (P1) and (P2), the PDEs which are coupled to the energy equation form a nearly diagonal system and admit a H\"older--estimate. A fundamentally different occurs for the case that $v \neq 0$, as the H\"older--continuity result is still not known for the Navier-Stokes equations with temperature-dependent coefficients.
%
%STATEMENTS

{\bf Limitations of the continuity result.} While it applies to important examples like {\bf (P1)} and {\bf (P2)}, the regularity result is at first not valid for fluids with energy or entropy balance with temperature dependent viscosity. Here the H\"older continuity for Navier-Stokes equations with state dependent viscosity is lacking. Hence, in very complex multiphysics systems the notion of distributional weak solutions cannot be improved substantially in high dimensions by the method of the present paper. Let us mention in these contexts also the weaker notion of variational entropy solutions, where the balance of entropy is relaxed to an inequality allowing to avoid defect measures, see for example \cite{MR3916818} and Chapter 2 of \cite{MR3729430}.

Another direction which we do not follow in this paper is overcoming the critical growth via partial decoupling of the equations and/or estimates of higher-order derivatives in smooth domains with pure Dirichlet or Neumann conditions, \cite{MR1446209,MR2837572} or other special features.
\\

{\bf Structure of the paper.} This paper is structured as follows. In the section \ref{ellipt} we formulate our results from the more abstract viewpoint of quasilinear elliptic systems. Then we prove an estimate for the $L^1$ worst-case scenario (Section \ref{Lun}), that we improve in Section \ref{SecHold} by showing the H\"older continuity result. In the Section \ref{ReduSec} we prove higher integrability results for the gradients of weak solution for the case of a linear principal part. Finally, in the Sections \ref{Thermite} and \ref{NPPSEC} we exhibit some consequences for the regularity issue in Problems {\bf (P1)} and {\bf (P2)}.

\section{Quasilinear elliptic system with an equation of energy type.\linebreak
	Notation and results}\label{ellipt}
For unknowns $\rho = (\rho_1, \ldots,\rho_N), \, u$ defined on $\Omega$, with $N \in \mathbb{N}$, we consider the equations \eqref{systemcomp}.
Throughout the paper, we shall use the identification $(\rho_1, \ldots,\rho_N, \, u) = w $ and then $w$ is a map defined on $\Omega$ with values in $\mathbb{R}^{N+1}$. Similarly, $(\nabla \rho_1, \ldots,\nabla \rho_N, \, \nabla u) = \nabla w $ maps from $\Omega$ into $\mathbb{R}^{(N+1) \times n}$. 

\subsection{Notations and basic definitions}

Generic points of the compound $\Sigma := \Omega \times \mathcal{O} \times \mathbb{R}^{(N+1) \times n}$ shall be denoted by $(x, \, w, \, z)$. We say that a map $f: \, \Sigma \rightarrow \mathbb{R}$ satisfies the Caratheodor\`y conditions (with respect to $x$ and $(w, \, z)$) if $x\mapsto f(x, \, w,z)$ is Lebesgue measurable on $\Omega$ for all $(w,z) \in \mathcal{O} \times \mathbb{R}^{(N+1) \times n}$ and $(w,z) \mapsto f(x, \, w,z)$ is continuous on $\mathcal{O} \times \mathbb{R}^{(N+1) \times n}$ for $\mathscr{L}_n$-almost all $x \in \Omega$. Here $\mathscr{L}_n$ denotes the $n$ dimensional Lebesgue measure. For $\mathcal{O} \subseteq \mathbb{R}^{N+1}$ open and ${
\bf p} := (p_1,\ldots,p_{N+1}) \in [1,+\infty]^{N+1}$ we introduce the notation\begin{align*}
 W^{1,{\bf p}}(\Omega; \, \mathcal{O}) := \{w \in W^{1,{\bf p}}(\Omega; \, \mathbb{R}^{N+1}) \, : \, w(x) \in \overline{\mathcal{O}} \text{ for almost all } x \in \Omega \} \, .
\end{align*}
Note that $W^{1,{\bf p}}(\Omega; \, \mathcal{O})$ is a closed subset of $ W^{1,{\bf p}}(\Omega; \, \mathbb{R}^{N+1})$ which is a well known Sobolev space. If $S \subseteq \partial \Omega$ is a $\mathscr{H}_{n-1}$ measurable subset with positive measure, we use the notation $W^{1,p}_S(\Omega)$ for the subspace of $W^{1,p}(\Omega)$ characterised by zero trace on $S$ and correspondingly $W^{1,{\bf p}}_S(\Omega; \, \mathbb{R}^{N+1})$ if the trace of all components vanishes on $S$. Here $\mathscr{H}_{n-1}$ is the standard surface measure.

The coefficients of the quasilinear system are functions of the form $f(x, \, w, \, z)$ defined on $\Sigma$. Growth conditions have to be formulated. Typically so-called natural conditions like $$|f(x,w,z)| \leq f_0(x) + \sum_{i=1}^{N+1} (f_i(x) \, |w_i|^{r_i} + g_i(x) \, |z_i|^{s_i})$$ would be required, with suitable restrictions imposed for the integrability index of the function $f_0$, $f_1,\ldots,f_{N+1}$, $g_1,\ldots,g_{N+1}$ and for the exponents ${\bf r} = (r_1,\ldots,r_{N+1})$ and ${\bf s} = (s_1,\ldots,s_{N+1})$. However, this approach would oblige us to discuss in-depth the properties of Nemicki--operators, introducing a lot of growth parameters and algebraic restrictions. In this paper we only need to know of the coefficients $f(x, \, w, \, z)$ whether their Nemicki--operator transforms $W^{1,{\bf p}}-$ or $L^{\bf q} -$vector fields into $L^r-$functions or $L^{r,\alpha}_{\rm M}-$functions with $r$ and $\alpha$ appropriate. Here $L^{r}(\Omega)$ is the standard Lebesgue space for $1 \leq r \leq +\infty$, while the Morrey--space $L^{r,\alpha}_{\rm M}$ for $1 \leq r \leq +\infty$ and $0 < \alpha \leq n$ is the subspace of all $u \in L^r(\Omega)$ such that
\begin{align}\label{morreycondi}
 [u]_{p,\alpha}^{\rm M} := \sup \Big \{\frac{1}{r^{\alpha}} \, \int_{B_{r}(x^0)\cap \Omega} |u|^{p} \, dx \, : \, x^0 \in \overline{\Omega}, \, r > 0\Big\} < + \infty \, .
\end{align}
We thus introduce the following way of speaking.
\begin{defin}
Let $X(\Omega)$ and $Y(\Omega)$ be normed spaces of measurable functions over $\Omega$ ($X = W^{1,\bf p}$ or $X = L^{\bf q}$, $Y = L^r$ or $Y = L^{r,\alpha}_{\rm M}$). We call $f: \, \Sigma \rightarrow \mathbb{R}$ {\it bounded of class} $X(\Omega) \rightarrow Y(\Omega)$ if $f$ satisfies the Caratheodor\`y conditions with respect to $x$ and $(w, \, z)$ and if for all $M >0$
\begin{align*}
\sup \{ \|f(\cdot, \, w(\cdot), \, \nabla w(\cdot))\|_{X(\Omega)} \, : \,  w(\Omega) \subseteq \mathcal{O}, \, \|w\|_{X(\Omega)}   \leq M\}< +\infty  \, .
\end{align*}
\end{defin}
Next we define weak (distributional) solutions to the system \eqref{systemcomp} in the spirit of elliptic theory for equations with $p-$growth and right-hand side measure.
\begin{defin}\label{solandenergy}
Assume that for all $i=1,\ldots,N+1$ and $j = 1,\ldots,n$, the coefficient functions $f = A_{ij}$, $f = B_i, \, \Pi$ are bounded of class $W^{1,{\bf p}}(\Omega; \, \mathcal{O})\cap L^{\bf q}(\Omega; \, \mathcal{O}) \rightarrow L^1(\Omega)$. We call $ w = (\rho_1, \ldots,\rho_N, \, u) \in W^{1,{\bf p}}(\Omega; \, \mathcal{O})\cap L^{\bf q}(\Omega; \, \mathcal{O}) $ a {\it distributional solution} to \eqref{systemcomp}, \eqref{BC}, \eqref{BC2} if $w = (\rho^S,  \, u^S)$ in the sense of traces on $S$ and the integral identities 
\begin{align}
& \label{transportdistrib} \sum_{i=1}^N\int_{\Omega} \sum_{j=1}^n A_{ij}(x,\rho,u,\nabla \rho,\nabla u) \, \partial_{x_j}\zeta_i + B_i(x,\rho,u,\nabla \rho,\nabla u)\, \zeta_i\,  dx = 0 \, ,\\
& \label{heatdistrib} \int_{\Omega} \sum_{j=1}^n A_{N+1,j}(x,\rho,u,\nabla \rho,\nabla u) \, \partial_{x_j}\zeta_{N+1} + B_{N+1}(x,\rho,u,\nabla \rho,\nabla u) \, \zeta_{N+1} \, dx \nonumber\\
& \quad = \int_{\Omega} \Pi(x,\rho,u,\nabla \rho,\nabla u) \, \zeta_{N+1}\, dx \, ,
\end{align}
are valid for all $\zeta \in C_c^{\infty}(\Omega;\, \mathbb{R}^{N+1})$ or, in the case of mixed and conormal derivative problems, for all $\zeta \in C^{\infty}_c(\overline{\Omega} \setminus \overline{S};\, \mathbb{R}^{N+1})$.

Moreover, we say that the weak solution $w = (\rho, \, u)$ satisfies the \emph{energy dissipation property} if there exist extensions $w^S = (\rho^S, \, u^S) \in W^{1,{\bf p}}(\Omega;\, \mathcal{O})$ of the boundary data to $\Omega$ such that the identities \eqref{transportdistrib} are valid for $\zeta_i = \rho_i - \rho_i^S$, if for all $g \in C^{0,1}(\mathbb{R}) \cap L^{\infty}(\mathbb{R})$ satisfying $g(0) = 0$ and $g^{\prime} \geq 0$, the integral $\int_{\Omega} A_{N+1}(x, \, w,\, \nabla w) \cdot \nabla g(u-u^S) \, dx$ exists, and 
\begin{align}\label{renorm}
 \int_{\Omega} A_{N+1}(x, \, w,\, \nabla w) \cdot \nabla  g(u-u^S) +[B_{N+1}-\Pi](x,\, w,\, \nabla w)   \, g(u-u^S) \, dx \leq 0\, .  
\end{align}
\end{defin}

\subsection{Assumptions and statement of the results} 
In order to formulate the growth conditions, we introduce for $1 < p < \infty$ the conjugate exponent $p^\prime :=p/(p-1)$.  We denote by $p^*$ the Sobolev--embedding exponent for $W^{1,p}$, hence $p^* := np/(n-p)$ for $1 \leq p < n$, $p^* = + \infty$ for $p > n$, and $p^*$ is arbitrary large for $p = n$.
For the system \eqref{systemcomp} a minimal requirement is ellipticity. We moreover assume the $p-$growth of the principal part in each equation. In order to formulate these conditions, we shall assume that the "fluxes" functions $A_i$ are the sum of a monotone, $p_i-$coercive part denoted by $a_i$, and of a lower-order term denoted by $b_i$. We define
\begin{align}\label{qadd}
	A_i = a_{i} + b_{i} \quad  \text{ for }  \quad i=1,\ldots,N+1 \text{ and } (x,w, \, z) \in \Sigma  \, ,
\end{align}
with fields $a_{i}$ and $b_{i}$ defined on $\Sigma$ and subject to the Caratheodor\`{y}--conditions and moreover certain growth conditions formulated just hereafter. Even if higher technicality means that reading is less easy, we at first formulate rather general growth conditions. Notice that illustrations are given later in the sections \ref{Thermite} and \ref{NPPSEC} for {\bf (P1)} and {\bf (P2)}. 

At first, the conditions (A) concern growth of the leading coefficients and their ellipticity:
\begin{labeling}{(A00)}
\item[(A0)] 
For all $i=1,\ldots,N+1$, for almost all $x \in \Omega$, all $(w, \, y) \in \overline{\mathcal{O}}\times  \mathbb{R}^{(N+1) \times n}$ the map $z \mapsto a_i(x, \, w, \, y^1, \ldots,y^{i-1}, \, z, \, y^{i+1}, \ldots, y^N)$ is strict monotone on $\mathbb{R}^n$. 
\item[(A1)] There are $1 < p_1, \ldots,p_{N+1} < n$ 
%and  $2 - \frac{1}{n} < p_{N+1} < +\infty$ 
and for $i = 1,\ldots,N+1$ constants $\mu_i \geq \nu_i > 0$ such that
\begin{align*}
\nu_i \, |z_i|^{p_i} \leq \sum_{j=1}^n a_{ij}(x,w, \, z) \, z_{ij} \leq \mu_i \,  |z_i|^{p_i} \, \quad \text{ for all } \quad (x,w, \, z) \in \Sigma \, . 
\end{align*}
\end{labeling}
Then, we formulate conditions (B) for the lower-order terms and for the sources. For more ease we introduce for parameters $ p > 1$, $k > 1$ and $q \geq 0$ the threshold
\begin{align}\label{thresholds}
	\hat{r}_{\rho}(p,q,k) = \beta \, p^* + \Big(\frac{1}{k^{\prime}}-\beta\Big) \, q \quad \text{ with }\quad 
	\begin{cases} \beta = 0   & \text{ if } q \geq p^* \, ,\\
		0\leq \beta < \min\big\{\frac{1}{k^{'}}, \, \frac{p}{p^*}\big\} \text{ arbitrary}   & \text{ if } 0 \leq q < p^* \, .
	\end{cases} 
\end{align}
Certain of the next conditions (B) are formulated in dependence of a further vector of parameters ${\bf q} = (q_1,\ldots, q_{N+1}) \in [0, \, +\infty]^N$ which shall later stand for \emph{a priori} informations on the integrability of weak solutions.
\begin{labeling}{(A00)}
	\item[(B0)] 
	For all $i=1,\ldots,N$, there is a constant $\omega_i > 0$ such that
	\begin{align*}
		|\Pi(x,w, \, z)|\leq \sum_{i=1}^N \omega_i \,  |z_i|^{p_i} \, \quad \text{ for all } \quad (x,w, \, z) \in \Sigma \, . 
	\end{align*}
\item[(B1$(\bf q)$)] There are $w^S = (\rho^S, \, u^S) \in W^{1,\bf p}(\Omega; \, \mathcal{O})$ extending the boundary data and $\delta_0 > 0$, such that the functions
\begin{align*}%\label{PSI}
	\Psi_i :=& - b_i(x,w,z) \cdot z_i - B_i(x,w,z) \, (w_i-\rho_i^S(x)) + A_i(x,w,z) \cdot \nabla \rho_i^S(x) \, ,\\
	\Psi_{N+1} :=	& \delta_0 \,  \frac{\max\{A_{N+1}(x,w,z)\cdot \nabla u^S(x)-b_{N+1}(x,w,z)\cdot z_{N+1}, \, 0\}}{(1+|w_{N+1}-u^S(x)|)^{1+\delta_0}} \\
	& +\max\Big\{ -B_{N+1}(x,w,z) \, {\rm sign}(w_{N+1}-u^S(x)), \, 0\Big\}\,  
\end{align*}
satisfy for all $(x,w,\, z) \in \Sigma $ the upper bound
\begin{align*} 
	\sum_{i=1}^{N+1} \Psi_i(x,w,\, z) 	\leq & \psi_0(x) + \sum_{i=1}^{N+1} \psi_i(x) \, |w_i|^{r_i} + \sum_{i=1}^N \phi_i(x)\,  |z_i|^{s_i} +  \frac{\phi_{N+1}(x)\,  |z_{N+1}|^{s_{N+1}}}{1+|w_{N+1}-u^S(x)|^t} \, .
\end{align*}
with nonnegative functions $\psi_0 \in L^1(\Omega)$, $\psi_i \in L^{k_i}(\Omega)$ and $\phi_i \in L^{K_i}(\Omega)$ with $k_i\, K_i > 1$ for $i = 1,\ldots,N+1$ and nonnegative $r_i$, $s_i$ and $t$ subject to
\begin{gather*}
	0\leq r_i \leq \hat{r}_{\rho}(p_i,q_i,k_i)\quad  \text{ and }\quad 0 \leq s_i \leq p_i/K_i^{\prime} \quad \text{ for } i =1,\ldots,N\, ,\\
	0\leq r_{N+1} \leq \hat{r}_{u}(p_{N+1},q_{N+1},k_{N+1}) := \frac{1}{p_{N+1}^{\prime}} \, \hat{r}_{\rho}(p_{N+1},\, q_{N+1} \, p_{N+1}^{\prime}, \, k_{N+1})\, ,\\ 
	\begin{cases}
		0 \leq s_{N+1}  \leq p_{N+1}/K_{N+1}^{'} & \text{ if } t > 1/K_{N+1}^{\prime} \, ,\\
		0 \leq s_{N+1} < \gamma\,  p_{N+1}/K_{N+1}^{'} + (1-\gamma) \, t \, p_{N+1}, \, \gamma := \frac{\hat{r}_{u}(p_{N+1},q_{N+1},\infty)}{1+\hat{r}_{u}(p_{N+1},q_{N+1},\infty)} & \text{ otherwise.}
		\end{cases}
\end{gather*} 
% = p_i/(p_i-s_i)$ for $i=1,\ldots,N$ and $k_{N+1} > p_{N+1}-1/(p_{N+1}-1-r_{N+1})$, $K_{N+1} >p_{N+1}-1/(p_{N+1}-1-s_{N+1}+t)$ such that defining 
%
%$r_i,s_i \geq 0$ for $i=1,\ldots,N+1$ and $0\leq t \leq s_{N+1}$ such that  
%\begin{gather*}%\label{lesrs}
%	r_i,s_i < p_i \quad \text{ for } \quad i = 1,\ldots,N \quad \text{ and } \quad  r_{N+1}, \, s_{N+1}-t < p_{N+1} - 1  \, ,
%\end{gather*}
%and 
%	\item[(B2)] There are ${\bf q}\in [0, \, +\infty]^{N+1}$ and $1 \leq \sigma_0 < \frac{\max\{p_{N+1}^*/p^{\prime}_{N+1}, \, q_{N+1}\}}{1+\max\{p_{N+1}^*/p^{\prime}_{N+1}, \, q_{N+1}\} } \, p_{N+1}$ and, for $i=1,\ldots, N+1$, numbers $0< \theta_i \leq 1$ such that $b_i$ is of class $W^{1,(p_1,\ldots,p_N,\sigma_0)} \cap L^{\bf q} \rightarrow L^{p_i^{\prime},\alpha_i}_{\rm M}$ with $\alpha_i = n-p_i+p_i \, \theta_i$ and $B_i$ is of class $W^{1,(p_1,\ldots,p_N,\sigma_0)} \cap L^{\bf q} \rightarrow L^{(p_i^*)^{\prime},\alpha_i}_{\rm M}$ with $\alpha_i = \frac{(p_i^*)^{\prime}}{(p_i)^{\prime}} \, (n-p_i+p_i \, \theta_i)$.
	\item[(B2($({\bf q})$)] There is
	$1 \leq \sigma < \frac{\max\{p_{N+1}^*/p^{\prime}_{N+1}, \, q_{N+1}\}}{1+\max\{p_{N+1}^*/p^{\prime}_{N+1}, \, q_{N+1}\} } \, p_{N+1}$,
	such that $b_{N+1}$ is of class $W^{1,(p_1,\ldots,p_N,\sigma_0)} \cap L^{\bf q} \rightarrow L^{p_{N+1}^{\prime}}$ and $B_{N+1}$ is of class $W^{1,(p_1,\ldots,p_N,\sigma_0)} \cap L^{\bf q} \rightarrow L^{(p_{N+1}^*)^{\prime}}$;
	\item[(B3$({\bf q},i)$)] For $i=1, \ldots, N+1$, the function $b_i$ is of class $W^{1,{\bf p}} \cap L^{\bf q} \rightarrow L^{p_i^{\prime},n-p_i+p_i \, \theta_i^b}_{\rm M}$ and $B_i$ is of class $W^{1,{\bf p}} \cap L^{\bf q}\rightarrow L^{1,\alpha_i}_{\rm M}$ with $\alpha_i =  n-p_i+p_i \, \theta_i^B$.
%	\item[(B3)] For $i = 1, \ldots, N+1$ there is $p_i^S >n$ such that $w_i^S \in W^{1,p_i^S}(\Omega)$.
\end{labeling}
The condition {\rm (B3$({\bf q},i)$)} is weaker than formulating conditions in higher Lebesgue spaces. The proof of the following remark is a well-known exercise on H\"older's inequality.
\begin{rem}\label{morrey}
Assume that $b_i$ is of class $W^{1,{\bf p}} \cap L^{\bf q} \rightarrow L^{p_i^b}$ and $B_i$ is of class $W^{1,{\bf p}} \cap L^{\bf q} \rightarrow L^{p_i^B}$ where the parameters obey:
	%	\begin{itemize}
		%\item for $i \in \{1,\ldots, N\}$: 
		$q_i\in \Big[\max\Big\{p^*_i, \,\frac{n \, \max\{2-p_i,0\}}{1-p_i}\Big\}, \,  +\infty\Big]$ and
		\begin{align*}
			p_i^b > \begin{cases}
				\frac{n}{p_i-1} & \text{ for } p_i \geq 2\\
				\frac{nq_i}{q_i(p_i-1) - n(2-p_i)} & \text{ for } 1 < p_i < 2
			\end{cases}\, , \quad p_i^B > \frac{n}{p_i} \, .
		\end{align*}
		%\item for $i = N+1$: $q_{N+1} \geq p^*_{N+1}$, $p_{N+1}^B > n$ and we set $p_{N+1}^b$ arbitrary finite.
	%\end{itemize} 
	Then we can show that {\rm (B3$({\bf q},i)$)} is valid with $\theta_i^b := 1-n/[p_i^b \, (p_i-1)] >0$ and $\theta_i^B := 1-n/(p_i\, p_i^B) >0$.
	%$\theta_i^B := n(p_i/n-1/p_i^B)/(p_i-1) >0$
\end{rem}

As to the boundary conditions, we always assume at least that $\Omega$ is a domain of class $\mathscr{C}^{0,1}$ and $S \subset \partial \Omega$ is measurable in the sense of the surface measure $\mathscr{H}_{n-1}$, but we shall in places need more assumptions for the boundary data:
\begin{labeling}{(A00)}
	\item[(C1)] If $S \subset \Omega$ and $\Omega \setminus S$ have positive measure (mixed b.-v.\ problem), then $S$ is open and the the relative boundary $\partial S$ in $\partial \Omega$ is a $n-2$ dimensional manifold of class $\mathscr{C}^{0,1}$. Moreover there is $R_0(\partial S) >0$ such that for all $x \in \partial S$, we find a Lipschitz continuous diffeomorphism $T^x_S$ such that $T^x_S(x) = 0$ and such that:
	\begin{enumerate}[(a)]
		\item \label{Sun} the image under $T^x_S$ of $B_{R_0^S}(x) \cap S$ is contained the plane $x_n = 0$;
		\item\label{Sdeux} the image under $T^x_S$ of $B_{R_0^S}(x) \cap \partial \Omega \setminus \overline{S}$ is contained in the plane $x_{n-1} = 0$;
		\item \label{Strois} $T^x_S(\Omega\cap B_{R_0^S}(x)) \subset \{x \in \mathbb{R}^n \, : \, x_{n}, \, x_{n-1} < 0\}$;
		\item \label{Squatre} the transformations are bi-Lipschitzian, uniformly for $x \in \partial S$: 
		\begin{align*}%\label{LM}
			M := \sup_{x\in\partial S} \|T^x_S\|_{C^{0,1}(B_{R_0}(x))}, \, L :=  \sup_{x\in\partial S} \|(T^x_S)^{-1}\|_{C^{0,1}(T^x_S(B_{R_0}(x)))} <+\infty\, .
		\end{align*}  
	\end{enumerate}
%	Also in this case, we assume that
		\item[(C2$(i)$)] $w_i^S \in L^{\infty}(S)$; 
\item[(C3$(i)$)] There is $p_i^S >n$ such that $w_i^S \in W^{1,p_i^S}(\Omega)$. 
\end{labeling}
\begin{rem}\label{stampa}
The additional properties {\rm (C1a) -- (C1d)} are needed for the H\"older continuity of solutions to the mixed boundary-value problem for elliptic equations. These conditions were formulated in \cite{MR0126601}, see the example 13.\ therein. In fact, it is shown there that the assumption of $\partial S$ being $n-2$ dimensional of class $\mathscr{C}^{0,1}$ is sufficient for \eqref{Sun} -- \eqref{Squatre}. For quasilinear equations with quadratic growth of the leading coefficients, the latter assumptions are used in the work \cite{MR2271364}. For more recent studies on H\"older continuity for {\it linear equations} with mixed boundary conditions, in which more general constellations for $S \subset \partial \Omega$ are admitted, we refer to \cite{MR1022556,MR1827089,ter2014holder}.
\end{rem}
%Roughly speaking, it is required that for a neighbourhood of $ x \in \partial S$, the surfaces $S$ and $\partial \Omega\setminus \overline{S}$ are mapped into perpendicular planes. 
%Technically, we require that there is $R_0^{S}(\partial S) >0$ such that for all $x \in \partial S$, we find a Lipschitz continuous diffeomorphism $T^x_S$ such that $T^x_S(x) = 0$ and:
%\begin{enumerate}[(a)]
%	\item \label{Sun} the image under $T^x_S$ of the intersection $B_{R_0^S}(x) \cap \overline{S}$ is contained the plane $x_n = 0$;
%	\item\label{Sdeux} the image under $T^x_S$ of the intersection $B_{R_0^S}(x) \cap \partial \Omega \setminus \overline{S}$ is contained in the plane $x_{n-1} = 0$;
%	\item\label{Strois} $T^x_S(\Omega\cap B_{R_0^S}(x)) \subset \{x_{n}, \, x{n-1} < 0\}$.
%\end{enumerate}
%Also in this case, we assume that the transformation is uniformly bi-Lipschitzian:
%\begin{align}\label{LM}
%	M := \sup_{x\in\partial S} \|T^x_S\|_{B_{R_0^S}(x)}, \, L :=  \sup_{x\in\partial S} \|(T^x_S)^{-1}\|_{C^{0,1}(T^x_S(B_{R_0^S}(x)))} <+\infty\, .
%\end{align}  
%
%We next formulate our main results. First result is the estimate for weak solutions in Sobolev norms.
% ?????
% Moreover, for $ p > 1$, $k > 1$ and $q \geq 0$, we define
%\begin{align}\label{thresholds}
%	\hat{r}_{\rho}(p,q,k) = \beta \, p^* + \Big(\frac{1}{k^{\prime}}-\beta\Big) \, q \quad \text{ with }\quad 
%	\beta =\begin{cases}  0   & \text{ if } q \geq p^* \, ,\\
%		\min\big\{\frac{1}{k^{'}}, \, \frac{p}{p^*}\big\} \text{ arbitrary}   & \text{ if } 0 \leq q < p^* \, .
%	\end{cases} 
%\end{align}

The first main result is an estimate for the H\"older norm of weak solutions. 
%For more ease we introduce for $ p > 1$, $k > 1$ and $q \geq 0$ the threshold
%\begin{align}\label{thresholds}
%	\hat{r}_{\rho}(p,q,k) = \beta \, p^* + \Big(\frac{1}{k^{\prime}}-\beta\Big) \, q \quad \text{ with }\quad 
%	\begin{cases} \beta = 0   & \text{ if } q \geq p^* \, ,\\
%		0\leq \beta < \min\big\{\frac{1}{k^{'}}, \, \frac{p}{p^*}\big\} \text{ arbitrary}   & \text{ if } 0 \leq q < p^* \, .
%	\end{cases} 
%\end{align}
\begin{theo}\label{MAIN}
	Let $\Omega \subset \mathbb{R}^n$ be a bounded domain and $S \subseteq \partial \Omega$ be measurable. Consider a weak solution $w = (\rho,u) \in W^{1,{\bf p}}(\Omega; \, \mathcal{O})$ with the energy dissipation property to \eqref{systemcomp}, \eqref{BC}, \eqref{BC2}. We assume all conditions of {\rm (A)}, {\rm (B$(\bf q)$)} and {\rm (C)} and we let the parameters be subject to 
	${\bf q} \geq {\bf 0}$ if $\mathscr{H}_{n-1}(S) > 0$ and ${\bf q} > {\bf 0}$ if $\mathscr{H}_{n-1}(S) = 0$,
	\begin{align}\label{plower0}
		p_{N+1} \geq \min\{2-1/n,\, 1+1/q_{N+1}\} \quad \text{ and } \quad p_{N+1} \geq \max_{i=1,\ldots, N} p_i\, ,
	\end{align}
%	Moreover let the numbers $r_i$ and $s_i$ in {\rm (B1)} be subject to
%	\begin{gather*}
%		0\leq r_i \leq \hat{r}_{\rho}(p_i,q_i,k_i)\quad  \text{ and }\quad 0 \leq s_i \leq p_i/K_i^{\prime} \quad \text{ for } i =1,\ldots,N\, ,\\
%		0\leq r_{N+1} \leq \hat{r}_{u}(p_{N+1},q_{N+1},k_{N+1}) := \frac{1}{p_{N+1}^{\prime}} \, \hat{r}_{\rho}(p_{N+1},\, q_{N+1} \, p_{N+1}^{\prime}, \, k_{N+1})\, ,\\ \quad 0 \leq s_{N+1}  \leq \gamma\,  p_{N+1}/K_{N+1}^{'} + (1-\gamma) \, t \, p_{N+1} \quad \text{ with } \quad \gamma := \frac{\hat{r}_{u}(p_{N+1},q_{N+1},\infty)}{1+\hat{r}_{u}(p_{N+1},q_{N+1},\infty)} \, .
%	\end{gather*} 
	%moreover let $p_{N+1} \geq \max_{i=1,\ldots, N} p_i$. 
	 Then, there are constants $0< \lambda_1, \ldots,\lambda_{N+1}< 1$ and $M_0> 0$ such that every weak solution $w = (\rho,u) \in W^{1,{\bf p}}(\Omega; \, \mathcal{O}) \cap L^{\bf q}(\Omega; \, \mathcal{O})$ with the energy dissipation property satisfies $\sum_{i=1}^{N+1} \|w_i\|_{C^{\lambda_i}(\overline{\Omega})} \leq M_0$. The constants $\lambda_1,\ldots,\lambda_{N+1}$ and $M_0$ depend only on the data in the conditions {\rm (A)}, {\rm(B)} and {\rm (C)} and on $\|w\|_{L^{\bf q}(\Omega)}$.\footnote{For $0 \leq q < 1$, we write $u \in L^q(\Omega)$ if $\int_{\Omega} |u|^q \, dx$ is finite. For $q = 0$ the condition $u \in L^q(\Omega)$ is empty since $\Omega$ is bounded.}
\end{theo}
Remark that the quantity $\|w\|_{L^{\bf q}}$ is allowed to occur in the bounds. There are many situations in which one gets a control on some $L^q$-norm  for free, like systems with finite mass or finite energy, reaction terms etc.\ If such additional information is not available, one chooses ${\bf q} = {\bf 0}$. For the pure conormal problem (Neumann case) we must however require such information with ${\bf q} > {\bf 0}$ (componentwise).\\

%Note that here is an additional condition on ${\bf q}$ for the conormal problem, here we have to show that we control some $L^{\bf q}$ norm independently.

The second result concerns the situation where the leading terms in the system are linear in the gradient. Here the H\"older continuity of the variables can be combined with linear theory to obtain a gain of integrability for the gradient. We assume that the fields $a_i$, $b_i$ and the functions $B_i$ occurring in \eqref{qadd} satisfy
\begin{labeling}{(A00)}
	\item[(A1')] $a_{i}(x, \, w, \, z) = a_{i}(x, \, w) \, z_i$ on $\Sigma$, with a scalar function $a_i \in C^{\alpha}( \overline{\Omega \times \mathcal{O}})$, $0 <\alpha \leq 1$ for $i=1,\ldots,N+1$;
	\item[(B3'(i))] There are indices $p_i^b, \, p_i^B > 2$ such that $b_{ij}$ is of class $W^{1,\bf 2}\cap L^{\infty}(\Omega;\,\mathcal{O}) \rightarrow L^{p^b_i}(\Omega)$ and $B_i$ is of class $W^{1,2}(\Omega;\,\mathcal{O}) \rightarrow L^{\frac{np_i^B}{n+p_i^B}}(\Omega)$ for $i=1,\ldots,N+1$ and $j=1,\ldots,n$.
\end{labeling}
Moreover, we introduce a further parameter $p_\Omega >2$ depending only on the geometry:
\begin{labeling}{(A00)}
	\item[(C1')] We call $2 < p_\Omega \leq +\infty$ the supremum of all $p$ such that the Laplacian with mixed Dirichlet--Neumann boundary conditions on $(S,\, \partial\Omega \setminus \overline{S})$ as an operator from $W^{1,p}(\Omega)$ into $[W^{1,p^{\prime}}(\Omega)]^*$ is continuously invertible.  
\end{labeling}
The size of $p_\Omega$ for general domains is an important question in elliptic regularity theory. It can be shown for quite general constellations that $p_{\Omega}$ is not attained, that $p_{\Omega} >2$, see \cite{MR3573649} for recent discussions, and that $p_\Omega > 3$ for Lipschitz domains in the case of pure Neumann or Dirichlet conditions, see \cite{MR1331981}, Theorem 0.5 (a) for the Dirichlet case. For investigations of this question in the case of mixed boundary value problems, we refer to \cite{MR1147281}, \cite{MR2378088} and the literature therein.
\begin{theo}\label{reduced}
	We adopt all assumptions of Theorem \ref{MAIN}	Moreover, we let $p_1 = \ldots = p_{N+1} = 2$ and we assume {\rm(A1')} and {\rm (B3')}. With $p_{\Omega}$ from {\rm (C1')}, for $i=1,\ldots,N$, we let 
	\begin{gather*}
	t_i = \begin{cases}
		 \min\{p_i^b, p_i^B,p_i^S\} & \text{ if } \min\{p_i^b, p_i^B,p_i^S\} < p_{\Omega} \, ,\\
		2 \leq s < p_{\Omega} \text{ arbitrary} & \text{ otherwise,} 
		\end{cases}  \qquad t_{\rm min} :=\min_{i=1,\ldots,N} t_i \, ,\\
		t_{N+1} = \begin{cases}
\min\{p_{N+1}^b,\, p_{N+1}^B,p_{N+1}^S,
\frac{nt_{\rm min}}{(2n-t_{\rm min})^+}\}	& \text{ if }\min\{p_{N+1}^b,\, p_{N+1}^B,p_{N+1}^S,\frac{nt_{\rm min}}{(2n-t_{\rm min})^+}\}  < p_{\Omega} \, ,\\
		2 \leq s < p_{\Omega} \text{ arbitrary} & \text{ otherwise.} 
	\end{cases}
	\end{gather*}
	Then there is $M_1 > 0$ exhibiting the same dependencies as $M_0$ in Theorem \ref{MAIN} such that $\sum_{i=1}^{N+1}\|\nabla w_i\|_{L^{t_{i}}(\Omega)} \leq M_1$ for every weak solution such that the energy dissipation property is valid. 
\end{theo}
At last we can state a uniqueness result for regular weak solutions with small data.
\begin{coro}\label{unique}
	Adopting all assumptions of the Theorem \ref{reduced}, suppose moreover that the dependency of all coefficient functions in $(w,\, z)$ is continuously differentiable on $\mathcal{O} \times \mathbb{R}^{{N+1}\times n}$ and that 
\begin{align}\label{Piderivs}
	|\partial_w \Pi(x,w,z)| \leq c_1(w) \, |z|^2 \quad \text{ and } \quad |\partial_z \Pi(x,w,z)| \leq c_2(w) \, |z| \, , 
\end{align}
with continuous functions $c_i$ on $\mathbb{R}^{N+1}$ for $i=1,2$. If $t_{\min} > n$ and the quantity $Q_0 := \|\psi_0\|_{L^1} + \|(\psi_1,\ldots,\psi_{N+1})\|_{L^{\bf k}} + \|(\phi_1,\ldots,\phi_{N+1})\|_{L^{\bf K}}$ with the data of  condition {\rm (B1)} is sufficiently small, then there exists at most one weak solution.
\end{coro}
%NONROUGH PROBLEM
%allowing also, for many situations, for uniqueness results ''in the small''.

Applications of these results concern the existence and uniqueness of regular weak solutions to the thermistor problem {\bf (P1)} and the Nernst--Planck equations {\bf (P2)}. We refer to Lemma \ref{thermiste} for an example of regularity result that we can obtain in the context of {\bf (P1)}. Let us here restrict to formulate a Theorem that we can prove for {\bf (P2)}. We use the N+1 state variables $\rho_i = \mu_i/T$ for $i = 1,\ldots,N-1$, $\rho_N = \varphi$ and $u=T$. The domain of validity of the model is assumed to be 
\begin{align*}
	\mathcal{O} = \mathbb{R}^N \times \mathbb{R}_+ \, .
\end{align*}
Here $\mathbb{R}_+ = ]0, \, +\infty[$. To formulate the result and we assume the complete state dependence of the phenomenological coefficients $m_i$, $\varepsilon$ and $\kappa$.
%
%There are constants $0<\varepsilon_0 \leq \varepsilon_1$ such that 
%\begin{align}\label{atleastthat}
%	\varepsilon_0 < \varepsilon(\rho, \, u) \leq \varepsilon_1 \quad \text{ for all } \quad (\rho,u) \in \mathbb{R}^N \times \mathbb{R}_+\, .
%	%\quad \text{ and } \quad \kappa_0 \leq \kappa(\rho,\, u) \leq \kappa_1 \, .
%\end{align}
%Then we have the following remark which is essential for the analysis. The proof follows from the considerations in Lemma \ref{linfty} and Lemma \ref{contihold}.
%\begin{rem}\label{dissipphi}
%	Let $\varphi^S \in W^{1,p_N^S}(\Omega)$ with $p_N^S >n$, and assume that {\rm (C1)} is valid for the domain $\Omega$. For given measurable $(\rho, \, u)$ the weak solution $\varphi = \rho_N$ to the equation \eqref{poisson} with the mixed boundary conditions \eqref{BC} belongs to $W^{1,2}(\Omega)$ and to $C^{\lambda_N}(\overline{\Omega})$ for a $0 < \lambda_N \leq 1-n/p_N^S$ determined by the space dimension, the domain and the quotient $\varepsilon_1/\varepsilon_0$ only. Moreover, the measure $d\mu_{\varphi} = |\nabla \varphi|^2 \, dx$ satisfies for all $\rho \leq \rho_0={\rm diam}(\Omega)$ and all balls $$\mu_{\varphi}(B_{\rho} \cap \Omega) \leq c \, \|\nabla \varphi^S\|_{L^{p_N^S}}^2 \, \Big(\frac{\rho}{\rho_0}\Big)^{n-2+2\lambda_N} \,  .$$ 
%\end{rem}
%Then we can at first prove a regularity result under a condition on the values of weak solutions. 
\begin{theo}\label{NPPtheo}
Assume that there are constants $0<\varepsilon_0 \leq \varepsilon_1$ such that 
$\varepsilon_0 < \varepsilon(\rho, \, u) \leq \varepsilon_1$ for all $(\rho,u) \in \mathbb{R}^N \times \mathbb{R}_+$. Moreover, we assume that the boundary data $(\rho_1^S \ldots,\rho^S_N, u^S)$ belong to $W^{1,(p, \ldots,p)}(\Omega; \, \mathcal{O})$ with a $p > n$. 
\begin{enumerate}[(a)]
	\item\label{condireguNPP}{\rm  [Conditional regularity:]} If $m_1, \ldots, m_{N-1}$ and $\kappa$ belong to $C^1(\mathbb{R}^N \times \mathbb{R}_+)$, every weak solution $(\rho,u)$ to the equations \eqref{npp} (see also \eqref{ions}, \eqref{poisson}, \eqref{templie} below) with mixed boundary conditions \eqref{BC}, \eqref{BC2} which moreover satisfies 
	\begin{align*}
		\sup_{i=1,\ldots,N} \|\rho_i\|_{L^{\infty}(\Omega)} + \|u\|_{L^{\infty}(\Omega)} < +\infty \quad \text{ and }\quad  T_0 := \inf_{x \in \Omega} u(x) >0 \, .
	\end{align*}
	has all components H\"older continuous over $\overline{\Omega}$. If {\rm (C1')} is valid with $p_\Omega > n$, all components belong to $W^{1,q}(\Omega)$ with $q > n$. There is at most one such weak solution if $\sum_{i=1}^{N-1} \|\nabla \rho_i^S\|_{L^2} + \|\nabla \varphi^S\|_{L^2} + \|\nabla T^S\|_{L^2}$ is sufficiently small.
	\item\label{ExiNPP} {\rm [Existence under a compatibility condition:] } Suppose that $m_i(\rho,u) = \tilde{d}_i(u) \, e^{\rho_i}$, and there are positive constants $d_0^{\prime} \leq d^{\prime}_1$ and $\kappa_0 \leq \kappa_1$ such that $d_0^{\prime} \leq \tilde{d}_i(u) \leq d_1^{\prime}$ and $\kappa_0\leq \kappa(\rho,u) \leq \kappa_1$ for all $(\rho,u) \in \mathbb{R}^N \times \mathbb{R}_+$. Then there exist a positive function $\hat{\varphi}_0 \in C(\mathbb{R}^2_+ )$, non decreasing in the second argument, such that $\liminf_{b\rightarrow +\infty} \hat{\varphi}_0(a, b) = + \infty$ and  $\limsup_{b\rightarrow 0+} \hat{\varphi}_0(a, b) = 0$ for all $a \geq 0$ and if $$\|\nabla \varphi^S\|_{L^p(\Omega)} \leq \hat{\varphi}_0\Big(\sum_{i=1}^{N-1}\|\nabla \varrho_i^S\|_{L^p(\Omega)}, \, \inf_{x \in S} u^S(x) \Big) \, ,$$ the problem  \eqref{npp} with mixed boundary conditions \eqref{BC}, \eqref{BC2} possesses at least one solution with the regularity as in \eqref{condireguNPP}.
\end{enumerate} 
\end{theo}
Let us remark that in part \eqref{ExiNPP}, the compatibility condition for existence means that the applied voltage should not be too wild if the applied temperature is low, or reversely, the applied temperature should be large enough for a given applied voltage. These conditions are needed to guarantee that the temperature remains positive everywhere in the bulk, which is partly an artificial problem\footnote{The model is not valid for too low a temperature. For too large ion concentrations, the assumption that the momentum balance of the electrolyte as a whole is not affected by the motion of the ions becomes problematic.}. 
\begin{rem}\label{nppextend}
For the extension of the model \eqref{npp} to the Nernst--Planck--Poisson case, one considers for the electric potential the inhomogeneous equation
\begin{align*}
-\divv 	(\varepsilon \, \nabla\varphi) = q_0\, \sum_{i=1}^{N-1} z_i \, c_i \, .
\end{align*}
This introduces an additional coupling which is difficult to handle.
It is however readily verified that part \eqref{condireguNPP} of Theorem \ref{NPPtheo} remains valid for this case, see Section \ref{NPPSEC}. The verification of part \eqref{ExiNPP} is not possible in the present context, it would make this investigation too long. For an overview of recent results on the (isothermal) Nernst--Planck--Poisson model, we refer to the introduction of \cite{MR3160531}.
\end{rem}

%Moreover we discuss other examples.
%\subsection{Applications}

\section{Weak estimates for distributional solutions}\label{Lun}

In this section we provide an estimate for the gradients of solutions with minimal regularity.
In many relevant cases, like for instance the $B$ terms being of reaction type, we can hope that the data naturally provide, independently on the regularity of gradients, a control on the solutions in Lebesgue spaces. This is why we introduce parameters $q_1,\ldots,q_{N+1} \geq 0$ and allow a dependency on $\int_{\Omega}|w|^{\bf q} \, dx$ for the bounds.

Since we remain in the framework of solutions with integrable gradients for the $N+$first equation, a lower bound for $p_{N+1}$ is expected (see \cite{MR1025884,MR1097910}). Here we require
\begin{align}\label{plower}
p_{N+1} \geq \min\{2-1/n,\, 1+1/q_{N+1}\} \, ,
\end{align}
showing that an {\it a priori} information on the Lebesgue--regularity of $u$ might contribute to significantly improve the threshold of integrability for the gradient in the last equation. 
%
%For more ease we introduce for $ p > 1$, $k > 1$ and $q \geq 0$ the threshold
%\begin{align}\label{thresholds}
%	\hat{r}_{\rho}(p,q,k) = \beta \, p^* + \Big(\frac{1}{k^{\prime}}-\beta\Big) \, q \quad \text{ with }\quad 
%	\begin{cases} \beta = 0   & \text{ if } q \geq p^* \, ,\\
%		0\leq \beta < \min\big\{\frac{1}{k^{'}}, \, \frac{p}{p^*}\big\} \text{ arbitrary}   & \text{ if } 0 \leq q < p^* \, .
%		\end{cases} 
%\end{align}
\begin{lemma}\label{ENERGYEST}
Let $w = (\rho, \, u) \in W^{1,{\bf p}}(\Omega; \, \mathcal{O})$ be a weak solution for \eqref{systemcomp} with the energy dissipation property. Assume that there are $q_1,\ldots,q_{N+1} \geq 0$ such that $\|w\|_{L^{\bf q}} < +\infty$, and that $p_{N+1}$ and $q_{N+1}$ satisfy \eqref{plower}.
We assume {\rm (A0), (A1)}, {\rm (B0)}, {\rm (B1$(\bf q)$)}. 
%Moreover let the numbers $r_i$ and $s_i$ in {\rm (B1)} be subject to
%\begin{gather*}
%	0\leq r_i \leq \hat{r}_{\rho}(p_i,q_i,k_i)\quad  \text{ and }\quad 0 \leq s_i \leq p_i/K_i^{\prime} \quad \text{ for } i =1,\ldots,N\, ,\\
%	0\leq r_{N+1} \leq \hat{r}_{u}(p_{N+1},q_{N+1},k_{N+1}) := \frac{1}{p_{N+1}^{\prime}} \, \hat{r}_{\rho}(p_{N+1},\, q_{N+1} \, p_{N+1}^{\prime}, \, k_{N+1})\, ,\\ \quad 0 \leq s_{N+1}  \leq \gamma\,  p_{N+1}/K_{N+1}^{'} + (1-\gamma) \, t \, p_{N+1} \quad \text{ with } \quad \gamma := \frac{\hat{r}_{u}(p_{N+1},q_{N+1},\infty)}{1+\hat{r}_{u}(p_{N+1},q_{N+1},\infty)} \, .
%\end{gather*} 
Then for every 
$$1 \leq \sigma < \frac{\max\{p_{N+1}^*/p^{\prime}_{N+1}, \, q_{N+1}\}}{1+\max\{p_{N+1}^*/p^{\prime}_{N+1}, \, q_{N+1}\} } \, p_{N+1} \, ,$$
%1 \leq \sigma < (p_{N+1}-1) \, n/(n-1)
there is a function $C_{\sigma}$ depending on $\Omega$, on the parameters $n$, $N$, ${\bf p}$, ${\bf q}$, ${\bf r}$, ${\bf s}$, $t$, ${\bf k}$, ${\bf K}$ and on $\|\psi_0\|_{L^1(\Omega)} + \sum_{i=1}^{N+1} \|\psi_i\|_{L^{k_i}(\Omega)} + \sum_{i=1}^{N+1} \|\phi_i\|_{L^{K_i}(\Omega)}$ and $\int_{\Omega} |w|^{\bf q} \, dx$, such that $$\sum_{i=1}^N \|\rho_i\|_{W^{1,p_i}(\Omega)} + \|u\|_{W^{1,\sigma}} \leq C_{\sigma}\, .$$ If $q_{N+1} = + \infty$, the same is valid with $\sigma = p_{N+1}$. Moreover $C_{\sigma}$ tends to zero if $Q_0 :=  \|\psi_0\|_{L^1(\Omega)} + \sum_{i=1}^{N+1} \|\psi_i\|_{L^{k_i}(\Omega)} + \sum_{i=1}^{N+1} \|\phi_i\|_{L^{K_i}(\Omega)}$ tends to zero.
\end{lemma}
\begin{proof}
Due to \eqref{qadd}, (A1) and (B0), we can always bound $\Pi$ via
\begin{align}\begin{split}
	\label{CRHS}
	|\Pi(x,w,z)| \leq & \sum_{i=1}^N  \omega_i \, |z_i|^{p_i}\leq \left(\max_{i=1,\ldots,n} \frac{\omega_i}{\nu_i}\right) \, \sum_{i=1}^N \nu_i \, |z_i|^{p_i}\\
	 %\leq & \mu_i \, \gamma_i \, |z_i|^{p_i} + \gamma_i \, |b_i(x,z,w)| \, |z_i| \\
	\leq & L \, \sum_{i=1}^N \, (A_i-b_i)(x,w,z) \cdot z_i \quad \text{ with } \quad L := \max_{i=1,\ldots,n} \frac{\omega_i}{\nu_i}   \, .
\end{split}
\end{align}
According to the Definition \ref{solandenergy}, we might now choose $\zeta_i = \rho_i-\rho^S_i$ in \eqref{transportdistrib} to get
%\begin{gather*}
 %$ \int_{\Omega} A_i(x, \, w, \, \nabla w) \cdot \nabla ()\rho_i + B_i(x, \, w, \, \nabla w) \, \rho_i\,  dx = 0$ and, 
 in connection with the growth condition (A1) that
 \begin{align}\label{deux}
& \int_{\Omega} \nu_i \, |\nabla \rho_i|^{p_i} \, dx \leq  \int_{\Omega}  (A_i-b_i)(x, \, w, \, \nabla w) \cdot \nabla \rho_i \, dx \nonumber\\
= & \int_{\Omega}  A_i(x, \, w, \, \nabla w) \cdot \nabla (\rho_i-\rho_i^S) +A_i(x, \, w, \, \nabla w) \cdot \nabla \rho_i^S-b_i(x, \, w, \, \nabla w) \cdot \nabla \rho_i\, dx
 \nonumber\\
= &
- \int_{\Omega} B_i(x, \, w, \, \nabla w) \, (\rho_i-\rho_i^S) - A_i(x, \, w, \, \nabla w) \cdot \nabla \rho_i^S +b_i(x, \, w, \, \nabla w) \cdot \nabla \rho_i\,  dx \, .
 \end{align}
For any $0 < \delta < \min\{p_{N+1} - 1,\delta_0\}$, we next choose $ g(u-u^S(x)) := \sign(u-u^S(x)) \, (1 - 1/(1+|u-u^S(x)|)^{\delta})$ in the inequality \eqref{renorm}. Use of $|g| \leq 1$, of \eqref{CRHS} and of \eqref{deux} yields
\begin{align}\label{unweak}
& \int_{\Omega} g^\prime(u-u^S)\,  \Big([A_{N+1}-b_{N+1}] \cdot \nabla u \Big)  \, dx\\
 \leq &\int_{\Omega}g^{\prime}(u-u^S) \,  (A_{N+1}\cdot \nabla u^S-b_{N+1}\cdot\nabla u) +  (\Pi-B_{N+1}) \, g(u-u^S)  \, dx \nonumber\\
\leq & \int_{\Omega}g^{\prime}(u-u^S) \,(A_{N+1}\cdot \nabla u^S-b_{N+1}\cdot\nabla u) -B_{N+1}\, g(u-u^S)   + L \, \sum_{i=1}^N  \, (A_i-b_i) \cdot \nabla\rho_i  \, dx\nonumber\\
  =& \int_{\Omega} g^{\prime}(u-u^S) \,(A_{N+1}\cdot \nabla u^S-b_{N+1}\cdot\nabla u) -B_{N+1}\, g(u-u^S) \, dx \nonumber\\
 &  -  L\, \sum_{i=1}^N  \int_{\Omega}b_i\cdot\nabla\rho_i + B_i \, (\rho_i-\rho_i^S) - A_i\cdot\nabla \rho_i^S \, dx \,. \nonumber
\end{align}
%\leq \int_{\Omega} \delta \,   \Big(\frac{1}{1+|u|}\Big)^{1+\delta} \, \beta_h + \beta - \sum_{i=1}^NB_i \, \rho_i \,  dx\, . \nonumber
We sum up over $i=1,\ldots,N$ in \eqref{deux} and we add the result to \eqref{unweak}. Using the functions $\Psi_1,\ldots,\Psi_{N+1}$ of (B1) yields
\begin{align}\label{moule}
& \int_{\Omega} \sum_{i=1}^N \nu_i \,  |\nabla \rho_i|^{p_i} + \nu_{N+1}\,  g^\prime(u-u^S)\, |\nabla u|^{p_{N+1}} \, dx 
%& \leq  \int_{\Omega} \delta \,   \Big(\frac{1}{1+|u|}\Big)^{1+\delta} \, \beta_h + \, 2 \, \Big( \beta - \sum_{i=1}^N  B_i \, \rho_i\Big) - (h+\beta)  \,   \sign(u) \, \Big(1 - \frac{1}{(1+|u|)^{\delta}}\Big) \,  dx \nonumber \\
 \leq (L+1)\, \int_{\Omega} \sum_{i=1}^{N+1}\Psi_i(x, \, w, \, \nabla w) \,dx \leq \nonumber\\
 & c\, \int_{\Omega}\Big\{ \psi_0(x) +  \sum_{i=1}^N (\psi_i(x) \, |\rho_i|^{r_i} + \phi_i(x) \, |\nabla \rho_i|^{s_i}) + \psi_{N+1}(x) \, |u|^{r_{N+1}} +  \frac{\phi_{N+1}(x)\,|\nabla u|^{s_{N+1}}}{(1+|u-u^S|)^{t}}\Big\} \, dx\, .
\end{align}
Now, with $\tilde{\psi}_0 = \psi_0 + |\nabla u^S|^{p_{N+1}} +  \psi_{N+1} \, |u^S|^{r_{N+1}} + \phi_{N+1}\,|\nabla u^S|^{s_{N+1}}$ which belongs to $L^1$, we can obtain from the latter with $v = u-u^S$
\begin{align*}
	& \int_{\Omega} \sum_{i=1}^N \nu_i \,  |\nabla \rho_i|^{p_i} + \nu_{N+1}\,  g^\prime(v)\, |\nabla v|^{p_{N+1}} \, dx\\ 
	\leq 
	& c_1\, \int_{\Omega}\Big\{ \tilde{\psi}_0(x) +  \sum_{i=1}^N (\psi_i(x) \, |\rho_i|^{r_i} + \phi_i(x) \, |\nabla \rho_i|^{s_i}) + \psi_{N+1}(x) \, |v|^{r_{N+1}} +  \frac{\phi_{N+1}(x)\,|\nabla v|^{s_{N+1}}}{(1+|v|)^{t}}\Big\} \, dx\, .
\end{align*}
Hence, we loose no generality in assuming that $u^S = 0$ in \eqref{moule}, and we will finish the proof under this simplifying assumption.

Next we define $\alpha= \alpha(\delta) := 1-(1+\delta)/p_{N+1} $. Since we choose $0 < \delta < p_{N+1} -1$, it follows that $0 < \alpha < 1$. We then introduce the function 
\begin{align*}%\label{Tildeuprops0}
	\tilde{u} := \frac{1}{\alpha} \, \sign(u) \, \big((1+|u|)^{\alpha} - 1\big) \, ,
\end{align*}
which satisfies 
\begin{align}\label{Tildeuprops}
	|u|^{\alpha} \leq |\tilde{u}| \leq \alpha^{-1} \, |u|^{\alpha}  \quad \text{ and } \quad  |\nabla \tilde{u}|^{p_{N+1}} = \Big(\frac{1}{1+|u|} \Big)^{1+\delta} \, |\nabla u|^{p_{N+1}} \, .
\end{align}
Then, with $\tilde{w}_i = \rho_i$, $\tilde{\nu}_i = \nu_i$ and $\tilde{r}_i = r_i$ for $i=1,\ldots,N$ and with $\tilde{w}_{N+1} = \tilde{u}$, $\tilde{\nu}_{N+1} = \delta\nu_{N+1}$ and $\tilde{r}_{N+1} = r_{N+1}/\alpha$, the inequality \eqref{moule} implies that 
\begin{align}\label{moule2}
	& \int_{\Omega} \sum_{i=1}^{N+1} \tilde{\nu}_i \,  |\nabla \tilde{w}_i|^{p_i} \, dx 
	\leq c\, \int_{\Omega} \psi_0(x) + \sum_{i=1}^{N+1} \psi_i(x) \, |\tilde{w}_i|^{\tilde{r}_i}\,  dx \nonumber\\
	& + c\, \int_{\Omega}\sum_{i=1}^{N} \phi_i(x) \, |\nabla \tilde{w}_i|^{s_i} +\phi_{N+1}(x)\, |\nabla \tilde{w}_{N+1}|^{s_{N+1}} \, (1+|\tilde{w}_{N+1}|)^{\frac{s_{N+1} \, (1-\alpha) - t}{\alpha}} \, dx \,  .
\end{align}
From now the proof reduces to the technical verification that the exponents can be chosen as claimed, see section \ref{append}.
\end{proof}

\section{Hoelder continuity}\label{SecHold}

We will at first discuss the question how to obtain a H\"older bound for the weak solution $\rho = (\rho_1,\ldots,\rho_N)$ to the $N$ first equations in \eqref{systemcomp}. Here we shall basically rely on the regularity theory for equations. In order to apply these or similar results, the system of PDEs for $\rho_1,\ldots,\rho_N$ has to be nearly diagonal in the main order.
At second we discuss the H\"older continuity of the energy variable $u$. Here we shall follow closely the papers \cite{MR0998128}, \cite{MR1233190},  \cite{MR1290667} and \cite{MR2352517}, more particularly Lemma 3.2 and the proof of Theorem 1.14 in \cite{MR1887015}.

\subsection{H\"older continuity for the $N$ first variables}
Let us emphasise that the proof of H\"older continuity for the $N$ first rows of the system ought to be fairly standard under the imposed assumptions. However, we were not able to find references treating at the same time the following three features: Coefficients with $p-$growth, non-homogenous mixed boundary conditions and lower--order terms in Morrey spaces. Hence we prefer discussing the arguments in some details.

We introduce the following notations. For a weak solution $w = (\rho, \, u)$ to \eqref{systemcomp} and for $f: \, \Sigma \rightarrow \mathbb{R}$, $x\in \Omega$ and $z \in \mathbb{R}^n$ and $i = 1,\ldots,N$ we adopt the notation
\begin{align}\label{tildetrafo}
	\tilde{f}_i(x, \, z) = f\Big(x, \,  \rho(x),\, u(x), \, \sum_{j\neq i} \nabla \rho_j(x) \, e^j + \nabla u(x) \, e^{N+1} + z \, e^i\Big) \, .
\end{align}
Then, $\rho_i \in W^{1,p_i} \cap L^{q_i}$ also solves in the weak sense
\begin{gather}
	\label{rhoihoelder1}	- \divv (\tilde{a}_i(x,\nabla \rho_i) + b_i) + B_i = 0 \quad \text{ in } \quad \Omega\, ,\\
	\label{rhoihoelder2}	\rho_i = \rho_i^S \quad \text{ on } \quad S, \qquad \nu(x) \cdot (\tilde{a}_i(x,\nabla \rho_i) + b_i) = 0 \quad \text{ on } \quad \partial \Omega\setminus S \, ,
\end{gather}
where $b_i(x) = b_i(x,w(x),\nabla w(x))$, $B_i(x) = B_i(x,w(x),\nabla w(x))$, and $\nu(x)$ is a unit normal field at $x\in\partial \Omega$. We begin with proving $L^{\infty}-$bounds. 
\begin{lemma}\label{linfty}
	We adopt the assumptions of Lemma \ref{ENERGYEST} and moreover {\rm (B3,$({\bf q},1)$)}, $\ldots$, {\rm (B3,$({\bf q},N)$)} and {\rm (C2$(1)$), $\ldots$, (C2$(N)$)}. Then every weak solution $w = (\rho,u) \in W^{1,\bf p}(\Omega;\, \mathcal{O}) \cap  L^{\bf q}(\Omega;\, \mathcal{O})$ to \eqref{systemcomp} with the energy dissipation property satisfies $\max_{i=1,\ldots,N} \|\rho_i\|_{L^{\infty}(\Omega)} \leq M$,  where $M$ depends on the solution only via $\|w\|_{L^{\bf q}}$.
\end{lemma}
\begin{proof}
	 Since the weak estimates of Lemma \ref{ENERGYEST} apply, the condition (B3,$({\bf q},i)$) implies for $i = 1,\ldots,N$ that $b_i \in L^{p_i^{\prime},\alpha_i^b}_M$ and $B_i \in L^{1,\alpha_i^B}_M$ with the norm controlled independently of the weak solution. We shall at first prove an upper bound for $\rho_i$. In the weak form \eqref{rhoihoelder1} the testfunction $\rho_{i,L} = \max\{0,\rho_i-L\}$ with $L > \sup_{x \in S} \rho_i^S(x)$. For $\Omega_{L} := \{x \, : \, \rho_i(x) > L\}$ with characteristic function denoted by $\chi_{\Omega_L}$, we obtain that
	\begin{align}\label{laberling}
		\nu_i \, \int_{\Omega_{L}} |\nabla \rho_i|^{p_i} \, dx  \leq & \int_{\Omega_L} (A_i- b_i)(x, \, w, \, \nabla w) \cdot \nabla \rho_i \, dx = \int_{\Omega} (A_i(x, \, w, \, \nabla w)-b_i) \cdot \nabla \rho_{i,L} \, dx\nonumber\\
		%& \quad = \int_{\Omega} \Big(a(x, \, w, \, \nabla w) \,:\, \ \nabla \rho - h(x,w,\nabla w)\Big) \, u_L + \beta_h(x, \, w, \, \nabla w) \, \chi_{\Omega_L} \, dx \, \\
		= & - \int_{\Omega_{L}} b_i \cdot \nabla \rho_i + B_i \,  \, \rho_{i,L}  \, dx \,.
	\end{align}
	It also easily follows that
	\begin{align}\label{laberlingpr} 
		\frac{\nu_i}{2} \, \int_{\Omega_{L}} |\nabla \rho_i|^{p_i} \, dx  \leq c \, \int_{\Omega_{L}} |b_i|^{p_i^\prime}  \, dx + \int_{\Omega_L} |B_i| \, \rho_{i,L} \, dx \,.
	\end{align} 
	We consider the measures $d\calb = |b_i|^{p_i^\prime} \, dx$ and $d \calB = |B_i| \ dx$. Using the assumptions (B3,$i$) and the definitions of $\theta_i^b$ and $\theta_i^B$ contained therein, we can verify on every ball $B_{\rho}$ contained in $\Omega$ that
	\begin{align}
		\calb(B_{\rho}) \leq & [b_i]_{p_i^{\prime},\alpha_i^b}^{M} \, \rho^{n-p_i+\theta_i^bp_i} \quad \text{ and } \quad 
		\calB(B_{\rho}) \leq  [B_i]_{1,\alpha_i^B}^M \, \rho^{n-p_i + \theta_i^B \, p_i} \, .\label{Brate}
		%\calB(B_{\rho}) \leq & \|B_i\|_{L^{(p_i^*)^\prime}(B_{\rho})} \, [\lambda_n(B_{\rho})]^{1-\frac{1}{(p_i^*)^{\prime}}} \leq [B_i]_{(p_i^*)^\prime,\alpha_i^B} \, \rho^{\frac{\alpha_i^B}{(p^*_i)^{\prime}}} \, (\omega_1)^{1-\frac{1}{(p_i^*)^{\prime}}} \, \rho^{\frac{n}{p_i}-1}\nonumber\\
		%\leq & c \, [B_i]_{(p_i^*)^\prime,\alpha_i^B} \, \rho^{n-p_i + \theta_i^B \, (p_i-1)}
	\end{align}
	Here $\omega_1$ is the measure of the unit ball. With $\kappa_i^B = 1+\theta_i^B(p_i-1)/(n-p_i)$, we applying the measure-Sobolev--embedding theorem of Lemma \ref{lieberman} (see also Lemma 1.1 in \cite{MR1233190} and \cite{MR0287301}, Th.\ 1 for the original result\footnote{We follow \cite{MR1233190} for this attribution.}), and we obtain that
	\begin{align}\label{laberling12}
		\int_{\Omega_L} |B_i| \, \rho_{i,L} \, dx = \int_{\Omega_L}  \rho_{i,L} \, d\calB \leq \|\rho_{i,L}\|_{L^{\kappa_i^B \, p_i}(\Omega;d\calB)} \, [\calB(\Omega_L)]^{1-\frac{1}{\kappa_i^B \, p_i}} \leq c \, \|\nabla \rho_{i,L}\|_{L^{p_i}} \, [\calB(\Omega_L)]^{1-\frac{1}{\kappa_i \, p_i}} \,.
	\end{align}
	Hence, by means of Young's inequality, \eqref{laberling} further implies that
	\begin{align}\label{laberling2} 
		\frac{\nu_i}{4} \, \int_{\Omega_{L}} |\nabla \rho_i|^{p_i} \, dx  \leq c \, \int_{\Omega_{L}} |b_i|^{p_i^\prime}  \, dx + c  [\calB(\Omega_L)]^{(1-\frac{1}{\kappa_i^B \, p_i})p_i^{\prime}} \,.
	\end{align}
	Therefore, combining \eqref{laberling12} and \eqref{laberling2}, we get
	\begin{align*}
		\int_{\Omega_L}  \rho_{i,L} \, d\calB \leq \|\rho_{i,L}\|_{L^{\kappa_i^B \, p_i}} \, [\calB(\Omega_L)]^{1-\frac{1}{\kappa_i^B \, p_i}} \leq c \, \Big(\calb(\Omega_{L}) + [\calB(\Omega_L)]^{(1-\frac{1}{\kappa_i^B \, p_i})p_i^{\prime}}\Big)^{\frac{1}{p_i}} \, [\calB(\Omega_L)]^{1-\frac{1}{\kappa_i^B \, p_i}} \, .
	\end{align*}
	With similar arguments, with $\kappa_i^b = 1 + \theta_i^b \, p_i/(n-p_i)$, we see that
	\begin{align*}
		\int_{\Omega_L}  \rho_{i,L} \, d\calb \leq \|\rho_{i,L}\|_{L^{\kappa_i^b \, p_i}(\Omega;d\calb)} \, [\calb(\Omega_L)]^{1-\frac{1}{\kappa_i^b \, p_i}} \leq c \, \Big(\calb(\Omega_{L}) + [\calB(\Omega_L)]^{(1-\frac{1}{\kappa_i^B \, p_i})p_i^{\prime}}\Big)^{\frac{1}{p_i}} \, [\calb(\Omega_L)]^{1-\frac{1}{\kappa_i^B \, p_i}} \, .
	\end{align*}
	Using now that both $\kappa_i^B > 1$ and $\kappa_i^b > 1$, we can deduce
	\begin{align}\label{grundlage}
		\int_{\Omega_L}  \rho_{i,L} \, d(\calb+\calB) \leq [(\calb + \calB)(\Omega_L)]^{\sigma} \, ,\quad \sigma = 1+\frac{1}{p_i}-\frac{1}{\min\{\kappa_i^B,\kappa_i^b\} \, p_i} > 1 \, .
	\end{align}
	An upper bound for $\rho_i$ follows using the standard argument of Lemma 5.1 in \cite{MR0244627}, chapter 2. A lower-bound for $\rho_i$ is proven similarly.
	%The $L^{\infty}-$bound for $\rho_i$ is induced from well-known arguments, (cf.\ \cite{MR0l244627}, chapter 3, paragraph 13 for $p_i = 2$, and chapter 4, paragraph 7 for the general nonlinear case) in the form
	Overall
	\begin{align}\label{rhoibounded}
		\|\rho_i\|_{L^{\infty}(\Omega)} \leq & \|\rho_i^S\|_{L^{\infty}(S)} + c(\Omega, \alpha_i^b,\alpha_i^B,n,p_i,\|b_i\|_{L^{p_i,\alpha_i^b}_M},\|B_i\|_{L^{1,\alpha_i^B}_M}) \nonumber \\
		\leq & \|\rho_i^S\|_{L^{\infty}(S)} + c(\Omega, \alpha_i^b,\alpha_i^B,n,p_i,\|w\|_{W^{1,{\bf p}}} + \|w\|_{L^{\bf q}}) \, .
	\end{align}
\end{proof}
Next we discuss the H\"older continuity. Starting from \eqref{rhoihoelder1}, \eqref{rhoihoelder2}, we first homogenise the Dirichlet boundary conditions by defining the function $v^0(x) = \rho_i(x)- \rho_i^S(x) \in W^{1,p_i}(\Omega)$ which solves
\begin{gather}\label{0auxivhier}
	-\divv (a^0_i(x, \, \nabla v^0 ) + b_i) + B_i = 0 \text{ in } \Omega \,, \\  
	v^0 = 0 \text{ on } S \, ,\qquad
	\nu(x) \cdot (a^0_i(x, \, \nabla v^0) + b_i) = 0 \text{ on } \partial \Omega \setminus S\, ,\nonumber
\end{gather}
with the new vector field $a_i^0(x,z) := \widetilde{a}_i(x, \, z + \nabla \rho_i^S(x))$. (A1) implies that
\begin{align*}
	\nu_i \, |z + \nabla \rho_i^S(x)|^{p_i} \leq a_i^0(x, \, z) \cdot (z + \nabla \rho_i^S(x)) \leq 	 a_i^0(x, \, z) \cdot z + \mu_i \, |z + \nabla \rho_i^S(x)|^{p_i-1} \, |\nabla \rho_i^S(x)| \, .
\end{align*}	
Hence with the help of elementary inequalities we find a constant $C_0 = C_0(\nu_i,\mu_i,p_i) >0$ such that
\begin{align}\label{a1weakenedhere}
	& \frac{\nu_i}{2} \, |z|^{p_i} - C_0 \, |\nabla \rho_i^S(x)|^{p_i} \leq a_i^0(x, \, z) \cdot z\, , \\
	\label{a1weakenedhereupper}	& |a_i^0(x, \, z)|	\leq \mu_i \, (|z|^{p_i-1} +|\nabla \rho_i^S(x)|^{p_i-1}) \quad \text{ for all } z \in \mathbb{R}^n\, .
\end{align}
In a first point, let us briefly show how to reduce the discussion of H\"older continuity to interior estimates.
\begin{rem}\label{extend}
Suppose that $\Omega$ satisfies the condition {\rm (C1)}. Then there is $R_0 > 0$ such that for arbitrary $x^0 \in \partial \Omega$, we find extensions (same notation) of $v^0$ and the data $a_i^0$, $b_i$, $B_i$ in the problem \eqref{0auxivhier} to $B_{R_0}(x^0) \setminus \Omega$, such that $v^0 \in (W^{1,p_i} \cap L^{q_i})(B_{R_0}(x^0))$ solves $-\divv (a^0_i(x, \, \nabla v^0 ) + b_i) + B_i = 0$ in $B_{R_0}(x^0)$ in the weak sense.
\end{rem}
We shall show the procedure for a sphere centred at the intersection $\partial S$ of the Dirichlet and the Neumann part, which is the most demanding case.
We let $x^0 \in \partial S$ be fixed, and $R_0 =  R_0(\partial S) >0$ be as required by the condition (C1). We assume without loss of generality that the boundary is already in the reference configuration addressed in (C1). This means that $x^0 = 0$, that the hypersurfaces $S_{1,R_0}(x^0) = S \cap  B_{R_0}(x^0) \subset \{x_n = 0\}$ and $S_{2,R_0}(x^0) = (\partial \Omega \setminus \overline{S}) \cap B_{R_0}(x^0) \subset \{x_{n-1} = 0\}$ are planar and orthogonal to each other, while $\Omega_{R_0}(x^0) = B_{R_0}(x^0) \cap \Omega$ is the intersection of $B_{R_0}(0)$ with the two half spaces $x_n <0$ and $x_{n-1} < 0$.
%the hypersurfaces $S_{1,R_0}(x^0) = S \cap  B_{R_0}(x^0)$ and $S_{2,R_0}(x^0) = (\partial \Omega \setminus \overline{S}) \cap B_{R_0}(x^0)$ are flat and perpendicular to each other at the junction $\overline{S_{1,R_0}(x^0)} \cap \overline{S_{2,R_0}(x^0)} \subseteq \partial S$. 
Otherwise, due to the assumptions (C1) on the domain $\Omega$, there is a bi-Lipschitzian transformation $T^{x^0}_S$ satisfying the properties (C1a) -- (C1d) and mapping into the reference configuration and back. We sketch in the section \ref{append} how this additional transformation can be handled.

The function $v^0 \in W^{1,p_i}(\Omega)$ satisfies
\begin{gather}\label{0auxiv}
	-\divv (a_i^0(x, \, \nabla v^0)+b_i(x))  + B_i(x) = 0 \text{ in } \Omega_{R_0} \,, \\  
	v = 0 \text{ on } S_{1,R_0} \, ,\qquad
	e^{n-1} \cdot (a_i^0(x, \, \nabla v) + b_i(x)) = 0 \text{ on } S_{2,R_0}\, ,\nonumber
\end{gather}
in the weak sense. 

We at first extend the data of the problem \eqref{0auxiv} for $0 \leq x_n < R_0$ and $x_{n-1} < 0$. For $x = (x^{\prime},x_n)$ we define $\calR (x) = (x^{\prime},-x_n)$. Then, for $x_n >0$, we set $v^0(x) := - v^0(\calR (x))$. This preserves the continuity across $x_n = 0$ since $v^0$ vanishes on this plane. For a vector field $A: \Omega\times \mathbb{R}^n \rightarrow \mathbb{R}^n$ we define
\begin{align}\label{calRmap}
	A_{j}(x, \, z^{\prime},z_n) := \begin{cases}
		- A_{j}(\calR(x),-z^{\prime},z_n) & \text{ for } j = 1,\ldots,n-1 \, ,\\
		A_{n}(\calR(x),-z^{\prime},z_n) & \text{ for } j = n \, .
	\end{cases} 
\end{align}
We apply this procedure to $A = a_i^0(x,z)$ and to $A = b_i(x)$. Note that for $z^{\prime} = 0$ the extended normal component $A_n$ is continuous across the plane $x_n = 0$\footnote{Hence the map $A_{n}(x,\nabla v^0(x))$ is formally continuous across the plane $x_n = 0$}. We extend the function $B_i$ by means of $B_i(x) = - B_i(\calR(x))$.

For $x_n > 0$ we can verify that $A(x,z) \cdot z = A(\calR(x), \, -z^{\prime}, \, z_n) \cdot(-z^{\prime},\, z_n)$, and therefore \eqref{a1weakenedhere}, \eqref{a1weakenedhereupper} show that the extended data satisfy
\begin{align}\label{a1extended}
	\frac{\nu_i}{2} \, |z|^{p_i} - C_0 \, \phi_0^{p_i}(\calR(x)) \leq a_i^0(x, \, z) \cdot z, \quad |a^0_i(x, \, z)| \leq \mu_i \, (|z|^{p_i-1} + \phi_0^{p_i-1}(\calR (x))) \, , 
\end{align}
for $x_n > 0$ and $z \in \mathbb{R}^n$ arbitrary with $\phi_0 = |\nabla \rho_i^S|$. Now we can verify in $B_{R_0}^-(0) := \{x \in B_{R_0}(0) \, : \, x_{n-1} < 0\}$ that the extended function $v^0$ belongs to $W^{1,p_i}(B_{R_0}^-)$ and is a weak solution to the Neumann problem
\begin{gather}\begin{split}\label{00auxiv}
		&		-\divv (a_i^0(x, \, \nabla v^0)+b_i(x))  + B_i = 0 \text{ in } B_{R_0}^-(0) \,, \\
		& e^{n-1} \cdot (a_i^0(x, \, \nabla v^0)+b_i(x)) = 0 \text{ on } B_{R_0}(0) \cap \{ x_{n-1} = 0\}\, .
	\end{split}
\end{gather}
For the second step we extend the data of the Homogeneous conormal problem \eqref{00auxiv} to $B_{R_0}^+(0)$. Now for $x_{n-1} >0$ we let $\calQ(x) = (x_1, \ldots,x_{n-2}, - x_{n-1}, x_n)$. Moreover we define $v(x) = v(\calQ(x))$, preserving the continuity across $x_{n-1} = 0$, and for vector fields on $\Omega \times \mathbb{R}^n$
\begin{align}\label{calQmap}
	A_{j}(x, \, z^{\prime},z_n) := \begin{cases}
		A_{j}(\calQ(x),z^{\prime},-z_n) & \text{ for } j \neq n-1 \, ,\\
		- A_{n-1}(\calQ(x),z^{\prime},-z_n) & \text{ for } j = n-1 \, ,
	\end{cases} 
\end{align}
which preserves the vanishing normal component on the plane $x_{n-1} = 0$ for $z = \nabla v^0$. We further let $B(x) = B(\calQ(x))$. For $x_{n-1} > 0$ we have this time $A(x,z) \cdot z = A(\calQ(x), \, z^{\prime}, \, -z_n) \cdot(z^{\prime},-z_n)$, and therefore \eqref{a1extended} shows that the data satisfy
\begin{align}\label{a1extended2}
	\frac{\nu_i}{2} \, |z|^{p_i} - \phi^{p_i}_0(\calR\circ\calQ(x)) \leq a_i^0(x, \, z) \cdot z, \quad |a_i^0(x,z)| \leq \mu_i \, |z|^{p_i-1} + \phi_0^{p_i-1}(\calR\circ\calQ (x)) \, 
\end{align}
where $x_{n-1} > 0$ and $z \in \mathbb{R}^n$ arbitrary. At last we verify that the extended function $v^0 \in W^{1,p_i}(B_{R_0})$ is a weak solution to $-\divv (a^0_i(x, \, \nabla v^0)+b_i)  + B_i = 0$ in $B_{R_0}(0)$.\\

The technique of Remark \ref{extend} shall be applied next to reduce the problem of global H\"older estimates to the case of interior estimates. The extension maps $\calR$ and $\calQ$ are reflections through planes, therefore the Lebesgue and Lebesgue-Morrey norms of the extended functions are controlled by the corresponding norms of the original data. This is left to the readers as an exercise.

\begin{lemma}\label{contihold}
We adopt the assumptions of Lemma \ref{ENERGYEST} and moreover {\rm (B3,$({\bf q},1)$)}, $\ldots$, {\rm (B3,$({\bf q},N)$)}, {\rm (C1)} and {\rm (C3$(1)$)}, $\ldots$, {\rm (C3$(N)$)}  Then there are $0 < \lambda_1, \ldots,\lambda_N < 1$ and $M_1$ such that every weak solution $w = (\rho,u) \in W^{1,\bf p}(\Omega;\, \mathcal{O}) \cap  L^{\bf q}(\Omega;\, \mathcal{O})$ to \eqref{systemcomp} with the energy dissipation property satisfies $\max_{i=1,\ldots,N} \|\rho_i\|_{C^{\lambda_i}(\overline{\Omega})} \leq M_1$. Moreover, there is $C >0 $ depending on the solution only via $\|w\|_{L^{\bf q}}$ such that 
\begin{align}
\label{migio} \int_{B_{R}(x^0) \cap \Omega} |\nabla \rho_i|^{p_i} \, dx \leq C \, R^{n-p_i+\lambda_i \, p_i} \quad \text{ for all } \quad R>0 \text{ and } x^0 \in \overline{\Omega}\, .
\end{align}
\end{lemma}
\begin{proof}
The H\"older bound shall follow from the inequality \eqref{migio}. Therein we can restrict to considering four different types of balls:
\begin{enumerate}
	\item $x^0 \in \Omega$ such that $B_{R_0}(x^0) \subset\!\!\subset \Omega$,
	\item $x^0 \in \partial \Omega$ and $B_{R_0}(x^0) \cap \partial \Omega \subset\!\!\subset S$,
	\item $x^0 \in \partial \Omega$ and $B_{R_0}(x^0) \cap \partial \Omega \subset\!\!\subset \partial\Omega\setminus\overline{S}$,
	\item $x^0\in \partial S$. 
\end{enumerate}
We shall concentrate on the case $x^0 \in \partial S$ since the other cases can be treated with similar or simpler methods. We recall that $v^0(x) = \rho_i(x)- \rho_i^S(x) \in W^{1,p_i}(\Omega)$, which solves \eqref{0auxivhier}.
%\begin{gather}\label{0auxivhier}
%	-\divv (a^0_i(x, \, \nabla v^0 ) + b_i) + B_i = 0 \text{ in } \Omega \,, \\  
%	v^0 = 0 \text{ on } S \, ,\qquad
%	\nu(x) \cdot (a^0_i(x, \, \nabla v^0) + b_i) = 0 \text{ on } \partial \Omega \setminus S\, ,\nonumber
%\end{gather}
%with the new vector field $a_i^0(x,z) := \widetilde{a}_i(x, \, z + \nabla \rho_i^S(x))$. (A1) implies that
%\begin{align*}
%	\nu_i \, |z + \nabla \rho_i^S(x)|^{p_i} \leq a_i^0(x, \, z) \cdot (z + \nabla \rho_i^S(x)) \leq 	 a_i^0(x, \, z) \cdot z + \mu_i \, |z + \nabla \rho_i^S(x)|^{p_i-1} \, |\nabla \rho_i^S(x)| \, .
%\end{align*}	
%Hence with the help of elementary inequalities we find a constant $C_0 = C_0(\nu_i,\mu_i,p_i) >0$ such that
%\begin{align}\label{a1weakenedhere}
%	& \frac{\nu_i}{2} \, |z|^{p_i} - C_0 \, |\nabla \rho_i^S(x)|^{p_i} \leq a_i^0(x, \, z) \cdot z\, , \\
%\label{a1weakenedhereupper}	& |a_i^0(x, \, z)|	\leq \mu_i \, (|z|^{p_i-1} +|\nabla \rho_i^S(x)|^{p_i-1}) \quad \text{ for all } z \in \mathbb{R}^n\, .
%\end{align}
Applying the extension procedure \eqref{0auxiv}--\eqref{a1extended2} we see that $v^0$ solves 
\begin{align}\label{0auxivextend}
-\divv (a^0_i(x, \, \nabla v^0 ) + b_i) + B_i = 0 \quad\text{ in } \quad B_{R_0}(x^0) \, ,
\end{align}
with extended data in $B_{R_0}(x^0)$ satisfying bounds in the same norms as the original ones, and in particular \eqref{a1extended2},
%\begin{align}\label{a1weakenedhere2}
%\nu \, |z|^{p_i} - C \, \phi_0^{p_i}(x) \leq a_i^0(x, \, z) \cdot z, \quad
%|a_i^0(x, \, z)|	\leq \mu \, (|z|^{p_i-1} + \phi_0^{p_i-1}(x)) \, \, 
%\end{align}
%for all $z \in \mathbb{R}^n$, with positive constants $\nu$, $\mu$ ans $C$ depending only on $\nu_i$, $\mu_i$ and $p_i$ from (A1) and 
with $\phi_0 \in L^{p_i^S}(B_{R_0}(x^0))$ satisfying the bound $$\|\phi_0\|_{L^{p_i^S}(B_{R_0}(x^0))} \leq C\,  \|\nabla \rho_i^S\|_{L^{p_i^S}(\Omega_{R_0}(x^0))}$$ with $p_i^S > n$ from (C3$(i)$). Next we want to show that there is $\lambda_i >0$ such that for all $B_r(y) \subset B_{R_0/4}(x^0)$ and all $0< \rho \leq r$ 
\begin{align}\label{claim}
	\int_{B_{\rho}(y)} |\nabla v^0|^{p_i}\, dx \leq c_1\, \Big(\frac{\rho}{r}\Big)^{n-p_i+\lambda_i \, p_i}  \, \int_{B_r(y)} |\nabla v^0|^{p_i} \, dx + c_2 \, r^{n-p_i+\lambda_i \, p_i} \, , 
\end{align}
with constants $c_1$ and $c_2$ depending only on the data. To this aim, for $r >0$ and $y \in B_{R_0}(x^0)$ such that $B_r(y) \subset B_{R_0/4}(x^0)$ we first consider the weak solution $w\in W^{1,p_i}(B_{r}(y))$ to
\begin{align}\label{auxiv}
-\divv a_i^0(x, \, \nabla w) = 0 \text{ in } B_{r}(y), \quad w = v^0 \text{ on } \partial B_r(y) \, .
\end{align}
%Obviously, $\|v\|_{L^{\infty}(\Omega_R)} \leq \|\rho_i\|_{L^{\infty}}$ showing in connection to \eqref{rhoibounded} that $v$ is bounded by the data.
%We insert $w - v^0$ into the weak formulation of the latter problem and by means of \eqref{a1weakenedhere2} and H\"older's inequality, we at first obtain the estimate
%\begin{align}\label{vprops1}
%\int_{B_r(y)} |\nabla w|^{p_i} \, dx \leq & C \, \int_{B_r(y)} (|\nabla v^0|^{p_i} + \phi_0(x))\, dx \, .
%\end{align}
Due to the Lemma auxiliary Lemma \ref{Holder} after this proof, $w$ is H\"older continuous in $B_{r/2}(y)$ with exponent $\kappa_i > 0$, and the norm $\|w\|_{C^{\kappa_i}(\overline{B_{r/2}(y)})}$ is controlled by the data. We also find constants $c_k^{\prime} = c_k^{\prime}(\Omega,p_i,p_i^S,\|\rho_i^S\|_{W^{1,p_i^S}},\mu_i,\nu_i)$ such that
\begin{align}
 \label{vprops2} \int_{ B_{\rho}(y)} |\nabla w|^{p_i} \, dx \leq c_1^\prime \, \Big(\frac{\rho}{r}\Big)^{n-p_i + \kappa_i \, p_i}  \,  \int_{ B_{r}(y)} |\nabla w|^{p_i} \, dx + c_2^\prime \, \rho^{n-p_i + \kappa_i \, p_i} \, .
 %
 %
 %\int_{B_{\rho}(y)} |\nabla w|^{p_i} \, dx \leq &  
 %c_1\, \Big(\frac{\rho}{r}\Big)^{n-p_i+\kappa_i \, p_i}  \, \int_{B_r(y)} |\nabla w|^{p_i} \, dx + c_2 \, \rho^{n-p_i+\kappa_i \, p_i}\, \text{ for } 0 < \rho \leq r \, .
 %c \, \Big(\frac{\rho}{r}\Big)^{n-p_i+\kappa_i \, p_i}  \, .
\end{align}
Applying \eqref{a1extended} again, and using the monotonicity property of $z\mapsto a_i^0(x,z)$, we have
\begin{align*}
& \int_{B_{\rho}(y)} \nu \,  |\nabla v^0|^{p_i} - C\, \phi^{p_i}_0(x) \, dx \leq \int_{B_{\rho}(y)} a_i^0(x, \, \nabla v^0) \cdot \nabla v^0 \, dx \, \\
&  \quad  = \int_{B_{\rho}(y)} (a_i^0(x, \, \nabla v^0) -a_i^0(x,\nabla w) )\cdot \nabla (v^0 - w) \,  dx\\
& \quad \phantom{=} + \int_{B_{\rho}(y)}  a_i^0(x, \, \nabla v^0) \cdot \nabla w +  a_i^0(x, \, \nabla w) \cdot \nabla (v^0 - w)\, dx \, \\
& \quad  \leq \int_{B_{r}(y)} (a_i^0(x, \, \nabla v^0) - a_i^0(x,\nabla w) )\cdot \nabla (v^0 - w)\,  dx \\
& \quad \phantom{\leq} + \int_{B_{\rho}(y)}  a_i^0(x,\,\nabla v^0) \cdot \nabla w + a_i^0(x, \, \nabla w) \cdot \nabla v^0   - \nu\, |\nabla w|^{p_i} + C \, \phi_0^{p_i}(x) \, dx \, .
\end{align*}
Employing the definition of $w$, we find that $\int_{B_{r}(y)}a_i^0(x, \, \nabla w) \cdot \nabla (v^0 - w) \,  dx = 0$. Due to $v^0$ being a weak solution to \eqref{0auxivextend} in $B_{R_0}(x^0)$, we moreover get
\begin{align*}
 \int_{B_{r}(y)} a_i^0(x,\, \nabla v^0) \cdot \nabla (v^0 - w) \,  dx = \int_{B_{R_0}(x^0)} (A^0_i-b_i) \cdot \nabla ( \widehat{v^0-w}) \,  dx\\
 = - \int_{B_{R_0}(x^0)} b_i\cdot \nabla \, (\widehat{v^0-w}) + B_i \,( \widehat{v^0-w} )\, dx   = - \int_{B_r(y)} b_i\cdot \nabla \, (v^0-w)  + B_i \, (v^0-w)\, dx \, .
\end{align*}
Here $\widehat{v^0-w}$ is the extension by zero of $v^0-w$ outside of $B_r(y)$. Hence, it follows that
\begin{align}\label{premiere}
& \nu \, \int_{B_{\rho}(y)} |\nabla v^0|^{p_i} + |\nabla w|^{p_i} \, dx \leq - \int_{B_{r}(y)} b_i \cdot \nabla (v^0-w) + B_i \, (v^0-w)\, dx\nonumber\\
& \quad + \int_{B_{\rho}(y)} |a^0_i(x, \, \nabla  v^0)| \, |\nabla w| + |a_i^0(x, \, \nabla w)| \, |\nabla v^0| + C \, \phi^{p_i}_0(x) \, dx \, .
\end{align}
Invoking \eqref{a1extended2}, Young's generalised inequality, and \eqref{un}, we can estimate
\begin{align}\label{deuxi}
& \int_{B_{\rho}(y)} |a_i^0(x, \, \nabla  v^0)| \, |\nabla w| + |a_i^0(x, \, \nabla w)| \, |\nabla v^0| \, dx \nonumber\\
& \quad \leq \mu \, \int_{B_{\rho}(y)} (|\nabla v^0|^{p_i-1}+\phi_0^{p_i-1}) \, |\nabla w| + (|\nabla w|^{p_i-1}+\phi_0^{p_i-1}) \, |\nabla v^0|\, dx \nonumber\\
& \quad \leq \frac{\nu}{2} \, \int_{B_{\rho}(y)} |\nabla v^0|^{p_i} \, dx + c(p_i,\mu,\nu) \, \int_{B_{\rho}(y)} |\nabla w|^{p_i} + \phi_0^{p_i}(x) \, dx \nonumber\\
& \quad \leq \frac{\nu}{2} \, \int_{B_{\rho}(y)} |\nabla v^0|^{p_i} \, dx +c(p_i,\mu,\nu)\, \|\phi_0\|_{L^{p_i^S}}^{p_i} \, \omega_1^{1-p_i/p_i^S} \, \rho^{n \, (1-p_i/p_i^S)} \, \nonumber\\
& \quad \qquad + c_1 \, \Big(\frac{\rho}{r}\Big)^{n-p_i + \kappa_i \, p_i}  \,  \int_{ B_{r}(y)} |\nabla w|^{p_i} \, dx + c_2 \, \rho^{n-p_i + \kappa_i \, p_i}
\end{align}
Moreover, by the generalised Young inequality and the assumption (B3), we can show for $\epsilon > 0$ arbitrary that 
\begin{align}\label{troisi}
\Big| \int_{B_{r}(y)} b_i \cdot \nabla (v^0-w) \, dx \Big| \leq & \epsilon \, \int_{B_r(y)}   |\nabla  v^0|^{p_i} + |\nabla w|^{p_i} \, dx + c_{\epsilon} \, \int_{B_{r}(y)} |b_i|^{p^{\prime}_i} \, dx \, \nonumber \\
\leq & \epsilon \, \int_{B_r(y)}   |\nabla v^0|^{p_i} + |\nabla w|^{p_i} \, dx + c_{\epsilon} \, [b_i]_{p_i^{\prime},\alpha_i^b}^{\rm M} \,  r^{\alpha_i^b}\, .
% \nonumber\\
%\leq & \epsilon  \, \int_{B_r(y)}   |\nabla v^0|^{p_i} + |\nabla w|^{p_i}  \, dx + c_{\epsilon}  \,  r^{\alpha_i^B} 
\end{align}
%where we used \eqref{vprops1}.
Invoking in addition to these tools the Poincar\'e and Sobolev inequalities, we also obtain that
\begin{align*}%\label{quatri}
	\Big| \int_{B_{r}(y)} B_i \, (v^0-w) \, dx \Big| \leq & \|v^0-w\|_{L^{p_i\kappa_i^B}(B_r; \, d\calB)} \, [\calB(B_r)]^{1-\frac{1}{\kappa^B_i}}  \nonumber\\
	\leq & C \, \|B\|_{L^{1,\alpha_i^B}_M(\Omega)} \, \|\nabla(v^0-w)\|_{L^{p_i}(B_r)}\, [\calB(B_r)]^{1-\frac{1}{\kappa^B_i}}  \nonumber \\
	\leq & \epsilon \,  \,\int_{B_r(y)}   |\nabla v^0|^{p_i}+ |\nabla w|^{p_i}  \, dx + c_{\epsilon} \,\|B\|_{L^{1,\alpha_i^B}_M(\Omega)}^{p_i^\prime}  \, r^{p_i^\prime\, \alpha_i^B \,(1-\frac{1}{\kappa^B_i})} \, ,
\end{align*}
with a Poincar\'e--Sobolev constant $C$ of Lemma \ref{lieberman}. As $\kappa_i^B = 1 + \theta_i^B \, p_i/(n-p_i)$ and $\alpha_i^B = n - p_i + \theta_i^B \, p_i$, we can verify that 
\begin{align}\label{quatri}
	\Big| \int_{B_{r}(y)} B_i \, (v^0-w) \, dx \Big| \leq 
	 \epsilon \,  \,\int_{B_r(y)}   |\nabla v^0|^{p_i}+ |\nabla w|^{p_i}  \, dx + c_{\epsilon} \,\|B\|_{L^{1,\alpha_i^B}_M(\Omega)}^{p_i^\prime}  \, r^{n - p_i + \theta_i^B \, p_i^{\prime}\, p_i} \, ,
\end{align}
%\begin{align}\label{quatri}
%\Big| \int_{B_{r}(y)} B_i \, (v^0-w) \, dx \Big| \leq & \|v^0-w\|_{L^{p_i^*}(B_r)} \, \|B_i\|_{L^{(p_i^*)^{\prime}}(B_r)}  \nonumber\\
% \leq & c_{PS} \, \|\nabla(v^0-w)\|_{L^{p_i}(B_r)}\, \|B_i\|_{L^{(p_i^*)^{\prime}}(B_r)} \nonumber \\
%\leq & \epsilon \,  \,\int_{B_r(y)}   |\nabla v^0|^{p_i}+ |\nabla w|^{p_i}  \, dx + c_{\epsilon} \, \big([B_i]_{(p^*_i)^{\prime},\alpha_i^B}\big)^{\frac{p_i^{\prime}}{(p_i^*)^{\prime}}} \, r^{\alpha_i^B \, \frac{p_i^{\prime}}{(p_i^*)^{\prime}}} \, ,
%\end{align}
%with a Poincar\'e--Sobolev constant $C_{PS}$. 
Due to (B3) and the Lemma \ref{ENERGYEST}, the quantities $[b_i]_{p_i^{\prime},\alpha_i^b}$ and $[B_i]_{1,\alpha_i^B}$ are bounded independently of the solution.
With $\lambda_i := \min\{\kappa_i, \, \theta_i^b, \, \theta_i^B\}$, collection of the inequalities \eqref{premiere}, \eqref{deuxi}, \eqref{troisi} and \eqref{quatri} then yields
\begin{align*}
\int_{B_{\rho}(y)} |\nabla v^0|^{p_i}+ |\nabla w|^{p_i} \, dx \leq C\Big(\epsilon+\Big(\frac{\rho}{r}\Big)^{n-p_i + \kappa_i \, p_i}\Big)\,  \, \int_{B_r(y)} |\nabla v^0|^{p_i}+ |\nabla w|^{p_i}  \, dx + C_{\epsilon} \, r^{n-p_i+\lambda_i \, p_i} \, . 
\end{align*}
This inequality is valid for all $B_r(y) \subset B_{R_0/4}(x^0)$ and $\rho \leq r$. We finish the proof of \eqref{claim} using the Lemma 2.1 of Chapter III in \cite{MR0717034}. 

Similar arguments apply to balls near the Dirichlet or Neumann boundary and in the interior. The claim follows.
\end{proof}

\begin{lemma}\label{Holder}
	There is $0 < \kappa_i \leq 1-n/p_i^S$ and constants $c_k^\prime$, $k = 1,2$ depending only on $\Omega$, $p_i$, $\|\rho^S_i\|_{W^{1,p_i^S}}$ and the constants $\nu$ and $\mu$ in \eqref{a1extended2} such that the weak solution $w\in W^{1,p_i}(B_{r}(y))$ to \eqref{auxiv} satisfies, for all $0< \rho \leq r$, 
	\begin{align}\label{un}
%		\sup_{x,\tilde{x} \in B_{r/2}(y)} |w(x)-w(\tilde{x})| \leq c_0 \, \Big(\frac{|x-\tilde{x}|}{r}\Big)^{\kappa_i} \quad  \text{ for all } \quad 0< \rho \leq r \, , \\
		 \int_{ B_{\rho}(y)} |\nabla w|^{p_i} \, dx \leq c_1^\prime \, \Big(\frac{\rho}{r}\Big)^{n-p_i + \kappa_i \, p_i}  \,  \int_{ B_{r}(y)} |\nabla w|^{p_i} \, dx + c_2^\prime \, \rho^{n-p_i + \kappa_i \, p_i}  \,.
		\end{align}
	\end{lemma}
\begin{proof}
	The solution $w$ to \eqref{auxiv} obviously satisfies $\sup_{x \in B_r(y)} w \leq \sup_{\partial B_r(y)} v^0$ and $\inf_{B_r(y)} w \geq \inf_{\partial B_r(y)} v^0$. Thus, in view of Lemma \ref{linfty}, $w$ is bounded by the data. Moreover we can insert $w - v^0$ into the weak formulation of the latter problem and by means of \eqref{a1extended} and H\"older's inequality, we at first obtain the estimate
	\begin{align}\label{vprops1}
	\int_{B_r(y)} |\nabla w|^{p_i} \, dx \leq & C \, \int_{B_r(y)} (|\nabla v^0|^{p_i} + \phi^{p_i}_0(x))\, dx \, .
	\end{align}
	In order to obtain \eqref{un} with constants that do not depend on the radius $r$, we apply the transformation $\xi := (x-y)/r$. For $\xi \in B_1(0)$ and $z\in \mathbb{R}^n$, we introduce $F_i(\xi, \, z) := r^{p_i-1} \,  a^0_i(r \, \xi, \, r^{-1} \, z)$. We observe that
\begin{align}\label{mitch} 
 F(\xi, \, z) \cdot z = r^{p_i} \, a_i^0\big(r\, \xi, \, r^{-1} \, z\big) \cdot (r^{-1} \, z) \geq \nu\, |z|^{p_i} - C\, (\phi_0^1(\xi))^{p_i} \, ,
\end{align}
	and similarly $ |F(\xi, \, z)| \leq  \mu \, (|z|^{p_i-1} + (\phi_0^1(\xi))^{p_i-1})$ with $\phi_0^1(\xi) := r \, \phi_0(r\, \xi)$.
	%%Moreover, with $p_i^S > n$ from the condition ??? we have
	%%\begin{align*}
	%%	\|F_i\|_{L^{p_i^S}(B_1^-(0))} = 2^{\frac{1}{p_i^S}} \, R^{1-\frac{n}{p_i^S}} \, \|\nabla \rho_i^{S}\|_{L^{p_i^S}(\Omega_R(0))} \, .
	%%\end{align*}
	Then we let $w^1(\xi) := w(r \, \xi)$ and see that \eqref{auxiv}$_1$ is equivalent to $\divv_{\xi} F(\xi, \, \nabla_\xi w^1) = 0$ in $B_1(0)$. 
	Using standard results for the quasilinear elliptic problems, we obtain the H\"older continuity of $w^1$ with exponent $0<\kappa_i \leq 1-n/p_i^S$ in the set $\overline{B_{1/2}(0)}$. For instance we can follow the approach in \cite{MR0244627}, Chapter 2, paragraph 7, showing that $w^1$ belongs to a De Giorgi class described there on page 90, see the next remark \ref{EnnodeGiorgi}.
	%
	%%
	%%described briefly hereafter for the case of mixed boundary conditions.
	%%
	%%Considering arbitrary spheres $B_{\rho}(y^0)$ with centre $y^0 \in \overline{B_{1}^-(0)}$ and $\rho >0$, we denote $K_{\rho}(y^0) = B_1^-(0) \cap B_{\rho}(y^0)$. Consider arbitrary nonnegative, smooth $\zeta$ with support in $B_{\rho}(y^0)$ and real number $k$. Let $\calw(y) := v_1(y) - \rho_i^{\hat{S}_1}(y)$, $\calw_k := \max\{\calw-k,0\}$ for $k >0$ and $\calw_k = \max\{\calw-k,0\}+k$ if $k < 0$. Note that $\calw_k$ vanishes on $\hat{S}_1$. The weak formulation of the problem is multiplied with the testfunction $\calw_k \, \zeta^{p_i}$ to obtain with $A_{\rho,k} :=\{y \in K_{\rho} \, :\, \calw(y) >k  \}$ the relation
	%%\begin{align*}
	%%	\int_{A_{\rho,k}} \gamma(y,\, \nabla v_1) \cdot (\nabla(v_1-\rho_i^{\hat{S}_1}) \, \zeta^{p_i} + p_i\, \zeta^{p_i-1} \, \gamma(y,\, \nabla v_1)\nabla\zeta \, \calw_k) \, dy  = 0 \, .
	%%\end{align*}
	%%By means of \eqref{mitch} and elementary inequalities, we easily show that
	%%\begin{align*}
	%%	\int_{A_{\rho,k}} \zeta^{p_i} \, |\nabla \calw|^{p_i} \, dy \leq C \, \int_{ A_{\rho,k}} \zeta^{p_i} \, |\nabla \rho_i^{\hat{S}_1}|^{p_i}  + |w_k|^{p_i} \, |\nabla \zeta|^{p_i}\, dy  \, .
	%%\end{align*}
	%
%	Then, for $0< \rho < r/2$ and $x, \tilde{x} \in B_{\rho}(y)$ arbitrary, we notice that
%	\begin{align*}
	% |w(x) - w(\tilde{x})| = \Big| w^1\Big( \frac{x}{r}\Big) - w^1\Big( \frac{\tilde{x}}{r}\Big)\Big| \leq [\tilde{v}]_{C^{\kappa_i}(\overline{B_{1/2}(0)})} \, \Big(\frac{x-\tilde{x}}{r}\Big)^{\kappa_i} \, .
	%\end{align*}
	%This proves the inequality \eqref{un}$_1$.
	%
	In order to obtain the inequality \eqref{un}, we consider first $\rho \leq r/4$ and we multiply the equation \eqref{auxiv} with $(w - w_{\rho}) \, \zeta^{p_i}$. Here $w_{\rho} := \fint_{B_{2\rho}(y)} v(x) \, dx$ is the average and $\zeta \in C^1_0(B_{2\rho})$ is an appropriate cutoff function, identically equal to one on $B_\rho$, such that $|\nabla \zeta| \leq k_1 \, \rho^{-1}$. We easily obtain that
	\begin{align*}
	 \int_{B_r(y)} |\nabla w|^{p_i} \, \zeta^{p_i} \, dx \leq C \, \int_{B_r(y)} \, \phi^{p_i}_0(x) \, \zeta^{p_i} + |w-w_{\rho}|^{p_i} \, |\nabla \zeta|^{p_i} \, dx \, .
	\end{align*}
	Using the H\"older continuity of $w$ on $B_{r/2}$, that $2\rho < r/2$ and that $|\nabla \zeta| \leq k_1 \, \rho^{-1}$, we get
	\begin{align*}
	 \int_{B_r(y)} |\nabla w|^{p_i} \, \zeta^{p_i} \, dx \leq \|\phi_0\|_{L^{p_i^S}(B_{2\rho})}^{p_i} \, [\lambda_n(B_{\rho})]^{1-\frac{p_i}{p_i^S}} + [w]_{C^{\kappa_i}(B_{2\rho})}^{p_i}\, \lambda_n(B_{2\rho}) \, \rho^{(\kappa_i-1) \, p_i}\\
	 \leq \|\phi_0\|_{L^{p_i^S}(B_{2\rho})}^{p_i} \, \omega_1^{1-p_i/p_i^S} \, \rho^{n-p_i+p_i(1-\frac{n}{p_i^S})} + [w]_{C^{\kappa_i}(B_{2\rho})}^{p_i}\, \omega_1 \, 2^n \, \rho^{n-p_i + \kappa_i \, p_i} \, .
	\end{align*}
Since $\kappa_i \leq 1-n/p_i^S$ it follows that 
	\begin{align}\label{rhopetit}
		\int_{B_{\rho}(y)} |\nabla w|^{p_i} \, dx \leq c \, \rho^{n-p_i + \kappa_i \, p_i} \, ,
	\end{align}
	where we assumed that $\rho \leq r/4$. For $r/4 \leq \rho \leq r$ we have
	\begin{align}\label{rhogrand}
		\int_{B_{\rho}(y)} |\nabla w|^{p_i} \, dx \leq \int_{B_{r}(y)} |\nabla w|^{p_i} \, dx \leq 4^{n} \, \left(\frac{\rho}{r}\right)^{n-p_i+\kappa_i \, p_i} \, \int_{B_{r}(y)} |\nabla w|^{p_i} \, dx\, .
	\end{align}
Summing up \eqref{rhopetit} and \eqref{rhogrand}, the claim follows.
%	where we invoked \eqref{vprops1}.
	%If $B_{2r}(x^0)$ intersects $\partial \Omega$, we multiply \eqref{auxiv} with $(v-\rho_i^S) \, \zeta^{p_i}$. In this case we obtain that
	%\begin{align*}
	%	\int_{\Omega} |\nabla v|^{p_i} \, \zeta^{p_i} \, dx \leq \int_{\Omega} |\nabla \rho_i^S|^{p_i} \, \zeta^{p_i} \, dx +  \Big(\frac{p_i\mu_i}{\nu_i}\Big)^{p_i} \, \int_{\Omega} |v-\rho_i^S|^{p_i} \, |\nabla \zeta|^{p_i} \, dx \, .
	%\end{align*}
	%For some $x^1 \in B_{2r}(x^0) \cap \partial \Omega$, we have $|v(x) - \rho_i^S(x)| = |v(x) - \rho_i^S(x) - (|v(x^1) - \rho_i^S(x^1)| )|  \leq c \, (|x-x^1|/R)^{\kappa_i}$ and can argue just the same. 
	%
	%Thus the statement is proved for $r \leq R/2$. However, for $R \geq r \geq R/4$, the inequality is obvious with $c = 4^{n-p_i}$	
\end{proof}
\begin{rem}\label{EnnodeGiorgi}
	The auxiliary function $w^1$ in the proof of Lemma \ref{Holder} belongs to the De Giorgi class
$\boldsymbol{\calB}_{p_i}(B_1(0), \,  \|\rho_i\|_{L^{\infty}}, \, \gamma, \, +\infty, 1/p_i^S)$ for a $\gamma$ depending only on $\Omega$, $p_i$, $p_i^S$ and $\|\phi_0\|_{L^{p_i^S}(\Omega)}$.
\end{rem}
\begin{proof}
Considering arbitrary spheres $B_{s}(\xi^0)$ with centre $\xi^0 \in \overline{B_{1/2}(0)}$ and $0 < s \leq 1/2$. Consider arbitrary nonnegative, smooth $\zeta$ with support in $B_{s}(\xi^0)$ and real number $k$. Let $w_k := \max\{w^1-k,0\}$. The weak formulation of the problem is multiplied with the testfunction $w_k \, \zeta^{p_i}$ to obtain with $A_{s,k} :=\{\xi \in B_{s}(\xi^0) \, :\, w^1(\xi) >k  \}$ the relation
	\begin{align*}
		\int_{A_{s,k}} F(\xi,\, \nabla_\xi w^1) \cdot \nabla_\xi w^1 + p_i\, \zeta^{p_i-1} \, F(\xi,\, \nabla_\xi w^1)\cdot\nabla_\xi \zeta \, w_k \, d\xi  = 0 \, .
	\end{align*}
	By means of  elementary inequalities, we show that
	\begin{align*}
		\int_{A_{s,k}} \zeta^{p_i} \, |\nabla w^1|^{p_i} \, d\xi \leq C \, \int_{ A_{s,k}} \zeta^{p_i} \, (\phi_0^1)^{p_i}  + |w_k|^{p_i} \, |\nabla \zeta|^{p_i}\, d\xi  \, .
	\end{align*}
By means of H\"older's inequality
\begin{align*}
	\int_{ A_{s,k}} \zeta^{p_i} \, (\phi_0^1)^{p_i} \, d\xi \leq \|\phi_0^1\|_{L^{p_i^S}(B_1(0))}^{p_i} \, [\lambda_n(A_{s,k})]^{1-\frac{p_i}{p_i^S}}   = r^{1-\frac{n}{p_i^S}} \,  \|\phi_0\|_{L^{p_i^S}(B_{\rho}(y))}^{p_i} \, [\lambda_n(A_{s,k})]^{1-\frac{p_i}{p_i^S}}  \, ,
\end{align*}
where we applied the transformation formula and the definition of $\phi_0^1(\xi) = r \, \phi_0(r(y+\xi)$) to transform the norm.
Hence we reach
	\begin{align*}
	\int_{A_{s,k}} \zeta^{p_i} \, |\nabla w^1|^{p_i} \, d\xi \leq C \, (1+ r^{1-\frac{n}{p_i^S}} \, \|\phi_0\|_{L^{p_i^S}(B_{r}(y))} )^{p_i} \, \Big([\lambda_n(A_{s,k})]^{1-\frac{p_i}{p_i^S}} + \int_{ A_{s,k}}  |w_k|^{p_i} \, |\nabla \zeta|^{p_i}\, d\xi \Big) \, 
\end{align*}
which is (6.2) on page 82 of \cite{MR0244627}.
\end{proof}
\subsection{H\"older regularity in the energy equation}
Now the higher regularity of $u$ follows from known results of potential theory and from the following Lemma.
\begin{lemma}\label{measure}
Adopting the same assumptions as in Lemma \ref{contihold}, we suppose that the constants $p_1,\ldots,p_{N+1}$ moreover satisfy
\begin{align*}
 p_{N+1} > \max_{i=1,\ldots,N} (1-\lambda_i) \, p_i \, .
\end{align*}
We let $ \theta := 1- \frac{\max_{i=1,\ldots,N} (1-\lambda_i) \, p_i}{p_{N+1}}$. Then there $c>0$ depending only on the data and on the solution via $\| w\|_{L^{\bf q}}$ such that for all weak solutions $w = (\rho, \, u)$ with energy dissipation property for \eqref{systemcomp} we have
\begin{align*}
\int_{\Omega_R(x^0)} \left|\Pi(x,w, \, \nabla w) \right|  \, dx\leq c \, R^{n - p_{N+1} + \theta \, p_{N+1}} \quad \text{ for all } x^0 \in \overline{\Omega}, \,  0 < R \leq {\rm diam}(\Omega)\, .
\end{align*} 
\end{lemma}
\begin{proof}
By a direct application of (A1) and (B0) (cp.\ \eqref{CRHS}$_1$) and of Lemma \ref{contihold}, we get
\begin{align*}
 \int_{\Omega_{R}(x)} |\Pi| \, dx \leq &  L \, \sum_{i=1}^N \int_{\Omega_R(x^0)} |\nabla \rho_i|^{p_i}  \, dx\leq C \, R^{n-p_i +\lambda_i \, p_i} \, .
 %\\
%\int_{\Omega_R} |b_i \cdot \nabla \rho_i| \, dx \leq & \frac{1}{p_i^{\prime}} \, \int_{\Omega_{R}} |b_i|^{p_i^{\prime}} \, dx + \frac{1}{p_i} \, \int_{\Omega_R} |\nabla \rho_i|^{p_i} \, dx \leq C \, R^{n-p_i +\lambda_i \, p_i} \, .
\end{align*}
By the choice of $\theta$ we have $\min_{i=1,\ldots,N} n - p_i + \lambda_i \, p_i = n - p_{N+1} + \theta \, p_{N+1}$, and the claim follows.
\end{proof}
\begin{rem}
 Thanks to the property shown in Lemma \ref{measure} the measure
 \begin{align*}
  \calm(\phi) := \int_{\Omega} \Pi(x, \, w, \, \nabla w) \, \phi \, dx 
 \end{align*}
extends to an element of $[W^{1,p_{N+1}}_S(\Omega)]^*$ with the norm controlled independently on $(\rho, \, u)$, see a.o.\ Lemma 1.1 in \cite{MR1233190}.
\end{rem}

\begin{coro}\label{coroN+1}
 Under the assumptions of Lemma \ref{measure} and {\rm (B2$({\bf q})$)}, there is a constant $M_2$ such that $\|u\|_{W^{1,p_{N+1}}(\Omega)} \leq M_2$ for all weak solutions $(\rho,u) \in W^{1,\bf p} \cap L^{\bf q}$ with the energy dissipation property. If also {\rm (B3$({\bf q},N+1)$)} and {\rm (C2$(N+1)$)} are valid, then $\|u\|_{L^{\infty}(\Omega)} \leq M_2^{\prime}$.
\end{coro}
\begin{proof}
For ease of writing, we drop the $N+1$ indices throughout the proof, hence writing $p =p_{N+1}$, $\nu = \nu_{N+1}$, etc. Invoking the energy dissipation property of Definition \ref{solandenergy}, we choose $g(u) = \sign(u) \, \min\{|u|, \, L\} = T_{L}(u)$ in \eqref{heatdistrib}. 
%Here we let $L > \|u^S\|_{L^{\infty}(S)}$. 
Thanks to the growth condition (A1), we obtain that
\begin{align*}
& \nu \, \int_{\Omega} T^{\prime}_L(u-u^S) \, |\nabla u|^{p} \, dx  \leq \int_{\Omega} [A-b](x, \, w, \, \nabla w) \cdot (\nabla T_L(u-u^S) + T^{\prime}_L(u-u^S) \, \nabla u^S) \, dx \\
& \quad \leq \int_{\Omega} \Big(\Pi - B\Big) \, T_{L}(u-u^S) + \, T_{L}^{\prime}(u-u^S)  \, (A\cdot \nabla u^S - b\cdot \nabla u )\, dx \, .
\end{align*}
Standard inequalities yield
\begin{align*}
 \frac{\nu}{2} \, \int_{\Omega} T^{\prime}_L(u-u^S) \, |\nabla u|^{p} \, dx \leq & C_{PS} \, (\|\calm\|_{[W^{1,p}_0]^*} + \|B\|_{L^{(p^*)^{\prime}}}) \, \|\nabla T_{L}(u-u_S)\|_{L^{p}(\Omega)}\\
 & + \int_{\Omega} c(p,\nu^{-1},\mu) \,  (|b|^{p^{\prime}} +|\nabla u^S|^p) \, dx \, .
\end{align*}
%which also yields
%\begin{align*}
%	\frac{\nu}{4} \, \int_{\Omega} T^{\prime}_L(u-u^S) \, |\nabla u|^{p} \, dx \leq & C \, (\|\calm\|_{[W^{1,p}_0]^*}^{p^{\prime}} + \|B\|_{L^{(p^*)^{\prime}}}^{p^{\prime}} + \|b\|_{L^{p^{\prime}}}^{p^{\prime}} +
%	\|\nabla u_S\|_{L^{p}(\Omega)}^p) \, ,
%\end{align*}
Now we can let $L \rightarrow \infty$ proving that $u \in W^{1,p}(\Omega)$, with the norm estimated by the quantities on the right-hand side which, in view of (B2), are bounded independently of the solution. 

Next, to prove an upper bound for $u$, we employ in the weak form \eqref{transportdistrib} the testfunction $u_L = \max\{0,u-L\}$ with $L > \sup_S u^S$. For $\Omega_{L} := \{x \, : \, u(x) > L\}$ with characteristic function denoted by $\chi_{\Omega_L}$, we obtain that
\begin{align}\label{laberlingu}
  \nu \, \int_{\Omega_{L}} |\nabla u|^{p} \, dx  \leq & \int_{\Omega_L} (A- b)(x, \, w, \, \nabla w) \cdot \nabla u \, dx \nonumber\\
   = & \int_{\Omega} A(x, \, w, \, \nabla w) \cdot \nabla u_L - b(x, \, w, \, \nabla w) \cdot \nabla u_L \, \chi_{\Omega_{L}} \, dx\nonumber\\
%& \quad = \int_{\Omega} \Big(a(x, \, w, \, \nabla w) \,:\, \ \nabla \rho - h(x,w,\nabla w)\Big) \, u_L + \beta_h(x, \, w, \, \nabla w) \, \chi_{\Omega_L} \, dx \, \\
= & \int_{\Omega} u_{L} \, d\calm - \int_{\Omega_{L}} b \cdot \nabla u + B \,  \, u_L  \, dx \,.
\end{align}
 It also easily follows that
\begin{align}\label{laberlingu2}
	\frac{\nu}{2} \, \int_{\Omega_{L}} |\nabla u|^{p} \, dx  \leq \int_{\Omega} u_{L} \, d\calm + \int_{\Omega_{L}} |b|^{p^\prime} - B \,  \, u_L  \, dx \,.
\end{align} 
Now the proof is essentially the same as in Lemma \ref{linfty}, but here the measure $\calB$ is defined as $d\calB = |B| \, dx + d\calm$. Comparing with \eqref{Brate} and using Lemma \ref{measure}, we have
\begin{align*}
\calB(B_{\rho} \cap \Omega) \leq  c \, [B]_{1,\alpha_{N+1}^B} \, \rho^{n-p + \theta_{N+1}^B \, p} + C\,  \rho^{n-p + \theta \, p} \, .
\end{align*}
Therefore we can argue as in Lemma \ref{linfty} to obtain that
	\begin{align}\label{grundlage2}
	\int_{\Omega_L}  u_{L} \, d(\calb+\calB) \leq [(\calb + \calB)(\Omega_L)]^{\sigma} \, ,\quad \sigma = 1+\frac{1}{p}-\frac{1}{\kappa_{N+1} \, p} > 1 \, ,
\end{align}
with $\kappa_{N+1} = 1 + \min\{\theta \, p, \theta_{N+1}^b \, p, \, \theta_{N+1}^B \, (p-1)\}/(n-p)$.
The claim follows using the standard argument of Lemma 5.1 in \cite{MR0244627}, chapter 2. A lower-bound for $u$ is proven similarly.
\end{proof}
\begin{coro}
 Under the assumptions of Corollary \ref{coroN+1}, {\rm (C1)} and {\rm (C3$(N+1)$)} we also obtain that $\|u\|_{C^{\lambda_{N+1}}(\overline{\Omega})} \leq M_3$ with $\lambda_{N+1 }> 0$ and $M_3$ dependent on the solution only via $\|w\|_{L^{\bf q}}$.
\end{coro}
\begin{proof}
We follow the argument of Lemma \ref{contihold}, where for $i = N+1$ the term $B_i$ is to be replaced by 	$B_{N+1} - \Pi$. Then, the only point in the proof that has to be investigated corresponds to the inequality \eqref{quatri}. For $B = B_{N+1}$ we apply the same argument and using (B3$({\bf q},N+1)$) we get
\begin{align}\label{quatriN+1}
	\Big| \int_{B_{r}(y)} B \, (v^0-w) \, dx \Big| \leq & \|v^0-w\|_{L^{\kappa p}(B_r)} \, [\calB(B_r)]^{1-\frac{1}{\kappa p}}  \nonumber\\
	\leq & C \, \|B\|_{L^{1,\alpha^B}_M(\Omega)} \, \|\nabla(v^0-w)\|_{L^{p}(B_r)}\,  [\calB(B_r)]^{1-\frac{1}{\kappa p}}\nonumber \\
	\leq & \epsilon \,  \,\int_{B_r(y)}   |\nabla v^0|^{p_i}+ |\nabla w|^{p_i}  \, dx + c_{\epsilon} \, \|B\|_{L^{1,\alpha_{N+1}^B}_M}^{p^{\prime}} \, r^{n-p+\theta_{N+1}^B \, p^{\prime} \, p} \, ,
\end{align}
For the term involving $\Pi$ we again rely on the Sobolev measure embedding theorem. With $\kappa = 1+\theta \, p/(n-p)$ we get
\begin{align}\label{quatripi}
	\Big| \int_{B_{r}(y)} \Pi \, (v^0-w) \, dx \Big| \leq & \|\nabla u-v\|_{L^{\kappa p}(B_r;d\calm)} \, [\calm(B_r)]^{1-\frac{1}{\kappa p}}  \nonumber\\
	\leq & c \, \|\nabla(v^0-w)\|_{L^{p}(B_r)}\, r^{(n-p+\theta \ p) \, (1-\frac{1}{\kappa p})} \nonumber \\
	\leq & \epsilon \,  \,\int_{B_r(y)}   |\nabla v^0|^{p} + |\nabla w|^p \, dx + c_{\epsilon} \, r^{p^{\prime} \, (n-p+\theta \ p) \, (1-\frac{1}{\kappa \, p})} \, ,
\end{align}
and we can verify that the latter exponent is nothing else but $n-p+\theta \, p^{\prime}p$.

Thus, we are able to prove by the same methods as in Lemma \ref{contihold} an estimate
\begin{align}
	\label{migio2} \int_{B_{R}(x^0) \cap \Omega} |\nabla u|^{p} \, dx \leq C \, R^{n-p+\lambda_{N+1} \, p} \quad \text{ for all } \quad R>0 \text{ and } x^0 \in \overline{\Omega}\, .
\end{align}
The H\"older continuity follows.
\end{proof}

\section{The quasilinear case with $p_1 = \ldots = p_{N+1} =2$}\label{ReduSec} 

Throughout this section we assume that the fields $a_i$, $b_i$ and $B_i$ satisfy, in addition to (A), (B), the conditions (A1') and (B1'). As a consequence of (A1'), we choose $p_1 = \ldots = p_{N+1} = 2$. Moreover, we assume that the domain $\Omega$ and the partition of its boundary satisfy the condition (C1').
\begin{lemma}\label{reducedlem}
	We adopt the assumptions and notations of Theorem \ref{reduced}. Then every weak solution $w = (\rho,u) \in W^{1,{\bf p}} \cap L^{\bf q}$ satisfies $\|\nabla \rho\|_{L^{(t_1, \ldots,t_N)}(\Omega)} \leq M_3$ with $M_3$ depending on $w$ only via $\|w\|_{L^{\bf q}}$.  
\end{lemma}
\begin{proof}
	We employ well-known localisation arguments for linear theory\footnote{A detailed proof is to find for instance in Lemma 28 of \cite{MR4350569}.}.
	We let $x^0 \in \Omega$ be arbitrary and define $a_i^0 := a_i(x^0, \, w(x^0))$ for $i=1,\ldots,N+1$. For $i=1,\ldots,N$, the functions $\rho_i$ satisfies
\begin{align}
	-\divv (a_i^0 \, \nabla \rho_i) = \divv[(a_i - a_i^0) \, \nabla \rho_i + b_i] - B_i \, . 
\end{align}
If $(\rho, \, u)$ is a weak solution to the system \eqref{systemcomp}, it is  H\"older--continuous with exponent $\lambda_0 := \min {\bf \lambda}$ according to the Theorem \ref{MAIN}. Thus, the function $x \mapsto a_i(x, \, w(x))$ is H\"older--continuous with exponent $\alpha\, \lambda_0$ over $\overline{\Omega}$ with the $\alpha$ of condition (A1'), and near $x^0$ we have
\begin{align*}
|a_i(x,w(x)) - a_i(x^0,w(x^0))| \leq [a_i]_{\alpha} \, (|x-x^0|^2 + |w(x)-w(x^0)|^2)^{\frac{\alpha}{2}}\leq c \, |x-x^0|^{\alpha\lambda_0} \, .
\end{align*} 
Localising the problem near $x^0$ by means of a smooth cut-off function $\zeta$ we obtain in a certain neighbourhood $\Omega_R(x^0)$ for the function $\zeta \, \rho_i$ an estimate in $W^{1,\tilde{t}_i}(\Omega_R(x^0))$. Here we observe that $t_i < p_{\Omega}$, by definition, and we cannot obtain estimates for indices larger than $p_\Omega$. 

Then, we glue the local estimates by means of well-known arguments. The final estimate reads 
\begin{align*}
	\|\nabla \rho_i\|_{L^{t_i}(\Omega)} \leq c \, (1+[a_i(\cdot,w(\cdot))]_{C^{\alpha\lambda_0}})^{\frac{1}{\alpha\lambda_0}}\, (\|\nabla \rho_i\|_{L^{2}}+ \|b_i\|_{L^{p_i^b}} +\|B_i\|_{L^{np_i^B/(n+p_i^B)}} + \|\rho_i^S\|_{W^{1,p_i^S}}) \, .
\end{align*}
Since we assume (B3'$(i)$) and (C3$(i)$), all quantities on the right-hand side are controlled by the claimed quantities.
\end{proof}
\begin{rem}\label{reducedremark}
A first obvious corollary of Lemma \ref{reducedlem} and (B0) is that $\Pi $ belongs to $L^{\frac{\min_{i=1,\ldots,N} t_i}{2}}(\Omega)$ with the norm estimated by a function of $M_3$.
\end{rem}
\begin{proof}
	As for \eqref{CRHS} we find that $|\Pi| \leq \max\omega \, \sum_i|\nabla \rho_i|^2$ and the claim follows.
\end{proof}
We next obtain the higher regularity of $\nabla u$.
\begin{lemma}\label{reducedu}
	We adopt the assumptions and notations of the Theorem \ref{reduced}. Then there is $M_4 > 0$ such that $\|\nabla u\|_{L^{t_{N+1}}(\Omega)} \leq M_4$ for every weak solution such that the energy dissipation property is valid.
\end{lemma}
\begin{proof}
	This is the same proof as for Lemma \ref{reducedlem}, but now the right-hand side is $-B_{N+1} + \Pi$ which, taking into account (B3') and the Remark \ref{reducedremark}, belongs to $L^s$ with $s$ being the minimum of $np^B_{N+1}/(n+p^B_{N+1})$ and $t_{\rm min}/2$.
\end{proof}

\paragraph{\bf Uniqueness in the small}

%Suppose that all assumptions of Theorem \ref{reduced} are valid, and $t_{\min} > n$. 
In this paragraph we prove the uniqueness of the solution under the conditions of Corollary \ref{unique}.
%if the data of the problem are not too large. The relevant data are the size of the domain, the structure of the boundary and the norms of source terms, which size is expressed by the quantity $Q_0 := \|\psi_0\|_{L^1} + \|(\psi_1,\ldots,\psi_{N+1})\|_{L^{\bf k}} + \|(\phi_1,\ldots,\phi_{N+1})\|_{L^{\bf K}}$ with the data of  condition (B1). 
Recall that Lemma \ref{ENERGYEST} states that the weak norm $\|\rho\|_{W^{1,\bf p}} + \|u\|_{W^{1,\sigma_0}}$ tends uniformly to zero for each weak solution as $Q_0$ tends to zero. Hence, the fact that $\min_{i=1,\ldots,N+1} t_i >n$, we can interpolate between $L^{(p_1,\ldots,p_N,\sigma_0)}$ and $L^{\bf t}$ to obtain that
\begin{align}\label{DATA}
\sum_{i=1}^{N+1} \|w_i\|_{W^{1,n}(\Omega)} \rightarrow 0 \quad \text{ for } \quad Q_0 \rightarrow 0 
\end{align}
uniformly for every weak solution of class $W^{1,\bf p} \cap L^{\bf q}$. 
%
%For the uniqueness in the small, we assume that the dependence in $(w, \, z)$ is of class $C^1$ for all coefficient functions $A_{ij}$, $B_i$ and $\Pi$. We moreover assume that
%\begin{align*}
%|\partial_w \Pi(x,w,z)| \leq c_1(w) \, |z|^2 \quad \text{ and } \quad |\partial_z \Pi(x,w,z)| \leq c_2(w) \, |z| \, , 
%\end{align*}
%with continuous functions $c_i$ on $\mathbb{R}^{N+1}$ for $i=1,2$.
For simplicity, we shall expose only the standard proof idea for the case that $b_i = 0$, that $B_i$ is independent on $z$ and that the matrix $B_w$ is positive semidefinite. The general procedure is to linearise the coefficients as
\begin{align*}
& A(x,w,z) -A(x,\hat{w},\hat{z})\\
=& \int_0^1 \partial_z A(x, w,\, \theta\, z+(1-\theta)\hat{z}) \, d\theta \, (z-\hat{z}) +  \int_0^1 \partial_w A(x, \, \theta\, w+(1-\theta)\hat{w}, \, \hat{z}) \, d\theta \, (w-\hat{w}) \\
=:& \widetilde{\partial_z A}(x,w,z,\hat{z}) \, (z-\hat{z}) + \widetilde{\partial_w A}(x,w,\hat{w},\hat{z}) \, (w-\hat{w}) \, . 
\end{align*}
For ease of writing, we will adopt the notation $\widetilde{\partial_z A}$ for $\widetilde{\partial_z A}(x,w,z,\hat{z})$, etc.\

We next assume that $w$ and $\hat{w}$ are two weak solutions with the energy dissipation property for \eqref{systemcomp}, we obtain by means of standard arguments identities
\begin{align*}
&	\int_{\Omega} \widetilde{A_{z}} \nabla (w-\hat{w}) \, :\,  \nabla (w-\hat{w}) + \widetilde{A_{w}}\,  (w-\hat{w}) \cdot \nabla (w-\hat{w}) +  \widetilde{B_{w}}\,  (w-\hat{w}) \, (w-\hat{w}) \, dx \, \\
&\quad	= \int_{\Omega} \delta_{N+1} \, (\widetilde{\Pi_w} \, (w-\hat{w}) + \widetilde{\Pi_z} \nabla (w-\hat{w})) \, (u-\hat{u}) \, dx 
\end{align*}
With the H\"older and Sobolev-Poincar\'e inequalities, we estimate
\begin{align*}
 \Big|	\int_{\Omega}  \widetilde{A_{w}}\,  (w-\hat{w}) \cdot \nabla (w-\hat{w}) \, dx\Big|\leq& \|\widetilde{A_{w}}\|_{L^n} \, \|\nabla (w-\hat{w}) \|_{L^2}\, \|w-\hat{w} \|_{L^{2^*}}  \\
 \leq& c_{PS} \, \|\widetilde{A_{w}}\|_{L^n} \, \|\nabla (w-\hat{w}) \|_{L^2}^2\, ,
 %\\
%\Big|\int_{\Omega}  \widetilde{B_{w}}\,  (w-\hat{w}) \, (w-\hat{w}) \, dx\Big| \leq& \|\widetilde{B_{w}}\|_{L^{\frac{3}{2}}} \,  \|w-\hat{w} \|_{L^6}^2 \leq c_{PS}^2 \, \|\widetilde{B_{w}}\|_{L^{\frac{3}{2}}} \, \|\nabla (w-\hat{w}) \|_{L^2}^2 \, ,
\end{align*}
and for the critical term similarly
\begin{align*}
\Big|\int_{\Omega} \widetilde{\Pi_w} \, (w-\hat{w})  \, (u-\hat{u})) \, dx \Big| \leq c_{PS}^2 \, \|\widetilde{\Pi_{w}}\|_{L^{\frac{n}{2}}} \, \|\nabla (w-\hat{w}) \|_{L^2}^2\, ,\\
\Big|\int_{\Omega} \widetilde{\Pi_z} \, \nabla (w-\hat{w})  \, (u-\hat{u})) \, dx \Big| \leq c_{PS}\, \|\widetilde{\Pi_{z}}\|_{L^n} \, \|\nabla (w-\hat{w}) \|_{L^2}^2\,.
\end{align*}
Hence, exploiting that $B_w$ is a positive matrix, we can easily derive an inequality
\begin{align*}
\min \nu \, \int_{\Omega} |\nabla (w-\hat{w})|^{ 2} \, dx \leq & c_{PS} \, \Big[ \|\widetilde{A_{w}}\|_{L^n} +\|\widetilde{\Pi_{z}}\|_{L^n} + c_{PS}  \, \|\widetilde{\Pi_{w}}\|_{L^\frac{n}{2}}\Big] \, \times \\
& \times \,  \int_{\Omega} |\nabla (w-\hat{w})|^{2} \, dx\, .
\end{align*}
If $w-\hat{w} \neq 0$, it then must hold that
\begin{align*}
	\min \nu \, \leq c_{PS} \, \Big[ \|\widetilde{A_{w}}\|_{L^n} +\|\widetilde{\Pi_{z}}\|_{L^n} + c_{PS}  \,  \|\widetilde{\Pi_{w}}\|_{L^\frac{n}{2}}\Big] \, .
\end{align*}
For the simple growth assumptions \eqref{Piderivs}, we have
\begin{gather*}
\|\widetilde{A_{w}}\|_{L^n} \leq \max_{|w| \leq R } \|a_w(\cdot,w)\|_{L^\infty(\Omega)} \, (\|\nabla w\|_{L^n}+\|\nabla \hat{w}\|_{L^n}) \,,\\
\|\widetilde{\Pi_{z}}\|_{L^n} \leq (\max_{|w| \leq R } c_2(w))^{\frac{1}{n}} \, (\|\nabla w\|_{L^n}+\|\nabla \hat{w}\|_{L^n}), \, \|\widetilde{\Pi_{w}}\|_{L^\frac{n}{2}} \leq (\max_{|w| \leq R } c_1(w))^{\frac{2}{n}} \, (\|\nabla w\|_{L^n}+\|\nabla \hat{w}\|_{L^n})^2 \, .
\end{gather*}
Here $R$ is an upper bound for the values of weak solutions. Now we see that if there are two different weak solutions, they satisfy
\begin{align*}
0< c_0(\nu,C_{PS},R) \leq \|\nabla w\|_{L^3}+\|\nabla \hat{w}\|_{L^3} \, ,
\end{align*}
which cannot occur for sufficiently small $Q_0$ due to \eqref{DATA}.

\section{Illustration I: The thermistor case with a heat source}\label{Thermite}

For the thermistor problem $\Omega$ represents some ohmic conductor. We choose $N = 1$ and re-denote $\rho_1 := \varphi$ representing the electric potential, while $u := T$ is the relative temperature (that is, $T = T^{\rm abs} - T_0$ with $T_0 > 0$ constant). For the simplest case, we assume scalar positive temperature-dependent electrical and heat conductivity coefficients $\sigma = \sigma(T)$ and $\kappa = \kappa(T)$, and then, the thermistor system \eqref{systemcomp} assumes the form
\begin{align}
 \label{Thermist1} -\divv (\sigma(T) \, \nabla \varphi) = & 0 \\
\label{Thermist2} -\divv(\kappa(T) \, \nabla T) = & \sigma(T) \, |\nabla \varphi|^2 + h(x) \, .
\end{align}
For this simple variant, we choose the variables $\rho_1 = \rho := \varphi$ and $u = k(T)$ where $k$ is the primitive function of $\kappa$ vanishing in zero. With $\hat{\sigma} := \sigma \circ k^{-1}$, we make the following choices 
\begin{align}\label{thermiteconsti}
&	a_1(x, \, \rho, \, u, \, \nabla \rho, \, \nabla u) = \hat{\sigma}(u) \, \nabla \rho \, ,\quad b_1 \equiv 0 \,, \quad B_1 \equiv 0 \, , \nonumber\\
&	a_2(x, \, \rho, \, u, \, \nabla \rho, \, \nabla u) =  \nabla u, \, \quad b_2 = 0, \, \quad B_2 = -h(x) \, ,\nonumber\\
& \Pi(x,, \, \rho, \, u, \, \nabla \rho, \, \nabla u) = \hat{\sigma}(u) \, |\nabla \rho|^2 \, .
\end{align}
Then, we easily verify that \eqref{systemcomp} is an equivalent form of \eqref{Thermist1}, \eqref{Thermist2}. We consider the mixed boundary conditions \eqref{BC}, \eqref{BC2}. The boundary data for $u$ are defined as $u^S = k^{-1}(T^S)$.
\begin{lemma}\label{thermiste}
We assume that $\sigma \in C^{1}(\mathbb{R})$ and that there are $0< \sigma_0 < \sigma_1$ such that $\sigma_0 \leq \sigma(u) \leq \sigma_1$ for all $u \in \mathbb{R}$. We moreover assume that $h \in L^p(\Omega)$ with $p > n/2$. We consider Dirichlet data $\varphi^S, \, k^{-1}(T^S) \in W^{1,p^S}(\Omega)$ with $p^S > n$. If $\Omega$ and $S$ are subject to the condition {\rm (C1)} and $p_\Omega > n$ in {\rm (C1')}, every weak solution to the mixed boundary value problem \eqref{systemcomp}, \eqref{BC}, \eqref{BC2} for \eqref{thermiteconsti} belongs to $W^{1,q}(\Omega)$ for a $q > n$. If moreover $\|h\|_{L^1} + \|\nabla \varphi^S\|_{L^2} + \|\nabla k^{-1}(T^S)\|_{L^2}$ is sufficiently small, the weak solution is unique.
\end{lemma}
\begin{proof}
The claims follows from the Theorem \ref{reduced} and the Corollary \ref{unique} if we can verify the assumptions {\rm (A)} and {\rm (B)}. The remaining assumptions {\rm (A')}, {\rm (B')}, {\rm (C)} and \eqref{Piderivs} are obviously valid.

We choose $p_1 = 2 = p_2$. Since $\hat{\sigma}$ is uniformly positive and bounded, the conditions (A1) and (B0) are valid.
We verify (B1) with ${\bf q} = 0$. Note that in the present context, the functions $\Psi_i(x,w,z)$ have the structure
\begin{align*}%\label{PSI}
	\Psi_1 =  \hat{\sigma}(w_2) z_1 \cdot \nabla \rho_i^S(x) \, , \quad
	\Psi_2 :=	& \delta_0 \,  \frac{[z_2 \cdot \nabla u^S(x)]^+}{(1+|w_2-u^S(x)|)^{1+\delta_0}} + [h(x) \, {\rm sign}(w_2-u^S(x))]^+  \, .  
\end{align*}
By choosing $\psi_0 = h$, and for $i = 1,2$, $\psi_i = 0$ and $\phi_i(x) = |\nabla w_i^S(x)|$, we verify easily the estimate of (B1$(\bf 0)$) $K_i = 2$, $s_i = 1$ and $t = 1+\delta_0 > 1/K_2^{\prime}$. In order to satisfy (B2) and (B3), we moreover assume that $h \in L^p(\Omega)$ with $p > n/2$. Then (B3') is also valid, and the claim follows.
\end{proof}
We remark that Lemma \ref{thermiste} gives an illustration for the simplest constitutive equations. Much more general choices than \eqref{thermiteconsti} are admissible for the regularity in the thermistor system.

\section{Stationary Nernst-Planck equations with variable temperature: Existence of continuous weak solutions}\label{NPPSEC}

In this section we discuss the proof of Theorem \ref{NPPtheo}.

In this case, $\Omega$ represents a container filled with an electrolyte. In the simplest case, this is a fluid mixture of a fluid solvent, a salt and ions of $N-1 \in \mathbb{N}$ different species dissolved therein. In the pure Nernst--Planck model, the fluid solvent is assumed non-participating and is not present in the model. Moreover the electrolyte is assumed electrically neutral.

For the isothermal Nernst--Planck model the main variables are often chosen to be the concentrations $c_1, \ldots,c_{N-1}$ of the ions and the electrical potential $\varphi$. For the non-isothermal case, we prefer replacing the ion concentrations by the quotients of chemical potential over temperature, which in the thermodynamic framework are the derivatives $\partial_{\bf c} s$ of the entropy density function (entropic variables). 

To address this model within our setting, for $i=1,\ldots,N-1$, we let $\rho_i := \ln c_i + g_i^0(T)/T$, with given concave twice differentiable functions $g_i^0$ of temperature. The quantities $T \, \rho_i$, $i=1, \ldots,N-1$ correspond to the chemical potential of the species for an ideal mixture\footnote{The chemical potential or Gibbs variables are accordingly of the form $T \, \ln c_i + g_i^0(T)$.}. We further choose $\rho_N:= \varphi$ to be the electrical potential and $u = T$ is the absolute temperature. 

For $i = 1,\ldots,N-1$, we denote by $ m_i(\rho,\, u)$ the thermodynamic diffusivity (Fick--Onsager coefficient) and by $\upsilon_i = z_iq_0$ the constant charge number of an ion species. Further we let $\varepsilon(\rho,\, u)$ be the permittivity of the electrolyte\footnote{This means that $\varepsilon(\rho,\, u) = \epsilon_0 \, (1+\chi(\rho,u))$ with the constant vacuum permittivity $\epsilon_0$ and the state-dependent dielectric susceptibility $\chi$.} and $\kappa(\rho,\, u)$ the heat conductivity coefficient of the electrolyte. We consider the PDE system
\begin{gather}
	\label{ions} -\divv \Big[ m_i(\rho,\, u) \,\Big( \nabla \rho_i + \frac{\upsilon_i}{u} \, \nabla \rho_N\Big)\Big] =  0 \, \quad \text{ for } i = 1,\ldots,N-1 \, ,\\
	\label{poisson} -\divv (\varepsilon(\rho,\, u) \, \nabla \rho_N) =  0 \, ,\\
	\label{templie}-\divv(\kappa(\rho,\, u) \, \nabla u) =  \sum_{i=1}^{N-1} \upsilon_i \, m_i(\rho,\, u) \, \Big(\nabla \rho_i + \frac{\upsilon_i}{u} \, \nabla \rho_N\Big) \cdot \nabla \rho_N \, .
\end{gather}
The right-hand side of the equation \eqref{templie} expresses the product $J \cdot E$ with the electric current due to free charges $J = -\sum_{i=1}^{N-1} \upsilon_i \, m_i(\rho,\, u) \, (\nabla \rho_i + (\upsilon_i/u) \, \nabla \rho_N)$ and $E = - \nabla \rho_N$ the electric field. We shall discuss this model under the assumption of uniformly positive and bounded permittivity: There are constants $0<\varepsilon_0 \leq \varepsilon_1$ such that 
\begin{align}\label{atleastthat}
\varepsilon_0 < \varepsilon(\rho, \, u) \leq \varepsilon_1 \quad \text{ for all } \quad (\rho,u) \in \mathbb{R}^N \times \mathbb{R}_+\, .
%\quad \text{ and } \quad \kappa_0 \leq \kappa(\rho,\, u) \leq \kappa_1 \, .
\end{align}
Then we have the following remark which is essential for the analysis. The proof follows from the considerations in Lemma \ref{linfty} and Lemma \ref{contihold}.
\begin{rem}\label{dissipphi}
Let $\varphi^S \in W^{1,p_N^S}(\Omega)$ with $p_N^S >n$, and assume that {\rm (C1)} is valid for the domain $\Omega$. For given measurable $(\rho, \, u)$ the weak solution $\varphi = \rho_N$ to the equation \eqref{poisson} with the mixed boundary conditions \eqref{BC}, \eqref{BC2} belongs to $W^{1,2}(\Omega)$ and to $C^{\lambda_N}(\overline{\Omega})$ for a $0 < \lambda_N \leq 1-n/p_N^S$ determined by the space dimension, the domain and the quotient $\varepsilon_1/\varepsilon_0$ only. Moreover, the measure $d\mu_{\varphi} = |\nabla \varphi|^2 \, dx$ satisfies for all $\rho \leq \rho_0={\rm diam}(\Omega)$ and all balls $$\mu_{\varphi}(B_{\rho} \cap \Omega) \leq c \, \|\nabla \varphi^S\|_{L^{p_N^S}}^2 \, \Big(\frac{\rho}{\rho_0}\Big)^{n-2+2\lambda_N} \,  .$$ 
\end{rem}
Then we can at first prove a regularity result under a condition on the values of weak solutions. This corresponds to part \ref{condireguNPP} of Theorem \ref{NPPtheo}. 
\begin{lemma}\label{NPPregu}
We let the assumptions of Remark \eqref{dissipphi} be valid. Moreover, we assume that $m_1, \ldots, m_{N-1}$ and $\kappa$ belong to $C^1(\mathbb{R}^N \times \mathbb{R}_+)$ and that the boundary data $\rho_i^S$ and $T^S$ belong to $W^{1,p}(\Omega)$ with a $p > n$. Let $(\rho, \, u)$ be a weak solution to the equations \eqref{ions}, \eqref{poisson}, \eqref{templie} with mixed boundary conditions \eqref{BC}, \eqref{BC2} and satisfying moreover the additional assumptions 
\begin{align}\label{zusatz2}
\sup_{i=1,\ldots,N} \|\rho_i\|_{L^{\infty}(\Omega)} + \|u\|_{L^{\infty}(\Omega)} < +\infty \quad \text{ and }\quad  T_0 := \inf_{x \in \Omega} u(x) >0 \, .
\end{align}
Then all components of $(\rho, \, u)$ are H\"older continuous on $\overline{\Omega}$. If {\rm (C1')} is valid with $p_\Omega > n$, all components belong to $W^{1,q}(\Omega)$ with $q > n$. There is at most one weak solution with \eqref{zusatz2} if $\sum_{i=1}^{N-1} \|\nabla \rho_i^S\|_{L^2} + \|\nabla \varphi^S\|_{L^2} + \|\nabla T^S\|_{L^2}$ is sufficiently small.
\end{lemma}
\begin{proof}
Since the coefficient functions $m_i$ and $\kappa$ are positive and of class $C^1$, the assumptions \eqref{zusatz2} imply that there are positive constants $0 < m^\prime_0 \leq m^\prime_1$ such that $0 < m^\prime_0 \leq m_i(\rho,u) \leq m^\prime_1$ for $i = 1, \ldots, N-1$, and likewise constants $0< \kappa_0 < \kappa_1$ such that $\kappa_0 \leq \kappa(\rho,\, u) \leq \kappa_1$.	
	
For $i=1,\ldots,N$, we describe the problem \eqref{ions}, \eqref{poisson}, \eqref{templie} by means of the the constitutive equations $B_i \equiv 0$ and
\begin{gather*}
 a_i(\rho, \, u, \, \nabla \rho) = m_i(\rho,u) \, \nabla \rho_i  \text{ for } i = 1,\ldots,N-1, \quad a_N(u, \, \nabla \rho) = \varepsilon(\rho,u) \, \nabla \rho_N \, ,\\
 b_{i}(x, \, u) = \frac{\upsilon_i}{u} \, \nabla \varphi(x) \text{ for } i = 1,\ldots,N-1\, , \quad b_N \equiv 0 \, .
\end{gather*}
For the energy equation we let $b_{N+1} \equiv 0$, $B_{N+1} \equiv 0$ and
\begin{align*}
 a_{N+1}(u, \, \nabla u) = \kappa(\rho, u) \, \nabla u \, ,
 \quad \Pi(\rho, \, u,\nabla \rho) :=  \sum_{i=1}^{N-1}\, m_i(\rho,u) \, \Big(\nabla \rho_i + \frac{\upsilon_i}{u} \, \nabla \rho_N\Big) \, \upsilon_i \, \nabla \rho_N  \, .
\end{align*}
%Here $ m_i(\rho,u)$ denote the mobility functions of the ions and $z_i \in \mathbb{Z}$ are the constant charge numbers. These choices give rise to the PDE system 
%\begin{gather}
%\label{ions} -\divv \Big[ m_i(u,\rho) \,\Big( \nabla \rho_i + \frac{z_i}{u} \, \nabla \rho_N\Big)\Big] =  0 \, \quad \text{ for } i = 1,\ldots,N-1 \, ,\\
% \label{poisson} -\divv (\epsilon(u,\rho) \, \nabla \rho_N) =  0 \, ,\\
% \label{templie}-\divv(\kappa(u) \, \nabla u) =  \sum_{i=1}^{N-1} z_i \, m_i(u,\rho) \, \Big(\nabla \rho_i + \frac{z_i}{u} \, \nabla \rho_N\Big) \cdot \nabla \rho_N \, .
% \end{gather}
% The electric potential is moreover subject to the Poisson equation
% \begin{align}\label{poisson}
% -\divv (\epsilon(T,\rho) \, \nabla \phi) = 0 \quad \text{ in } \quad \Omega \, ,
% \end{align}
%with the dielectric function $\epsilon$ of the electrolyte. 
We set $p_i = 2$ for $i = 1,\ldots,N+1$, and turn to verify the conditions (A), (B) and (A'), (B'). 

At first notice that the verification of (A1) is trivial, since we have shown that the mobilities in all equations are positive and bounded for the bounded weak solutions.
%In order to verify (A1) for $i < N$, the assumption that the mobilities are bounded functions seems unnatural. In fact, the usual choice for the Nernst--Planck equations would be $m_i(T,\ln c) = d_i(T) \, c_i$, with the diffusivity $d_i$ being a positive bounded function of temperature. Hence, we will at first restrict to proving the results under the assumption that bounds from below and above are already known for the species concentrations: We assume that there are $0< n_0 < n_1$ such that the weak solution satisfies
%\begin{align}\label{zusatz}
% n_0 \leq c_i(x) = e^{\rho_i(x)-\mu_i^0(u(x))}\leq n_1 \quad \text{ for all } \quad x \in \Omega \, . 
%\end{align} 
%Then we have also $0 < m^\prime_0 \leq m_i(\rho,u) \leq m^\prime_1$ for appropriate numbers $m^\prime_0$ and $m^\prime_1$ which depends on the thresholds of the diffusivities and on the numbers $n_0$ and $n_1$ of \eqref{zusatz}.

In order to next check (B0), we employ the Young inequality to see that
\begin{align*}
|	\Pi | \leq \sum_{i=1}^{N-1} \Big(\frac{u \, m_i}{2} \, |\nabla \rho_i|^2 + \frac{3m_i \, \upsilon_i^2}{2\,u} \, |\nabla \rho_N|^2\Big) \, .
\end{align*}
Hence we can choose $\omega_i = \|u\|_{L^{\infty}} \, m_1^\prime/2$ for $i =1,\ldots,N-1$ and $\omega_{N} = 3/(2T_0)\, m^\prime_1 \, \sum_{i=1}^{N-1}\, \upsilon_i^2$ to satisfy (B0). 

Next we want to check (B1) and (B2). 
%In addition to the conditions \eqref{zusatz}, we here need to assume that the temperature is strictly bounded from below
%\begin{align}\label{zusatz2}
%T_0 := \inf_{x \in \Omega} T(x) >0 \, .
%\end{align}
Here we make use of the simple structure of the equation \eqref{poisson} and the Remark \ref{dissipphi}. 
%Here $b_N = 0 = B_N$ so that standard methods allow to show that there is $\lambda_N \in ]0,1[$ such that $\|\rho_N\|_{C^{\lambda_N}} + \|\nabla \rho_N\|_{L^{2,\alpha_N}_{\rm M}}$ are bounded independently, with $\alpha_N := n -2 \, (1-\lambda_N)$.
To verify (B1(${\bf q}$)) with ${\bf q} = +\infty$, we have to consider
\begin{gather*}
 \Psi_i(x,w,\nabla w) = - \frac{m_i(\rho,u)\, \upsilon_i}{u} \, \nabla \varphi(x)\cdot \nabla \rho_i+ m_i(\rho,u) \, (\nabla \rho_i +\frac{\upsilon_i}{u} \nabla \rho_N) \cdot\nabla \rho_i^S(x) \, ,\\
 \Psi_N(x,w,\nabla w) = \varepsilon(\rho,u) \, \nabla \rho_N \cdot \nabla \varphi^S(x)\, ,\\
  \Psi_{N+1}(x,w,\nabla w) = \kappa(\rho,u) \, \delta_0 \,\frac{[ \nabla u \cdot \nabla u^S(x)]^+}{(1+|u-u^S(x)|)^{1+\delta_0}} \, .
 \end{gather*}
We let $\phi_i(x) := m_1^{\prime} \, (|\upsilon_i|/T_0 \, |\nabla \varphi(x)|+|\nabla\rho_i^S(x)|)$ for $i=1,\ldots,N-1$, $\phi_N \equiv 0$ and $\phi_{N+1}(x) = \kappa_1 \, \delta_0 \, |\nabla u^S(x)|$. Moreover we define $$\psi_0(x) = m_1^\prime \, |\upsilon| \, T_0^{-1} \, |\nabla \varphi(x)|\, \sum_{i=1}^{N-1} |\nabla\rho_i^S(x)| + \varepsilon_1 \, |\nabla \varphi(x)|\, |\nabla \varphi^S(x)|\, .$$ Choosing $s_i = 1$ and $K_i = 2$ for $i = 1,\ldots,N+1$, the condition (B1($\bf \infty$)) is established with $t = 0$. Moreover, the parameters $K_i$ and $s_i$ satisfy the assumptions of the Lemma \ref{ENERGYEST}, where we can choose $q_1,\ldots, q_{N+1} = + \infty$.

%Similarly, we have
%\begin{align*}
%	B_{N+1} \leq |\nabla T| \,  |\nabla \phi(x)|\,\sum_{i=1}^{N-1} m_i \, |z_i| \, |\partial_T\mu_i^0(T)|  \, . 
%\end{align*}
%Hence with the choice $s_{N+1} = 1$ and $\phi_{N+1}(x) = |\nabla \varphi(x)|\,\sum_{i=1}^{N-1} m_i \, |z_i| \, |\partial_T\mu_i^0(T)|$. Suppose now that
%\begin{align*}
%|\partial_T\mu_i^0(T)| \leq C_0 \, (1+T)^{-t} \quad \text{ with } \quad t > \frac{1}{2} \, .
%\end{align*}
%Then we attain $|B_{N+1}| \leq \phi_{N+1}(x)  \, |\nabla T|/(1+T)^t$. We choose $K_{N+1} = 2 > 1/t$. Overall, we can verify (B1) with $\psi_i = 0$ and all $\phi_i$ bounded in $L^2$.

In order to verify (B2(${\bf +\infty}$)) and (B3(${\bf +\infty},i$)) we apply directly the Morrey--regularity of $\varphi$ in Remark \ref{dissipphi} and the fact that $|b_i| \leq m_1^{\prime}|\upsilon_i|/T_0 \, |\nabla \varphi|$. We obtain (B3(${\bf +\infty},i$)) with $\theta_i^b = \lambda_N$ for $i =1,\ldots,N$. 
Hence, the Theorem \ref{MAIN} now guarantees that $(\rho,\, u)$ is H\"older continuous. 

Next we return in the equation \eqref{poisson}, and suppose that (C1') is valid with $p_{\Omega} > n$. Since the coefficient $\varepsilon(x) =\varepsilon(\rho(x),u(x))$ is H\"older continuous, we obtain a bound for $\nabla \rho_N$ in $L^{q}$, with $n < q \leq p_N^S$ and $q < p_{\Omega}$ arbitrary. This gives next $b_i \in L^{q}$ for $i =1,\ldots,N-1$. Hence (B3'$(i)$) is satisfied as well, with $p_i^b >n$. Theorem \ref{reduced} now guarantees that $\nabla \rho_i\in L^{\tilde{q}}$ with $\tilde{q} = \min\{q,p_i^S\} > n$. The right-hand side of the energy equation \eqref{templie} is now found to belong to $L^{\tilde{q}/2}$. Thus, the temperature gradient is also integrable to the power $\min\{\tilde{q}, p_{N+1}^S\}$, and the uniqueness principle is valid in the small.
\end{proof}
As to the proof of the additional Remark \ref{nppextend}, replacing the homogeneous equation \eqref{poisson} by 
\begin{align}\label{poisson2}
	-\divv( \varepsilon(\rho,u) \, \nabla \rho_N) = \sum_{i=1}^N \upsilon_i \, c_i = \sum_{i=1}^{N-1} \upsilon_i \, \exp(\rho_i - g_i^0(u)/u) =: \varrho^F \, , 
\end{align}
does not affect qualitatively the result of Lemma \ref{NPPregu}. Still under the condition \eqref{atleastthat}, the estimate in Remark \ref{dissipphi} must now reads
$$\mu_{\varphi}(B_{\rho} \cap \Omega) \leq c \, (\|\nabla \varphi^S\|_{L^{p_N^S}}^2+ \| \varrho^F\|_{L^{1,n-2+2\lambda_N}_M}) \Big(\frac{\rho}{\rho_0}\Big)^{n-2+2\lambda_N} \,  ,$$ 
and the assumption \eqref{zusatz2} allows to estimate
\begin{align*}
	\| \varrho^F\|_{L^{1,n-2+2\lambda_N}_M} \leq c(\Omega, \, T_0, \, \|u\|_{L^{\infty}}, \, \sup_{i=1,\ldots,N-1} \|\rho_i\|_{L^\infty}) \, .
\end{align*}
Thus, we obtain the H\"older continuity exactly as in the pure Nernst--Planck case.\\

Next we would like to discuss the validity of the conditions \eqref{zusatz2} and to prove part \eqref{ExiNPP} of Theorem \ref{NPPtheo}. For this aim we need more structure assumptions on the coefficients. Concerning the Fick--Onsager coefficients the direct assumption that $m_i$ is bounded and uniformly positive would be unusual in the context of the Nernst-Planck equations. Instead we will consider the case that $m_i = d_i(T) \, c_i$ where the diffusivity $d_i$ is a positive function of temperature. Hence, in our notations,
\begin{align}\label{mob1}
m_i(\rho,u) = d_i(u) \, \exp(\rho_i - g_i^0(u)/u) = \tilde{d}_i(u) \, e^{\rho_i} \, .
\end{align}
\begin{lemma}\label{boundvalues}
Let $(\rho, u)$ be a weak solution. Suppose that $m_i$ obeys \eqref{mob1} with $0< d_0^\prime\leq \tilde{d}_i(u)\leq d_1^\prime$ uniformly positive and bounded for all $u \in \mathbb{R}_+$, then $\|\nabla \rho_i\|_{L^2(\Omega)} + \|\rho_i\|_{L^{\infty}(\Omega)}< \infty$ for $i = 1,\ldots,N-1$.
\end{lemma}
\begin{proof}
{\bf Bounds in $L^2$ norms.} At first we show a bound for $\rho_i$ in $L^2$.
Let $n_i := e^{\rho_i}$ and $\delta_0 >0$. We test the equation \eqref{ions} with $ 1-(n_i^S+\delta_0)/(n_i+\delta_0)$. Invoking also Young's inequality we first show that 
\begin{align*}
	\int_{\Omega} \tilde{d}_i(u) \, (n_i^S+ \delta_0)\, |\nabla \ln(n_i+\delta_0)|^2 \, dx \leq 2 \, \int_{\Omega} \tilde{d}_i(u) \, \Big(|\nabla \sqrt{n_i^S+\delta_0}|^2 + \frac{|\upsilon_i|^2}{u^{2}} \, (n_i^S+\delta_0) \, |\nabla \varphi|^2\Big) \, dx \, .
\end{align*}
Thus, since $n^S_i \geq e^{-\|\rho_i^S\|_{L^{\infty}}}$ it follows that
\begin{align}\label{crunch}
	\|\nabla \ln(n_i+\delta_0)\|^2_{L^2} \leq 2 \, e^{\|\rho_i^S\|_{L^{\infty}}} \, \frac{d_1^\prime}{d_0^\prime}\,  \Big(\|\nabla \sqrt{n_i^S+\delta_0}\|_{L^2}^2+|\upsilon_i|\, T_0^{-2} \, \|n_i^S+\delta_0\|_{L^{\infty}} \, \|\nabla\varphi\|_{L^2}^2\Big) \, .
\end{align}
Consider $g = \ln[(n_i+\delta_0)/(n^S_i + \delta_0)$ which belongs to $W^{1,2}_S(\Omega)$ (subspace of $W^{1,2}(\Omega)$ with trace zero on $S$). By means of the Sobolev and Poincar\'e Lemma we obtain  that $\|g\|_{L^2} \leq c_1 \, \|\nabla g\|_{L^2}$ and hence
\begin{align*}
	\|\ln(n_i+\delta_0)\|^2_{L^2} \leq 	& \|\ln(n_i^S+\delta_0)\|^2_{L^2}\\
	& + 2c_{1} \, e^{\|\rho_i^S\|_{L^{\infty}}}\,  \frac{d_1^\prime}{d_0^\prime}\,  \Big(\|\nabla \sqrt{n_i^S+\delta_0}\|_{L^2}^2+\frac{|\upsilon_i|}{T_0^{2}} \, \|n_i^S+\delta_0\|_{L^{\infty}} \, \|\nabla\varphi\|_{L^2}^2\Big) \, .
\end{align*}
We let $\delta_0 \rightarrow 0$ and using the definition of $n_i^S$ we obtain that
\begin{align}\label{rhointegrable}
	\|\rho_i\|^2_{L^2} \leq &
	\|\rho_i^S\|^2_{L^2} +  2c_{1} \, e^{\|\rho_i^S\|_{L^{\infty}}} \,   \frac{d_1^\prime}{d_0^\prime}\,  \Big(\|\nabla \sqrt{n_i^S}\|_{L^2}^2+\frac{|\upsilon_i|}{T_0^{2}} \, \|n_i^S\|_{L^{\infty}} \, \|\nabla\varphi\|_{L^2}^2\Big)\nonumber\\
\leq &	\|\rho_i^S\|^2_{L^2} +  2c_{1} \, e^{2\|\rho_i^S\|_{L^{\infty}}} \,   \frac{d_1^\prime}{d_0^\prime}\,  \Big(\|\nabla \rho_i^S\|_{L^2}^2+\frac{|\upsilon_i|}{T_0^{2}} \, \|\nabla\varphi\|_{L^2}^2\Big)
	 \, .
\end{align}
Note that using \eqref{crunch} and Lemma \ref{lieberman} with $\theta = \lambda_N$, and $g = \ln[(n_i+\delta_0)/(n^S_i + \delta_0)]$ we can infer in the same way the inequality
	\begin{align}\label{rhomuintegrable}
		\|\rho_i\|^2_{L^2_{d\mu_\varphi}} \leq \|\rho_i^S\|^2_{L^2_{d\mu_\varphi}} + C\,  [\mu_{\varphi}(\Omega)]^{1-\frac{1}{\kappa}} \, \|\nabla \varphi^S\|_{L^{p^S_N}}^{\frac{2}{\kappa}} \, e^{2\|\rho_i^S\|_{L^{\infty}}} \,   \frac{d_1^\prime}{d_0^\prime}\,  \Big(\|\nabla \rho_i^S\|_{L^2}^2+\frac{|\upsilon_i|}{T_0^{2}} \, \|\nabla\varphi\|_{L^2}^2\Big)\, .
		% \frac{d_1^\prime}{d_0^\prime}\,  \Big(\|\nabla \sqrt{n_i^S}\|_{L^2}^2+\frac{|z_i|}{T_0^{2}} \, \|n_i^S\|_{L^{\infty}} \, \|\nabla\varphi\|_{L^2}^2\Big) \, .
	\end{align}
A useful consequence is that for $\Omega_L = \{x \, : \, \rho_i > L\}$ we have
\begin{align}\label{smallmeasure}
	L^2\, \mu_\varphi(\Omega_L) \leq \|\rho_i^S\|^2_{L^2_{d\mu_\varphi}} +
	 C\,  [\mu_{\varphi}(\Omega)]^{1-\frac{1}{\kappa}} \, \|\nabla \varphi^S\|_{L^{p^S_N}}^{\frac{2}{\kappa}} \, e^{2\|\rho_i^S\|_{L^{\infty}}} \,   \frac{d_1^\prime}{d_0^\prime}\,  \Big(\|\nabla \rho_i^S\|_{L^2}^2+\frac{|\upsilon_i|}{T_0^{2}} \, \|\nabla\varphi\|_{L^2}^2\Big)\, .
%	 c_{PS}(\mu_{\phi}) \,  \frac{d_1^\prime}{d_0^\prime}\,  (\|\nabla \sqrt{n_i^S}\|_{L^2}^2+|z_i|\, T_0^{-2} \, \|n_i^S\|_{L^{\infty}} \, \|\nabla\varphi^S\|_{L^2}^2) \, .
\end{align}

Next we test \eqref{ions} with $n_i-n_i^S$, and this time we get
\begin{align}\label{gradl2un}
	\int_{\Omega} |\nabla n_i|^2 \, dx \leq 2\frac{d_1^\prime}{d_0^{\prime}} \, \int_{\Omega} \, \Big(|\nabla n_i^S|^2 + |\upsilon_i|^2u^{-2} \, n_i^2 \, |\nabla \varphi|^2\Big) \, dx \, .
\end{align}
With $L \geq \sup_{x \in S} \rho_i^S(x) =: L_0$ we estimate $n_i \leq 2\, ( e^{2L}+\max\{(n_i-e^L),0\}^2)$ and it follows that
\begin{align*}
	\int_{\Omega} |\upsilon_i|^2u^{-2} \, n_i^2 \, |\nabla \varphi|^2 \, dx \leq \frac{2 \,  |\upsilon_i|^2}{T_0^2} \, \Big(e^{2L} \, \|\nabla \varphi\|_{L^2}^2 + \int_{\Omega_{L}} (n_i-e^L)^2 \, d\mu_\varphi\Big) \, .
\end{align*}
On the other hand, with $n_{i,L} = \max\{n_i-e^L,0\}$, H\"older's inequality and the Lemma \ref{lieberman} imply that
\begin{align*}
	\int_{\Omega_{L}} (n_i-e^L)^2 \, d\mu_\varphi  \leq [\mu_\varphi(\Omega_{L})]^{1-\frac{1}{\kappa}} \, \|n_{i,L}\|_{L^{2\kappa}_{d\mu_\phi}}^2 \leq C \, \|\nabla \varphi^S\|_{L^{p^S_N}}^{\frac{2}{\kappa}}\, [\mu_\varphi(\Omega_{L})]^{1-\frac{1}{\kappa}} \,	\int_{\Omega_{L}} |\nabla n_i|^2 \, dx \, . 
\end{align*}
Thus, we obtain that
\begin{align*}
	\int_{\Omega} |\upsilon_i|^2u^{-2} \, n_i^2 \, |\nabla \varphi|^2 \, dx \leq \frac{2 \, |\upsilon_i|}{T_0^2} \, e^{2L} \, \|\nabla \varphi\|_{L^2}^2 + \frac{2  \, |\upsilon_i|^2 \, C \,  \|\nabla \varphi^S\|_{L^{p^S_N}}^{\frac{2}{\kappa}}}{T_0^2} \, [\mu_\varphi(\Omega_{L})]^{1-\frac{1}{\kappa}} \,	\int_{\Omega_{L}} |\nabla n_i|^2 \, dx\, ,
\end{align*}
and under the condition
\begin{align}\label{CONDI}
	\frac{2 \,  d_1^{\prime} \, |\upsilon_i|^2 \, C \, \|\nabla \varphi^S\|_{L^{p^S_N}}^{\frac{2}{\kappa}}}{d_0^{\prime} \, T_0^2} \, [\mu_\varphi(\Omega_{L})]^{1-\frac{1}{\kappa}} \leq \frac{1}{2} \, ,
\end{align}
we reach with the help of \eqref{gradl2un} the inequality
\begin{align}\label{gradl2deux}
	\int_{\Omega} |\nabla n_i|^2 \, dx \leq \frac{4d_1^{\prime}}{d_0^{\prime}}\, \Big(\int_{\Omega} |\nabla n_i^S|^2 \, dx  + \frac{|\upsilon_i|e^{2L}}{T_0^2} \, \|\nabla \varphi\|_{L^2}^2\Big)\, .
\end{align}
In order to satisfy \eqref{CONDI}, we recall the inequality \eqref{smallmeasure} and that $\ln n_i = \rho_i$.
% which imply that
%\begin{align*}
%	L^2\, \mu_\varphi(\Omega_L) \leq \|\rho_i^S\|^2_{L^2_{d\mu_\varphi}} + c_{PS} \,  \frac{d_1^\prime}{d_0^\prime}\,  (\|\nabla \sqrt{n_i^S}\|_{L^2}^2+|z_i|\, T_0^{-2} \, \|n_i^S\|_{L^{\infty}} \, \|\nabla\varphi^S\|_{L^2}^2) \, .
%\end{align*}
Hence, the condition \eqref{CONDI} is valid for $L = \max\{L_0,L_1\} $, where
\begin{align}\label{elleeins}
	L_1 = & \Big(4\frac{d_1^\prime|\upsilon_i|^2 \, C \, \|\nabla \varphi^S\|_{L^{p^S_N}}^{\frac{2}{\kappa}}}{d_0^\prime T_0^2}\Big)^{\frac{\kappa}{2(\kappa-1)}} \\ & \times \Big(\|\rho_i^S\|^2_{L^2_{d\mu_\varphi}} +
	C\,  [\mu_{\varphi}(\Omega)]^{1-\frac{1}{\kappa}} \, \|\nabla \varphi^S\|_{L^{p^S_N}}^{\frac{2}{\kappa}} \, e^{2\|\rho_i^S\|_{L^{\infty}}} \,   \frac{d_1^\prime}{d_0^\prime}\,  \Big(\|\nabla \rho_i^S\|_{L^2}^2+\frac{|\upsilon_i|}{T_0^{2}} \, \|\nabla\varphi\|_{L^2}^2\Big)\Big)\, .\nonumber
	%
	%\|\rho_i^S\|^2_{L^2_{d\mu_\varphi}} + c_{PS}\,  \frac{d_1^\prime}{d_0^\prime}\,  \Big(\|\nabla \sqrt{n_i^S}\|_{L^2}^2+\frac{|z_i|}{T_0^{2}} \, \|n_i^S\|_{L^{\infty}} \, \|\nabla\varphi\|_{L^2}^2\Big)\Big)^{\frac{1}{2}}  \,  .
\end{align}
%
%
%We let $\psi(x) = \exp(\rho_i(x)/2)$, and employing the definition \eqref{mob1} we get
%\begin{align*}
%	2 \, \int_{\Omega} \tilde{d_i}(u) \, |\nabla \psi|^2 \, dx \leq \int_{\Omega} \tilde{d_i}(u) \, \psi^2 \, (|\nabla \rho_i^S|^2 + |z_i|^2u^{-2} \, |\nabla \varphi|^2) \, dx \, .
%\end{align*}

{\bf Upper bound:} We test the equation \eqref{ions} with $\rho_{i,L} = \max\{\rho_i-L,0\}$ for $L \geq L_0 = \sup_{x\in S} \rho_i^S(x)$. We easily obtain for $\psi = \exp(\rho_i/2)$ that
\begin{align*}
4 \,  \int_{\Omega_L} \tilde{d}_i(u) \,|\nabla \psi|^2 \, dx = -2 \, \upsilon_i  \int_{\Omega_L} \frac{\tilde{d}_i(u)}{u}\, \psi \,\nabla \varphi \cdot \nabla \psi \, dx \, 
\end{align*}
with $\Omega_L := \{x \, : \, \rho_i(x) > L\}$, and this implies that
\begin{align*}
d_0^\prime \, \int_{\Omega_L} |\nabla \psi|^2 \, dx \leq \, \frac{|\upsilon_i|^2}{4T^2_0} \, d_1^\prime  \, \int_{\Omega_L} \psi^2 \, |\nabla \varphi|^2 \, dx \, . 
\end{align*}
With $L^\prime = \exp(L/2) $, this yields
\begin{align*}
 \int_{\Omega_{L}} |\nabla \psi|^2 \, dx \leq \, \frac{|\upsilon_i|^2 d_1^\prime}{2d_0^\prime T^2_0}  \, \int_{\Omega_{L}} ((\psi-L^\prime)^2 + (L^\prime)^2) \, |\nabla \varphi|^2 \, dx \, . 
\end{align*}
With $\kappa = 1 + 2\lambda_N/(n-2) > 1$, we have by the embedding Lemma \ref{lieberman}
\begin{align*}
	\int_{\Omega_{L}} (\psi-L^\prime)^2 \, |\nabla \varphi|^2 \, dx \leq& \left(\int_{\Omega}  |\psi_{L^\prime}|^{2\kappa} \, d\mu_{\varphi}\right)^{\frac{1}{\kappa}} \, [\mu_{\varphi}(\Omega_{L})]^{1-\frac{1}{\kappa}}\\
	\leq& C\, \|\varphi^S\|_{L^{p_N^S}}^{\frac{2}{\kappa}}  \, [\mu_{\varphi}(\Omega_{L})]^{1-\frac{1}{\kappa}}\, \int_{\Omega_{L}} |\nabla \psi|^2 \, dx \, 
\end{align*}
Hence, choosing $L$ so large that
\begin{align}\label{condiforsmall}
\frac{|\upsilon_i|^2d_1^\prime}{2d_0^{\prime}T_0^2} \, C\, \|\varphi^S\|_{L^{p_N^S}}^{\frac{2}{\kappa}} \, [\mu_{\varphi}(\Omega_{L^\prime})]^{1-\frac{1}{\kappa}} \leq \frac{1}{2} \, ,
\end{align}
we obtain that
\begin{align}\label{controlinfty}
\int_{\Omega_L} |\nabla \psi|^2 \, dx \leq \, \frac{|\upsilon_i|^2d^\prime_1}{d_0^\prime T_0^2}  \, (L^\prime)^2 \, \mu_{\varphi}(\Omega_{L}) \, . 
\end{align}
In order to satisfy \eqref{condiforsmall}, we use the same argument as for \eqref{CONDI}. In fact, we can choose $L > \max\{L_1/2^{\frac{\kappa}{\kappa-1}}, L_0\}$ with $L_1$ from \eqref{elleeins}.
% the inequality \eqref{rhomuintegrable} which implies that
%	\begin{align*}
%	L^2\, \mu_\varphi(\Omega_L) \leq \|\rho_i^S\|^2_{L^2_{\mu_\varphi}} + c_{PS} \,  \frac{d_1^\prime}{d_0^\prime}\,  (\|\nabla \sqrt{n_i^S}\|_{L^2}^2+|z_i|\, T_0^{-2} \, \|n_i^S\|_{L^{\infty}} \, \|\nabla\varphi^S\|_{L^2}^2) \, .
%\end{align*}
%implying that
%	\begin{align}\label{rhomuintegrable2}
%	4 \, (\ln L^\prime)^2\, \mu_\varphi(\Omega_{L^\prime}) \leq \|\rho_i^S\|^2_{L^2_{\mu_\varphi}} + c_{PS}^\prime \, ([\nabla \varphi]^M_{2,n-2+2\lambda_N})^2 \,  \frac{d_1^\prime}{d_0^\prime}\,  (\|\nabla \sqrt{n_i^S}\|_{L^2}^2+|z_i|\, T_0^{-2} \, \|n_i^S\|_{L^{\infty}} \, \|\nabla\varphi\|_{L^2}^2) \, .
%\end{align}
%Hence, the condition \eqref{condiforsmall} is valid for $L > L_1$, where
%\begin{align*}
%	L_1 = \Big(\frac{d_1^\prime|z_i|c_{PS}}{d_0^\prime T_0}\Big)^{\frac{\kappa}{2(\kappa-1)}} \, \Big(\|\rho_i^S\|^2_{L^2_{\mu_\varphi}} + c_{PS}\,  \frac{d_1^\prime}{d_0^\prime}\,  \Big(\|\nabla \sqrt{n_i^S}\|_{L^2}^2+\frac{|z_i|}{T_0^{2}} \, \|n_i^S\|_{L^{\infty}} \, \|\nabla\varphi\|_{L^2}^2\Big)\Big)^{\frac{1}{2}}  \,  .
%\end{align*}
Now, from \eqref{controlinfty}, we infer after a few steps (cf.\ \eqref{laberling2} -- \eqref{rhoibounded})
\begin{align*}
\int_{A_{L^{\prime}}} (\psi-L^{\prime}) \, d\mu \leq C^{\frac{1}{2}} \, \|\nabla \varphi^S\|_{L^{p_N^S}}^{\frac{1}{\kappa}} \, \Big(\frac{|\upsilon_i|^2d^\prime_1}{d_0^\prime T_0^2}\Big)^{\frac{1}{2}} \,  L^{\prime}\,  [\mu_{\varphi}(A_{L^{\prime}})]^{1+\frac{1}{2}-\frac{1}{2\kappa}}, \quad \kappa = 1 + \frac{2\lambda_N}{n-2} \, ,
\end{align*}
and $A_{L^{\prime}} = \{x \, : \, \psi(x) > L^{\prime}\}$ for $L^{\prime} > \exp(\max\{L_0,L_1\}/2)$. The upper bound follows from the Lemma 5.1 of Chapter II in \cite{MR0244627}. More precisely, with $L_0^{\prime} = \max\{e^{L_0}, e^{L_1}\}$ and $\sigma = (1-1/\kappa)/2$ we get
\begin{align}\label{psiabove}
	\psi(x) \leq 2^{1+\frac{1}{\sigma}} \, (L_0^{\prime} + [C^{\frac{1}{2}} \, \|\nabla \varphi^S\|_{L^{p_N^S}}^{\frac{1}{\kappa}} \, \Big(\frac{|\upsilon_i|^2d^\prime_1}{d_0^\prime T_0^2}\Big)^{\frac{1}{2}}]^{\frac{1}{\sigma}} \, \|\psi\|_{L^1}) \, . 
\end{align}
Note that $\|\psi\|_{L^1} \leq \|n\|_{L^1}^{\frac{1}{2}}$ which we control by means of \eqref{gradl2deux}.
\\

{\bf Lower bound:} We choose in the weak form of \eqref{ions} a testfunction of the form $e^{-\rho_i} \, \eta$, $\eta \in C^1_c(\Omega \setminus \overline{S})$. After a few calculation we obtain that
\begin{align*}
\int_\Omega \tilde{d}_i(u) \, (\nabla \rho_i + \frac{\upsilon_i}{u} \nabla \varphi) \cdot \nabla \eta \, dx = \int_{\Omega} \tilde{d}_i(u) \, \eta \,  (\nabla \rho_i + \frac{\upsilon_i}{u} \nabla \varphi) \cdot \nabla \rho_i \, dx\, .
\end{align*}	
After an approximation argument, we insert $\eta = \rho_{i,L} = \min\{\rho_i-L,0\}$ with $L\leq L_0 = \inf_{x\in S} \rho_i^S(x)$. Observe also that the right-hand is less than $|z_i|^2 \, \int_{\Omega} (\tilde{d}_i(u)/2u^2) \, |\eta| \,  |\nabla \varphi|^2 \, dx$ for $\eta$ nonpositive. Thus, it is readily shown again that
\begin{align*}
	\int_{\Omega_L} \tilde{d}_i(u) \, |\nabla \rho_i|^2 \, dx \leq \int_{\Omega_L}  \frac{|\upsilon_i|^2\tilde{d}_i(u) }{u^2} \, |\nabla \varphi|^2 \, (1+|\rho_{i,L}|) \, dx 
	\end{align*}
where $\Omega_{L} = \{x\, :\, \rho_i(x) < L\}$. The lower bound follows from similar arguments as above. 
\end{proof}
\begin{rem}\label{STRUC}
We can also state the $L^{\infty}$ bound in the following form, which later will be of use: There is a continuous function $\Phi = \Phi(a_1,a_2,a_3)$, nondecreasing in its arguments such that
\begin{align*}
	\|\rho_i\|_{L^{\infty}} \leq  \Phi(\|\rho^S_i\|_{L^{\infty}} + \|\nabla \rho_i^S\|_{L^2}, \, T_0^{-1}, \, \|\nabla \varphi_S\|_{L^{p_N^S}}) \, .
\end{align*}	
and with the additional property that $\Phi(a_1, a_2, 0) = \Phi_0(a_1)$ is independent of the second argument and $\liminf_{a_2 \rightarrow +\infty} \Phi(a_1,a_2,a_3) = +\infty$ if $a_1,a_3 >0$.
\end{rem}
We prove this remark by checking carefully throughout the proof of the estimates (see \eqref{rhomuintegrable}, \eqref{gradl2deux} and \eqref{elleeins}, and \eqref{psiabove}) that each occurrence of $T_0^{-1}$ is multiplied by a positive power of $\|\nabla \varphi^S\|_{L^{p^S_N}}$.\\

Thus, in view of Lemma \ref{boundvalues} the lower and upper bounds for the temperature are the decisive question to obtain the property \eqref{zusatz2} in Lemma \ref{NPPregu} and the regularity of weak solutions. We will discuss this point under the standard assumption that the heat conduction coefficient is a uniformly positive bounded function
\begin{align}\label{kappaass}
	0 < \kappa_0 \leq \kappa(\rho,u) \leq \kappa_1 \quad \text{ for all } \quad (\rho, \, u) \in \mathbb{R}^{N} \times \mathbb{R}_+ \, . 
\end{align}
\begin{lemma}
	We adopt the assumptions of Lemma \ref{boundvalues} and \eqref{kappaass}. Then, for every weak solution $(\rho, \, u)$ such that $T_0 >0$, we find that $\sup_{x\in \Omega} u(x) < +\infty$. In particular, every weak solution to \eqref{ions}, \eqref{poisson}, \eqref{templie} with mixed boundary conditions \eqref{BC}, \eqref{BC2} is regular in the sense of Lemma \ref{NPPregu} provided that $T_0$ is strictly positive. 
\end{lemma}
\begin{proof}
Due to the Lemma \ref{boundvalues}, for $i = 1,\ldots, N-1$ we find
\begin{gather*}
	0 < m_0^{\prime} \leq m_i(\rho(x),\, u(x)) \leq m_1^\prime  \quad \text{ in } \Omega \, ,\\
	m_0^\prime :=d_0^\prime \, \exp(\inf_{x\in \Omega} \rho_i(x)), \quad  m_1^\prime :=d_1^\prime \, \exp(\sup_{x\in \Omega} \rho_i(x)) \, ,
\end{gather*}
Then, we can apply the Lemma \ref{contihold} an obtain that $\rho_i \in C^{\lambda_i}(\overline{\Omega})$ with $0 <\lambda_i\leq \min\{1-n/p_i^S, \, \lambda_N\}$ depending on $m_0^{\prime}$ and $m_1^{\prime}$. Moreover, there is $c$ depending on the same quantities such that for $i = 1, \ldots, N-1$
\begin{align}\label{iondiffus}
	\|\rho_i\|_{C^{\lambda_i}(\overline{\Omega})} + \|\nabla \rho_i\|_{L^{2,n-2+2\lambda_i}_M(\Omega)} \leq c \, (\|\nabla \rho_i^S\|_{L^{p_i^S}(\Omega)} + \|\nabla \varphi^S\|_{L^{p_N^S}(\Omega)} )\, . 
\end{align}
Thus, the measure $d\calm = \Pi \, dx$ is diffusive. For $\Omega_{\rho} = \Omega\cap B_{\rho}$ with $B_{\rho}$ arbitrary we have
\begin{align}\label{Pidissipative}
	\int_{\Omega_\rho} \, |\Pi| \, dx \leq \int_{\Omega_{\rho}} (m(\rho,u) \nabla \rho \cdot \nabla \rho)^{\frac{1}{2}} \,  (m(\rho,u) \upsilon \cdot \upsilon)^{\frac{1}{2}} \, |\nabla \rho_N| + \frac{1}{u} \, (m(\rho,u) \upsilon \cdot \upsilon) \, |\nabla \rho_N|^2\, dx\nonumber\\
	\leq m_1^{\prime} \, |\upsilon| \,  \sum_{i=1}^{N-1} \, \Big(\int_{B_{\rho}} |\nabla \rho_i|^2 \, dx\Big)^{\frac{1}{2}} \, \Big(\int_{B_{\rho}} |\nabla \rho_N|^2 \, dx\Big)^{\frac{1}{2}} + \frac{m_1^\prime \, |\upsilon|^2}{T_0} \,  \int_{B_{\rho}} |\nabla \rho_N|^2 \, dx\nonumber\\
		\leq m_1^{\prime} \, |\upsilon| \,  \sum_{i=1}^{N-1} \, [\nabla \rho_i]_{2,\alpha_i}  \, [\nabla \rho_N]_{2,\alpha_N} \, \rho^{\frac{\alpha_i+\alpha_N}{2}} + \frac{m_1^\prime \, |\upsilon|^2}{T_0} \,   [\nabla \rho_N]_{2,\alpha_N}^2 \, \rho^{\alpha_N} \, .
\end{align}
We combine with \eqref{iondiffus} and with Remark \ref{dissipphi} and we find that $\int_{\Omega_\rho} \, |\Pi| \, dx \leq c \, \rho^{n-2+2\min_{i=1,\ldots,N} \lambda_i}$. Hence we obtain the upper bound for $u$ by the methods of Corollary \ref{coroN+1}.
\end{proof}
Thus, the problem of regularity is reduced to obtaining a strict lower bound for the temperature. We cannot show this in general as the measure $\Pi$ is not necessarily positive. However, we can prove that regular solutions always exist under suitable compatibility conditions for the data. To prove this result we note that the previous results, in particular the Remark \ref{STRUC} and the estimate \eqref{Pidissipative}, allow us with $\lambda_0 = \min_i \lambda_i$ to state that
\begin{align}\label{pistruc}
[\Pi]_{1,\lambda_0}^M \leq \|\nabla \varphi_S\|_{L^{p_N^S}}\, 	\Phi_1(\|\rho^S_i\|_{L^{\infty}} + \|\nabla \rho_i^S\|_{L^2}, \,  T_0^{-1}, \, \|\nabla \varphi_S\|_{L^{p_N^S}}) \, ,
\end{align}
where $\Phi_1$ is a function with the same structure as in Remark \ref{STRUC}.
\begin{lemma}
Let $\Phi_1$ be the function of \eqref{pistruc}. For $N_0, \, \varphi_0 >0$ and a constant $c_0$ given hereafter define $$T_0^S = T_0^S(N_0,\varphi_0) = \inf\{t > 0 \, :\,  t + c_0\, \varphi_0\,\Phi_1(N_0,1/t,\varphi_0)\} \, .$$ Then, there is $t_0 >0 $ such that $T_0^S = t_0 + c_0 \, \varphi_0\,  \Phi_1(N_0,1/t_0,\varphi_0)$.
Moreover, suppose that $\inf_{x\in S} u^S(x) \geq T_0^S$, that $\|\nabla \varphi^S\|_{L^{p_N^S}} \leq \varphi_0$ and that $\sum_{i=1}^{N-1} \|\nabla \rho^S_i\|_{L^{p_i^S}} \leq N_0$. Then there exists a regular weak solution satisfying $\inf_{x\in \Omega} u(x) \geq t_0$.
\end{lemma}
\begin{proof}
As a preliminary, we remark that the function $t \mapsto  t + c_0\, \varphi_0\, \Phi_1(N_0,1/t,\varphi_0)$ tends to $+\infty$ for $t \rightarrow\{0, \, +\infty\}$. Since it is continuous, it admits at least one global minimum in some $t_0 \in ]0,+\infty[$ depending on $(N_0,\varphi_0)$.\\

To prove the existence claim, we apply a fixed point argument. We let $\hat{\rho}_1,\ldots,\hat{\rho}_{N}, \, \hat{u} \in L^2(\Omega)$ be any functions such that $\hat{u} \geq t_0$ almost everywhere in $\Omega$. We first let $\varphi = \rho_N$ be the unique solution to \eqref{poisson} with $\varepsilon$ frozen at $(\hat{\rho},\hat{u})$. Then for $i=1,\ldots,N-1$ we let $\rho_i$ be the unique weak solutions to \eqref{ions} with data $\hat{u}$ and $\varphi$, and coefficients frozen at $(\hat{\rho},\hat{u})$. Finally, we let $u$ be the unique solution to 
\begin{align}\label{iter}
-\divv(\kappa(\rho,\hat{u}) \, \nabla u) = 	\sum_{i=1}^{N-1} \upsilon_i \, m_i(\rho,\, \hat{u}) \, \Big(\nabla \rho_i + \frac{\upsilon_i}{\hat{u}} \, \nabla \rho_N\Big) \cdot \nabla \rho_N =: \Pi\, .
\end{align}
If we choose data $\|\rho_i^S\|_{W^{1,p_i^S}(\Omega)} \leq N_0$ and $\|\varphi^S\|_{W^{1,p_N^S}(\Omega)} \leq \varphi_0$, then in view of \eqref{pistruc}, we can verify that $\|\Pi\|_{L^{1,n-2+2\lambda_0}_M(\Omega)} \leq \varphi_0\,\Phi_1(N_0,t_0^{-1},\varphi_0)$. Using \eqref{iter} and the maximum principle arguments exposed several times in the paper, we get 
\begin{align*}
	\inf_{x\in\Omega} u(x) \geq \inf_{x\in S} u^S(x) - c(\lambda_0, \,\kappa_0,\kappa_1)  \, \|\Pi\|_{L^{1,n-2+2\lambda_0}_M}\, , 
\end{align*}
which with $c_0 = c(\lambda_0, \,\kappa_0,\kappa_1)$ we can also state in the form 
\begin{align*}
	\inf_{x\in\Omega} u(x) \geq \inf_{x\in S} u^S(x) -c_0 \,\varphi_0\,  \Phi_1(N_0,t_0^{-1},\varphi_0) \geq T_0^S -c_0\,  \varphi_0\, \Phi_1(N_0,t_0^{-1},\varphi_0) = t_0 \, .
\end{align*} 
Now it is not difficult to show that for $R_0$ large enough the map $(\hat{\rho}, \, \hat{u}) \mapsto (\rho, \, u)$ is a self map for the set
\begin{align*}
	\{(\rho_1,\ldots,\rho_N,u) \in L^2(\Omega; \, \mathbb{R}^{N+1}) \, : \, \|(\rho,u)\|_{L^2} \leq R_0, \, \inf_{x\in\Omega} u(x) \geq t_0\} \, . 
\end{align*}
The compactness of the map follows from the properties of Sobolev-embeddings, so that the Schauder fixed-point theorem can be employed to prove the existence claim. 
\end{proof}
To finish the proof of Theorem \ref{NPPtheo}, \eqref{ExiNPP} it remains to observe that $T_0^S(N_0,\varphi_0)$ is increasing in $\varphi_0$, tends to zero for $\varphi_0 \rightarrow 0$ and to $+\infty$ for $\varphi_0 \rightarrow + \infty$. We define $\hat{\varphi}_0 (N_0,T_0)$ via the implicit solution of $T_0^S(N_0,\varphi_0) = T_0$, and are done.

%\bibliography{Refs}

\begin{thebibliography}{10}
	
	\bibitem{MR0287301}
	D.~R. Adams.
	\newblock Traces of potentials arising from translation invariant operators.
	\newblock {\em Ann. Scuola Norm. Sup. Pisa Cl. Sci. (3)}, 25:203--217, 1971.
	
	\bibitem{MR1888855}
	G.~Albinus, H.~Gajewski, and R.~H\"{u}nlich.
	\newblock Thermodynamic design of energy models of semiconductor devices.
	\newblock {\em Nonlinearity}, 15(2):367--383, 2002.
	\newblock \href {https://doi.org/10.1088/0951-7715/15/2/307}
	{\path{doi:10.1088/0951-7715/15/2/307}}.
	
	\bibitem{MR3397331}
	L.~Beck, M.~Bul\'{\i}\v{c}ek, and J.~Frehse.
	\newblock Old and new results in regularity theory for diagonal elliptic
	systems via blowup techniques.
	\newblock {\em J. Differential Equations}, 259(11):6528--6572, 2015.
	\newblock \href {https://doi.org/10.1016/j.jde.2015.07.030}
	{\path{doi:10.1016/j.jde.2015.07.030}}.
	
	\bibitem{MR1354907}
	P.~B\'{e}nilan, L.~Boccardo, T.~Gallou\"{e}t, R.~Gariepy, M.~Pierre, and J.~L.
	V\'{a}zquez.
	\newblock An {$L^1$}-theory of existence and uniqueness of solutions of
	nonlinear elliptic equations.
	\newblock {\em Ann. Scuola Norm. Sup. Pisa Cl. Sci. (4)}, 22(2):241--273, 1995.
	\newblock URL: \url{http://www.numdam.org/item?id=ASNSP_1995_4_22_2_241_0}.
	
	\bibitem{MR1944757}
	A.~Bensoussan and L.~Boccardo.
	\newblock Nonlinear systems of elliptic equations with natural growth
	conditions and sign conditions.
	\newblock volume~46, pages 143--166. 2002.
	\newblock Special issue dedicated to the memory of Jacques-Louis Lions.
	\newblock \href {https://doi.org/10.1007/s00245-002-0753-3}
	{\path{doi:10.1007/s00245-002-0753-3}}.
	
	\bibitem{MR1917320}
	A.~Bensoussan and J.~Frehse.
	\newblock {\em Regularity results for nonlinear elliptic systems and
		applications}, volume 151 of {\em Applied Mathematical Sciences}.
	\newblock Springer-Verlag, Berlin, 2002.
	\newblock \href {https://doi.org/10.1007/978-3-662-12905-0}
	{\path{doi:10.1007/978-3-662-12905-0}}.
	
	\bibitem{MR1452883}
	L.~Boccardo.
	\newblock Homogenization and continuous dependence for {D}irichlet problems in
	{$L^1$}.
	\newblock In {\em Partial differential equation methods in control and shape
		analysis ({P}isa)}, volume 188 of {\em Lecture Notes in Pure and Appl.
		Math.}, pages 41--52. Dekker, New York, 1997.
	
	\bibitem{MR1025884}
	L.~Boccardo and T.~Gallou\"{e}t.
	\newblock Nonlinear elliptic and parabolic equations involving measure data.
	\newblock {\em J. Funct. Anal.}, 87(1):149--169, 1989.
	\newblock \href {https://doi.org/10.1016/0022-1236(89)90005-0}
	{\path{doi:10.1016/0022-1236(89)90005-0}}.
	
	\bibitem{MR1409661}
	L.~Boccardo, T.~Gallou\"{e}t, and L.~Orsina.
	\newblock Existence and uniqueness of entropy solutions for nonlinear elliptic
	equations with measure data.
	\newblock {\em Ann. Inst. H. Poincar\'{e} C Anal. Non Lin\'{e}aire},
	13(5):539--551, 1996.
	\newblock \href {https://doi.org/10.1016/S0294-1449(16)30113-5}
	{\path{doi:10.1016/S0294-1449(16)30113-5}}.
	
	\bibitem{MR4348166}
	L.~Boccardo and L.~Orsina.
	\newblock An elliptic system related to the stationary thermistor problem.
	\newblock {\em SIAM J. Math. Anal.}, 53(6):6910--6931, 2021.
	\newblock \href {https://doi.org/10.1137/21M1420058}
	{\path{doi:10.1137/21M1420058}}.
	
	\bibitem{MR4350569}
	D.~Bothe and P.-E. Druet.
	\newblock Well-posedness analysis of multicomponent incompressible flow models.
	\newblock {\em J. Evol. Equ.}, 21(4):4039--4093, 2021.
	\newblock \href {https://doi.org/10.1007/s00028-021-00712-3}
	{\path{doi:10.1007/s00028-021-00712-3}}.
	
	\bibitem{MR3160531}
	D.~Bothe, A.~Fischer, M.~Pierre, and G.~Rolland.
	\newblock Global existence for diffusion-electromigration systems in space
	dimension three and higher.
	\newblock {\em Nonlinear Anal.}, 99:152--166, 2014.
	\newblock \href {https://doi.org/10.1016/j.na.2013.12.015}
	{\path{doi:10.1016/j.na.2013.12.015}}.
	
	\bibitem{MR2332419}
	H.~Brezis and J.~Van~Schaftingen.
	\newblock Boundary estimates for elliptic systems with {$L^1$}-data.
	\newblock {\em Calc. Var. Partial Differential Equations}, 30(3):369--388,
	2007.
	\newblock \href {https://doi.org/10.1007/s00526-007-0094-9}
	{\path{doi:10.1007/s00526-007-0094-9}}.
	
	\bibitem{MR0987900}
	G.~Cimatti.
	\newblock Remark on existence and uniqueness for the thermistor problem under
	mixed boundary conditions.
	\newblock {\em Quart. Appl. Math.}, 47(1):117--121, 1989.
	\newblock \href {https://doi.org/10.1090/qam/987900}
	{\path{doi:10.1090/qam/987900}}.
	
	\bibitem{MR3158841}
	G.~R. Cirmi and S.~Leonardi.
	\newblock Higher differentiability for the solutions of nonlinear elliptic
	systems with lower-order terms and {$L^{1,\theta}$}-data.
	\newblock {\em Ann. Mat. Pura Appl. (4)}, 193(1):115--131, 2014.
	\newblock \href {https://doi.org/10.1007/s10231-012-0269-7}
	{\path{doi:10.1007/s10231-012-0269-7}}.
	
	\bibitem{MR1446209}
	S.~Clain and R.~Touzani.
	\newblock A two-dimensional stationary induction heating problem.
	\newblock {\em Math. Methods Appl. Sci.}, 20(9):759--766, 1997.
	\newblock \href
	{https://doi.org/10.1002/(SICI)1099-1476(199706)20:9<759::AID-MMA879>3.3.CO;2-J}
	{\path{doi:10.1002/(SICI)1099-1476(199706)20:9<759::AID-MMA879>3.3.CO;2-J}}.
	
	\bibitem{MR1225511}
	R.~Coifman, P.-L. Lions, Y.~Meyer, and S.~Semmes.
	\newblock Compensated compactness and {H}ardy spaces.
	\newblock {\em J. Math. Pures Appl. (9)}, 72(3):247--286, 1993.
	
	\bibitem{MR1459011}
	L.~Consiglieri.
	\newblock Stationary weak solutions for a class of non-{N}ewtonian fluids with
	energy transfer.
	\newblock {\em Internat. J. Non-Linear Mech.}, 32(5):961--972, 1997.
	\newblock \href {https://doi.org/10.1016/S0020-7462(96)00087-X}
	{\path{doi:10.1016/S0020-7462(96)00087-X}}.
	
	\bibitem{MR1760541}
	G.~Dal~Maso, F.~Murat, L.~Orsina, and A.~Prignet.
	\newblock Renormalized solutions of elliptic equations with general measure
	data.
	\newblock {\em Ann. Scuola Norm. Sup. Pisa Cl. Sci. (4)}, 28(4):741--808, 1999.
	\newblock URL: \url{http://www.numdam.org/item?id=ASNSP_1999_4_28_4_741_0}.
	
	\bibitem{MR1147281}
	M.~Dauge.
	\newblock Neumann and mixed problems on curvilinear polyhedra.
	\newblock {\em Integral Equations Operator Theory}, 15(2):227--261, 1992.
	\newblock \href {https://doi.org/10.1007/BF01204238}
	{\path{doi:10.1007/BF01204238}}.
	
	\bibitem{MR1484710}
	G.~Dolzmann, N.~Hungerb\"{u}hler, and S.~M\"{u}ller.
	\newblock Non-linear elliptic systems with measure-valued right hand side.
	\newblock {\em Math. Z.}, 226(4):545--574, 1997.
	\newblock \href {https://doi.org/10.1007/PL00004354}
	{\path{doi:10.1007/PL00004354}}.
	
	\bibitem{dreyer2013overcoming}
	W.~Dreyer, C.~Guhlke, and R.~M{\"u}ller.
	\newblock Overcoming the shortcomings of the {N}ernst--{P}lanck model.
	\newblock {\em Physical Chemistry Chemical Physics}, 15(19):7075--7086, 2013.
	\newblock \href {https://doi.org/10.1039/C3CP44390F}
	{\path{doi:10.1039/C3CP44390F}}.
	
	\bibitem{MR2545657}
	P.-E. Druet.
	\newblock Existence for the stationary {MHD}-equations coupled to heat transfer
	with nonlocal radiation effects.
	\newblock {\em Czechoslovak Math. J.}, 59(134)(3):791--825, 2009.
	\newblock \href {https://doi.org/10.1007/s10587-009-0048-9}
	{\path{doi:10.1007/s10587-009-0048-9}}.
	
	\bibitem{MR2523256}
	P.-E. Druet.
	\newblock Existence of weak solutions to the time-dependent {MHD} equations
	coupled to the heat equation with nonlocal radiation boundary conditions.
	\newblock {\em Nonlinear Anal. Real World Appl.}, 10(5):2914--2936, 2009.
	\newblock \href {https://doi.org/10.1016/j.nonrwa.2008.09.015}
	{\path{doi:10.1016/j.nonrwa.2008.09.015}}.
	
	\bibitem{MR4452316}
	P.-E. Druet.
	\newblock Maximal mixed parabolic-hyperbolic regularity for the full equations
	of multicomponent fluid dynamics.
	\newblock {\em Nonlinearity}, 35(7):3812--3882, 2022.
	\newblock \href {https://doi.org/10.1088/1361-6544/ac5679}
	{\path{doi:10.1088/1361-6544/ac5679}}.
	
	\bibitem{MR2837572}
	P.-E. Druet, O.~Klein, J.~Sprekels, F.~Tr\"{o}ltzsch, and I.~Yousept.
	\newblock Optimal control of three-dimensional state-constrained induction
	heating problems with nonlocal radiation effects.
	\newblock {\em SIAM J. Control Optim.}, 49(4):1707--1736, 2011.
	\newblock \href {https://doi.org/10.1137/090760544}
	{\path{doi:10.1137/090760544}}.
	
	\bibitem{MR2506064}
	P.-E. Druet and J.~Naumann.
	\newblock On the existence of weak solutions to a stationary one-equation
	{RANS} model with unbounded eddy viscosities.
	\newblock {\em Ann. Univ. Ferrara Sez. VII Sci. Mat.}, 55(1):67--87, 2009.
	\newblock \href {https://doi.org/10.1007/s11565-009-0062-8}
	{\path{doi:10.1007/s11565-009-0062-8}}.
	
	\bibitem{drunau}
	P.-E. Druet, J.~Naumann, and J.~Wolf.
	\newblock A {M}eyers' type estimate for weak solutions to a generalized
	stationary navier-stokes system.
	\newblock Preprint 06 of the Institute for Mathematics of the
	Humboldt-University, Berlin, 2008.
	\newblock Available in pdf-format at
	{\verb|http://www.mathematik.hu-berlin.de/publ/pre/2008/p-list-08|}.
	
	\bibitem{MR3729430}
	E.~Feireisl and A.~Novotn\'{y}.
	\newblock {\em Singular limits in thermodynamics of viscous fluids}.
	\newblock Advances in Mathematical Fluid Mechanics. Birkh\"{a}user/Springer,
	Cham, second edition, 2017.
	
	\bibitem{MR0717034}
	M.~Giaquinta.
	\newblock {\em Multiple integrals in the calculus of variations and nonlinear
		elliptic systems}, volume 105 of {\em Annals of Mathematics Studies}.
	\newblock Princeton University Press, Princeton, NJ, 1983.
	
	\bibitem{MR1827089}
	J.~A. Griepentrog and L.~Recke.
	\newblock Linear elliptic boundary value problems with non-smooth data: normal
	solvability on {S}obolev-{C}ampanato spaces.
	\newblock {\em Math. Nachr.}, 225:39--74, 2001.
	\newblock \href
	{https://doi.org/10.1002/1522-2616(200105)225:1<39::AID-MANA39>3.3.CO;2-X}
	{\path{doi:10.1002/1522-2616(200105)225:1<39::AID-MANA39>3.3.CO;2-X}}.
	
	\bibitem{MR0990595}
	K.~Gr\"{o}ger.
	\newblock A {$W^{1,p}$}-estimate for solutions to mixed boundary value problems
	for second order elliptic differential equations.
	\newblock {\em Math. Ann.}, 283(4):679--687, 1989.
	\newblock \href {https://doi.org/10.1007/BF01442860}
	{\path{doi:10.1007/BF01442860}}.
	
	\bibitem{MR2271364}
	K.~Gr\"{o}ger and L.~Recke.
	\newblock Applications of differential calculus to quasilinear elliptic
	boundary value problems with non-smooth data.
	\newblock {\em NoDEA Nonlinear Differential Equations Appl.}, 13(3):263--285,
	2006.
	\newblock \href {https://doi.org/10.1007/s00030-006-3017-0}
	{\path{doi:10.1007/s00030-006-3017-0}}.
	
	\bibitem{MR3573649}
	R.~Haller-Dintelmann, A.~Jonsson, D.~Knees, and J.~Rehberg.
	\newblock Elliptic and parabolic regularity for second-order divergence
	operators with mixed boundary conditions.
	\newblock {\em Math. Methods Appl. Sci.}, 39(17):5007--5026, 2016.
	\newblock \href {https://doi.org/10.1002/mma.3484}
	{\path{doi:10.1002/mma.3484}}.
	
	\bibitem{MR2378088}
	R.~Haller-Dintelmann, H.-C. Kaiser, and J.~Rehberg.
	\newblock Elliptic model problems including mixed boundary conditions and
	material heterogeneities.
	\newblock {\em J. Math. Pures Appl. (9)}, 89(1):25--48, 2008.
	\newblock \href {https://doi.org/10.1016/j.matpur.2007.09.001}
	{\path{doi:10.1016/j.matpur.2007.09.001}}.
	
	\bibitem{MR2599927}
	D.~H\"{o}mberg, C.~Meyer, J.~Rehberg, and W.~Ring.
	\newblock Optimal control for the thermistor problem.
	\newblock {\em SIAM J. Control Optim.}, 48(5):3449--3481, 2009/10.
	\newblock \href {https://doi.org/10.1137/080736259}
	{\path{doi:10.1137/080736259}}.
	
	\bibitem{MR4108365}
	C.-Y. Hsieh, T.-C. Lin, C.~Liu, and P.~Liu.
	\newblock Global existence of the non-isothermal
	{P}oisson-{N}ernst-{P}lanck-{F}ourier system.
	\newblock {\em J. Differential Equations}, 269(9):7287--7310, 2020.
	\newblock \href {https://doi.org/10.1016/j.jde.2020.05.037}
	{\path{doi:10.1016/j.jde.2020.05.037}}.
	
	\bibitem{MR1331981}
	D.~Jerison and C.~E. Kenig.
	\newblock The inhomogeneous {D}irichlet problem in {L}ipschitz domains.
	\newblock {\em J. Funct. Anal.}, 130(1):161--219, 1995.
	\newblock \href {https://doi.org/10.1006/jfan.1995.1067}
	{\path{doi:10.1006/jfan.1995.1067}}.
	
	\bibitem{MR2765805}
	A.~J\"{u}ngel.
	\newblock Energy transport in semiconductor devices.
	\newblock {\em Math. Comput. Model. Dyn. Syst.}, 16(1):1--22, 2010.
	\newblock \href {https://doi.org/10.1080/13873951003679017}
	{\path{doi:10.1080/13873951003679017}}.
	
	\bibitem{MR3092288}
	A.~J\"{u}ngel, R.~Pinnau, and E.~R\"{o}hrig.
	\newblock Existence analysis for a simplified transient energy-transport model
	for semiconductors.
	\newblock {\em Math. Methods Appl. Sci.}, 36(13):1701--1712, 2013.
	\newblock \href {https://doi.org/10.1002/mma.2715}
	{\path{doi:10.1002/mma.2715}}.
	
	\bibitem{MR1290667}
	T.~Kilpel\"{a}inen.
	\newblock H\"{o}lder continuity of solutions to quasilinear elliptic equations
	involving measures.
	\newblock {\em Potential Anal.}, 3(3):265--272, 1994.
	\newblock \href {https://doi.org/10.1007/BF01468246}
	{\path{doi:10.1007/BF01468246}}.
	
	\bibitem{MR1887015}
	T.~Kilpel\"{a}inen and X.~Zhong.
	\newblock Removable sets for continuous solutions of quasilinear elliptic
	equations.
	\newblock {\em Proc. Amer. Math. Soc.}, 130(6):1681--1688, 2002.
	\newblock \href {https://doi.org/10.1090/S0002-9939-01-06237-2}
	{\path{doi:10.1090/S0002-9939-01-06237-2}}.
	
	\bibitem{MR0244627}
	O.~A. Ladyzhenskaya and N.~N. Ural'tseva.
	\newblock {\em Linear and quasilinear elliptic equations}.
	\newblock Academic Press, New York-London, 1968.
	\newblock Translated from the Russian by Scripta Technica, Inc, Translation
	editor: Leon Ehrenpreis.
	
	\bibitem{MR2719762}
	S.~Leonardi and J.~Star\'{a}.
	\newblock Regularity results for the gradient of solutions of linear elliptic
	systems with {VMO}-coefficients and {$L^{1,\lambda}$} data.
	\newblock {\em Forum Math.}, 22(5):913--940, 2010.
	\newblock \href {https://doi.org/10.1515/FORUM.2010.048}
	{\path{doi:10.1515/FORUM.2010.048}}.
	
	\bibitem{MR1418142}
	R.~Lewandowski.
	\newblock The mathematical analysis of the coupling of a turbulent kinetic
	energy equation to the {N}avier-{S}tokes equation with an eddy viscosity.
	\newblock {\em Nonlinear Anal.}, 28(2):393--417, 1997.
	\newblock \href {https://doi.org/10.1016/0362-546X(95)00149-P}
	{\path{doi:10.1016/0362-546X(95)00149-P}}.
	
	\bibitem{MR1022556}
	G.~M. Lieberman.
	\newblock Optimal {H}\"{o}lder regularity for mixed boundary value problems.
	\newblock {\em J. Math. Anal. Appl.}, 143(2):572--586, 1989.
	\newblock \href {https://doi.org/10.1016/0022-247X(89)90061-9}
	{\path{doi:10.1016/0022-247X(89)90061-9}}.
	
	\bibitem{MR1233190}
	G.~M. Lieberman.
	\newblock Sharp forms of estimates for subsolutions and supersolutions of
	quasilinear elliptic equations involving measures.
	\newblock {\em Comm. Partial Differential Equations}, 18(7-8):1191--1212, 1993.
	\newblock \href {https://doi.org/10.1080/03605309308820969}
	{\path{doi:10.1080/03605309308820969}}.
	
	\bibitem{MR3702856}
	H.~Meinlschmidt, C.~Meyer, and J.~Rehberg.
	\newblock Optimal control of the thermistor problem in three spatial
	dimensions, {P}art 1: {E}xistence of optimal solutions.
	\newblock {\em SIAM J. Control Optim.}, 55(5):2876--2904, 2017.
	\newblock \href {https://doi.org/10.1137/16M1072644}
	{\path{doi:10.1137/16M1072644}}.
	
	\bibitem{MR0159110}
	N.~G. Meyers.
	\newblock An {$L^{p}$}-estimate for the gradient of solutions of second order
	elliptic divergence equations.
	\newblock {\em Ann. Scuola Norm. Sup. Pisa Cl. Sci. (3)}, 17:189--206, 1963.
	\newblock URL: \url{http://www.numdam.org/item/ASNSP_1963_3_17_3_189_0.pdf}.
	
	\bibitem{MR2352517}
	G.~Mingione.
	\newblock The {C}alder\'{o}n-{Z}ygmund theory for elliptic problems with
	measure data.
	\newblock {\em Ann. Sc. Norm. Super. Pisa Cl. Sci. (5)}, 6(2):195--261, 2007.
	\newblock URL: \url{http://www.numdam.org/item/ASNSP_2007_5_6_2_195_0}.
	
	\bibitem{MR1296252}
	B.~Mohammadi and O.~Pironneau.
	\newblock {\em Analysis of the {$k$}-epsilon turbulence model}.
	\newblock RAM: Research in Applied Mathematics. Masson, Paris; John Wiley \&
	Sons, Ltd., Chichester, 1994.
	
	\bibitem{MR2491627}
	P.~B. Mucha and M.~Pokorn\'{y}.
	\newblock On the steady compressible {N}avier-{S}tokes-{F}ourier system.
	\newblock {\em Comm. Math. Phys.}, 288(1):349--377, 2009.
	\newblock \href {https://doi.org/10.1007/s00220-009-0772-x}
	{\path{doi:10.1007/s00220-009-0772-x}}.
	
	\bibitem{MR3916818}
	P.~B.~a. Mucha, M.~Pokorn\'{y}, and E.~Zatorska.
	\newblock Existence of stationary weak solutions for compressible heat
	conducting flows.
	\newblock In {\em Handbook of mathematical analysis in mechanics of viscous
		fluids}, pages 2595--2662. Springer, Cham, 2018.
	\newblock \href {https://doi.org/10.1007/978-3-319-13344-7\_64}
	{\path{doi:10.1007/978-3-319-13344-7\_64}}.
	
	\bibitem{MR2185227}
	J.~Naumann.
	\newblock Existence of weak solutions to the equations of stationary motion of
	heat-conducting incompressible viscous fluids.
	\newblock In {\em Nonlinear elliptic and parabolic problems}, volume~64 of {\em
		Progr. Nonlinear Differential Equations Appl.}, pages 373--390.
	Birkh\"{a}user, Basel, 2005.
	\newblock \href {https://doi.org/10.1007/3-7643-7385-7\_21}
	{\path{doi:10.1007/3-7643-7385-7\_21}}.
	
	\bibitem{MR1097910}
	J.-M. Rakotoson.
	\newblock Quasilinear elliptic problems with measures as data.
	\newblock {\em Differential Integral Equations}, 4(3):449--457, 1991.
	
	\bibitem{MR1201197}
	J.-M. Rakotoson.
	\newblock Resolution of the critical cases for problems with {$L^1$}-data.
	\newblock {\em Asymptotic Anal.}, 6(3):285--293, 1993.
	
	\bibitem{MR0998128}
	J.-M. Rakotoson and W.~P. Ziemer.
	\newblock Local behavior of solutions of quasilinear elliptic equations with
	general structure.
	\newblock {\em Trans. Amer. Math. Soc.}, 319(2):747--764, 1990.
	\newblock \href {https://doi.org/10.2307/2001263} {\path{doi:10.2307/2001263}}.
	
	\bibitem{MR2495071}
	T.~Roub\'{\i}\v{c}ek.
	\newblock Thermo-visco-elasticity at small strains with {$L^1$}-data.
	\newblock {\em Quart. Appl. Math.}, 67(1):47--71, 2009.
	\newblock \href {https://doi.org/10.1090/S0033-569X-09-01094-3}
	{\path{doi:10.1090/S0033-569X-09-01094-3}}.
	
	\bibitem{MR0126601}
	G.~Stampacchia.
	\newblock Problemi al contorno ellitici, con dati discontinui, dotati di
	soluzionie h\"{o}lderiane.
	\newblock {\em Ann. Mat. Pura Appl. (4)}, 51:1--37, 1960.
	\newblock \href {https://doi.org/10.1007/BF02410941}
	{\path{doi:10.1007/BF02410941}}.
	
	\bibitem{MR0251373}
	G.~Stampacchia.
	\newblock {\em \`Equations elliptiques du second ordre \`a coefficients
		discontinus}, volume No. 16 (\'Et\'e, 1965) of {\em S\'{e}minaire de
		Math\'{e}matiques Sup\'{e}rieures [Seminar on Higher Mathematics]}.
	\newblock Les Presses de l'Universit\'{e} de Montr\'{e}al, Montreal, QC, 1966.
	
	\bibitem{ter2014holder}
	A.~Ter~Elst and J.~Rehberg.
	\newblock H{\"o}lder estimates for second-order operators with mixed boundary
	conditions.
	\newblock Technical report, Berlin: Weierstra{\ss}-Institut f{\"u}r Angewandte
	Analysis und Stochastik, 2014.
	\newblock Preprint no 1950.
	
\end{thebibliography}
%\bibliographystyle{abbrvurl}
%Hence, we will at first restrict to proving the results under the assumption that bounds from below and above are already known for the species concentrations: We assume that there are $0< n_0 < n_1$ such that the weak solution satisfies
%\begin{align}\label{zusatz}
% n_0 \leq c_i(x) = e^{\rho_i(x)-\mu_i^0(u(x))}\leq n_1 \quad \text{ for all } \quad x \in \Omega \, . 
%\end{align} 
%Then we have also $0 < m^\prime_0 \leq m_i(\rho,u) \leq m^\prime_1$ for appropriate numbers $m^\prime_0$ and $m^\prime_1$ which depends on the thresholds of the diffusivities and on the numbers $n_0$ and $n_1$ of \eqref{zusatz}.

% 
% \section{More complicated examples}
% 
% More Examples: not working
% \subsection{Time-harmonic Maxwell}
% 
% TIME-Harmonic Maxwell equations. Here we choose $N=3$ and skew-symmetric $a^i_j$ via
% \begin{align*}
%  a_{ij}(T, \, z) = r(T) \,  \sign(i-j) \, z^i_j
% \end{align*}
% This generates a $\curl r(T) \curl $ bilinear form. System is not diagonal ????
% 
% \subsection{Navier-Stokes: Turbulence model}

\appendix

\section{Auxiliary results}\label{append}

We begin with finishing the proof of Lemma \ref{ENERGYEST}. We start from the inequality \eqref{moule2} and show that the right-hand side can be controlled if we choose the exponents $r_i$ and $s_i$ as claimed in the statement.

We let first $\tilde{q}_i = q_i $ for $i=1,\ldots,N$ and $\tilde{q}_{N+1}=q_{N+1}/\alpha$. The assumptions of the lemma now guarantee that $\tilde{r}_i \leq \hat{r}_{\rho}(p_i,\tilde{q}_i,k_i)=:\chi_i$ (cf.\ \eqref{thresholds}) if $i \leq N$. For $i = N+1$, we have assumed that $r_{N+1} \leq \frac{1}{p_{N+1}^{\prime}} \hat{r}_{\rho}(p_{n+1},q_{N+1}p_{N+1}^{\prime},k_{N+1})$, which implies that
\begin{align}\label{thresholds2}
	r_{N+1} \leq \beta_1 \, \frac{p^*_{N+1}}{p_{N+1}^{\prime}} + \Big(\frac{1}{k^{\prime}_{N+1}}-\beta_1\Big) \, q_{N+1} \quad \text{ with a }\quad  0\leq \beta_1 < m:= \min\big\{\frac{1}{k^{'}_{N+1}}, \, \frac{p_{N+1}}{p^*_{N+1}}\big\}\, .
\end{align}
In particular, the definition of $\hat{r}_{\rho}$ implies that $\beta_1 = 0$ for $q_{N+1}p_{N+1}^{\prime} \geq p^*_{N+1}$. 

We let $\xi := \alpha(\delta) \, p_{N+1}^{\prime} = 1-\delta/(p_{N+1}-1)$ and if $q_{N+1} < p^*_{N+1}/p_{N+1}^\prime$, we restrict $\delta$ via
\begin{align*}
	0<	\delta < (p_{N+1}-1) \, \Big(\frac{\beta_1}{m} + \big(1-\frac{\beta_1}{m}\big) \, \frac{p_{N+1}^\prime q_{N+1}}{p^*_{N+1}}-1\Big) \quad \text{if} \quad q_{N+1} < p^*_{N+1}/p_{N+1}^\prime \, ,
\end{align*}	
to guarantee that $\xi \, p^*_{N+1} - p^\prime_{N+1}q_{N+1} >0$. With $\beta_1$ from \eqref{thresholds2} we can compute that
\begin{align*}
	\tilde{r}_{N+1} \leq  \beta \, p^*_{N+1} + \Big(\frac{1}{k^{\prime}_{N+1}}-\beta\Big) \, \tilde{q}_{N+1}\, , \quad
	\beta = \begin{cases}
		0 & \text{ if } q_{N+1}p_{N+1}^{\prime} \geq p^*_{N+1}\\
		\beta_1 \, \frac{p^*_{N+1} - p^\prime_{N+1}q_{N+1}}{\xi \, p^*_{N+1} - p^\prime_{N+1}q_{N+1}} & \text{ otherwise}
	\end{cases}\, .
\end{align*}
%Restricting $\delta$ via
%\begin{align*}
%0<	\delta < (p_{N+1}-1) \, \Big(\frac{\beta_1}{m} + \big(1-\frac{\beta_1}{m}\big) \, \frac{p_{N+1}^\prime q_{N+1}}{p^*_{N+1}}-1\Big) \quad \text{if} \quad q_{N+1} < p^*_{N+1}/p_{N+1}^\prime
%\end{align*}	
We easily verify that $0\leq \beta < m$, and we define $\chi_{N+1} := \beta \, p^*_{N+1} +(1/k^{\prime}_{N+1}-\beta) \, \tilde{q}_{N+1}$.\\
%with the $\beta$ of \eqref{thresholds}, and $\delta \leq \beta p_{N+1}^*/ p_{N+1}$ guaranteees that $\tilde{r}_{N+1} \leq \chi_{N+1}$ as desired.
Next we consider arbitrary $i\in \{1,\ldots,N+1\}$. If we have $\chi_i = \tilde{q}_i/k_i^{\prime}$ then with the help of H\"older's inequality
\begin{align*}
	\int_{\Omega} \psi_i(x) \, |\tilde{w}_i|^{\chi_i} \, dx \leq \|\psi_i\|_{L^{k_i}} \, \|\tilde{w}_i\|_{L^{\tilde{q}_i}}^{\chi_i} \, .
\end{align*}
If instead $\chi_{i} := \beta \, p^*_{i} +(1/k^{\prime}_{i}-\beta) \, \tilde{q}_{i}$ with a $0<\beta < \min\{1/k_i^{\prime},p_i/p_i^*\}$, we invoke
\begin{align*}
	\int_{\Omega} \psi_i(x) \, |\tilde{w}_i|^{\chi_i} \, dx \leq \|\psi_i\|_{L^{k_i}} \, \|\tilde{w}_i\|_{L^{\tilde{q}_i}}^{(1/k_i^{\prime}-\beta)\tilde{ q}_i} \, \| \tilde{w}_i\|_{L^{p_i^*}}^{\beta\, p_i^*} \, .
\end{align*}
With the Sobolev-- and Poincar\'e--inequality, we further bound 
\begin{align*}
	\| \tilde{w}_i\|_{L^{p_i^*}} \leq c_S \, \| \tilde{w}_i\|_{W^{1,p_i}} \leq c_S \, (1+c_P) \, \|\nabla \tilde{w}_i\|_{L^{p_i}} \, .
\end{align*} 
Since $\beta\, p_i^* < p_i$, the generalised Young inequality allow with $\epsilon > 0$ arbitrary to bound 
\begin{align}\label{evidentes-1} 
	\int_{\Omega} \psi_i(x) \, |\tilde{w}_i|^{\chi_i} \, dx \leq 
	\epsilon \, \int_{\Omega} \tilde{\nu}_i\, |\nabla \tilde{w}_i|^{p_i} \, dx +c(p_i,\beta) \, \Big(\Big[\frac{c_S \, c_P}{(\epsilon\,\tilde{\nu}_i)^{\frac{1}{p_i}}}\Big]^{\beta\, p_i^*} \,\|\psi_i\|_{L^{k_i}} \, \|\tilde{w}_i\|_{L^{\tilde{q}_i}}^{(\frac{1}{k_i^{\prime}}-\beta) \tilde{q}_i}\Big)^{\frac{p_i}{p_i-\beta\,  p_i^*}} \, .
\end{align} 
Since $\tilde{r}_i \leq \chi_i$ implies that $|\tilde{w}_i|^{\tilde{r}_i} \leq 2^{\tilde{r}_i}\, (1+ |\tilde{w}_i|^{\chi_i})$, we hence can conclude that
\begin{align}\label{evidentes} 
	& \int_{\Omega} \psi_i(x) \, |\tilde{w}_i|^{\tilde{r}_i} \, dx \leq 
	\epsilon \, \int_{\Omega} \tilde{\nu}_i\, |\nabla \tilde{w}_i|^{p_i} \, dx + c(\epsilon,\tilde{r}_i,p_i,\tilde{q}_i, \, \|\psi_i\|_{L^{k_i}}, \, \|\tilde{w}_i\|_{L^{\tilde{q}_i}})\, .
\end{align}
By similar means, since $s_i \leq p_i/K_i^{\prime}$ we find that
\begin{align*}
	& \int_{\Omega} \phi_i(x) \, |\nabla \rho_i|^{s_i} \, dx \leq \epsilon \, \int_{\Omega} \nu_i\, |\nabla \rho_i|^{p_i} \, dx + c(\epsilon,s_i,p_i,\nu_i^{-1}) \, \|\phi_i\|_{L^{K_i}}^{K_i}
	\, . % \,  \int_{\Omega} \Big(\frac{C_0}{\nu_i}\Big)^{\frac{r_i}{p_i-r_i}} \, dx\, \, .
\end{align*}
Next we introduce for ease of writing the abbreviation $z:= \frac{s_{N+1} \, (1-\alpha) - t}{\alpha}$. Observe that for $t > s_{N+1}/p_{N+1}$ we have $z < 0$ for $\delta$ sufficiently small.

For $0 \leq s_{N+1} < p_{N+1}/K_{N+1}^{\prime}$ or for $s_{N+1} = p_{N+1}/K_{N+1}^{\prime}$ and $t > s_{N+1}/p_{N+1} = 1/K_{N+1}^{\prime}$ use of H\"older's inequality yields
\begin{align}\label{rsos}
	& \int_{\Omega}  \phi_{N+1}(x)\, |\nabla \tilde{w}_{N+1}|^{s_{N+1}} \, (1+|\tilde{w}_{N+1}|)^{z} \, dx
	\leq \|\phi_{N+1}\|_{L^{K_{N+1}}} \, \|\nabla w_{N+1}\|_{L^{p_{N+1}}}^{s_{N+1}} \, \|1+|\tilde{w}_{N+1}|\|_{L^{y}}^z\, ,
\end{align}
with $y = zp_{N+1}K_{N+1}^{\prime}/(p_{N+1}-s_{N+1})$, where for $t > s_{N+1}/p_{N+1}$ we can ignore the last term. Otherwise, if $\tilde{q}_{N+1} \geq p_{N+1}^*$, we require $y \leq \tilde{q}_{N+1}$ and with the help of Young's inequality we easily see that
\begin{align}\label{rsos0}
	\int_{\Omega}  \phi_{N+1}(x)\, |\nabla \tilde{w}_{N+1}|^{s_{N+1}} \, (1+|\tilde{w}_{N+1}|)^{z} \, dx
	\leq & \epsilon \, \int_{\Omega} \tilde{\nu}_{N+1} \,  |\nabla \tilde{w}_{N+1}|^{p_{N+1}} \, dx \nonumber\\ & +c(\epsilon,\delta,p_{N+1},s_{N+1},q_{N+1},\|\phi_{N+1}\|_{L^{K_{N+1}}},\|u\|_{L^{q_{N+1}}}) \, .
\end{align}
Invoking the definition of $y$ and $z$, the condition $y \leq \tilde{q}_{N+1}$ is equivalent with
\begin{align*}
	s_{N+1} \, ((1-\alpha) \, p_{N+1}K_{N+1}^\prime + q_{N+1}) \leq q_{N+1} \, p_{N+1} + t \, p_{N+1} \, K_{N+1}^{\prime} \, .
\end{align*}
Since $\alpha(\delta) \rightarrow 1/p_{N+1}^{\prime}$ for $\delta \rightarrow 0$, we can satisfy this condition by choosing $\delta$ small enough, provided that $s_{N+1}$ satisfies the conditions of the Lemma. The case $\tilde{q}_{N+1} < p_{N+1}^*$ is slightly more difficult to verify. We require $y \leq p^*_{N+1}$ and if $y \leq \tilde{q}_{N+1}$ we can argue as just seen. Otherwise we bound
\begin{align*}
	\|\tilde{w}_{N+1}\|_{L^y} \leq \|\tilde{w}_{N+1}\|_{L^{p_{N+1}^*}}^{\theta} \, \|\tilde{w}_{N+1}\|_{L^{\tilde{q}_{N+1}}}^{1-\theta}, \quad \theta := \frac{p^*_{N+1} (y-\tilde{q}_{N+1})}{y \, (p^*_{N+1}-\tilde{q}_{N+1})} \, .
\end{align*}
Then, invoking elementary inequalities and \eqref{rsos} yields
\begin{align*}%\label{rsos}
	& \int_{\Omega}  \phi_{N+1}(x)\, |\nabla \tilde{w}_{N+1}|^{s_{N+1}} \, (1+|\tilde{w}_{N+1}|)^{z} \, dx \\
	&	\leq \|\phi_{N+1}\|_{L^{K_{N+1}}} \, (|\Omega|^{\frac{1}{y}}+ \|\tilde{w}_{N+1}\|_{L^{p_{N+1}^*}}^{\theta} \, \|\tilde{w}_{N+1}\|_{L^{\tilde{q}_{N+1}}}^{1-\theta})^z\, \|\nabla w_{N+1}\|_{L^{p_{N+1}}}^{s_{N+1}} \\
	&	\leq c(\delta,p_{N+1},s_{N+1},q_{N+1},\|\phi_{N+1}\|_{L^{K_{N+1}}},\|u\|_{L^{q_{N+1}}},|\Omega|) \, ( \|\nabla w_{N+1}\|_{L^{p_{N+1}}}^{s_{N+1}}+ \|\nabla w_{N+1}\|_{L^{p_{N+1}}}^{s_{N+1}+z\theta} )\, .
\end{align*}
With the help of Young's inequality we can obtain the structure of estimate \eqref{rsos0} if $s_{N+1} + z\theta < p_{N+1}$, which can be verified for $\delta$ small enough under the assumptions of the Lemma.
% 
% 
% \epsilon \, \int_{\Omega} \tilde{\nu}_{N+1} \,  |\nabla \tilde{w}_{N+1}|^{p_{N+1}} \, dx\nonumber\\
% & \qquad + c\,\int_{\Omega} \Big(\frac{(1+|\tilde{w}_{N+1}|)^{\frac{s_{N+1} \, (1-\alpha) - t}{\alpha}}}{\epsilon\tilde{\nu}_{N+1}}  \, \Big)^{\frac{p_{N+1}}{p_{N+1}-s_{N+1}}} \, \phi_{N+1}^{\frac{p_{N+1}}{p_{N+1}-s_{N+1}}}\, dx\nonumber\\
% \leq& \epsilon \, \int_{\Omega} \tilde{\nu}_{N+1} \,  |\nabla \tilde{w}_{N+1}|^{p_{N+1}} \, dx +c(\epsilon,\tilde{\nu}_{N+1}) \, \|\phi_{N+1}\|_{L^{K_{N+1}}}^{\frac{p_{N+1}}{p_{N+1}-s_{N+1}}}\,  \Big(\int_{\Omega} (|\tilde{w}_{N+1}|^{y} + 1) \, dx\Big)^{1-\frac{p_{N+1}}{K_{N+1}(p_{N+1}-s_{N+1})}}  \, ,
%\end{align*}
%where we let $y(\delta) := \frac{[(1+\delta) \, s_{N+1}-t\, p_{N+1}]}{\alpha(\delta)\, (p_{N+1}-s_{N+1}-p_{N+1}/K_{N+1})}$. If $y(\delta) \leq \hat{r}_{\rho}(p_{N+1}, \, \tilde{q}_{N+1}, \, \infty)$, we can argue as in \eqref{evidentes} to control $\int_{\Omega} (|\tilde{w}_{N+1}|^{y} + 1) \, dx$ and 
We spare the trivial verification, and we thus finally obtain that
\begin{align}\label{quatre}
	&\int_{\Omega}  \phi_{N+1}(x)\, |\nabla \tilde{w}_{N+1}|^{s_{N+1}} \, (1+|\tilde{w}_{N+1}|)^{\frac{s_{N+1} \, (1-\alpha) - t}{\alpha}} \, dx
	\leq \epsilon \, \int_{\Omega} \tilde{\nu}_{N+1} \,  |\nabla \tilde{w}_{N+1}|^{p_{N+1}} \, dx\nonumber\\
	& \qquad +  c(\epsilon,\delta,p_{N+1},s_{N+1},q_{N+1},\|\phi_{N+1}\|_{L^{K_{N+1}}},\|u\|_{L^{q_{N+1}}})  \, .
\end{align}
%We can verify that $y(0) < \hat{r}_{\rho}(p_{N+1}, \, q_{N+1}\, p^{\prime}_{N+1}, \, \infty)$ if $s_{N+1}$ satisfies the conditions in the statement. Hence we can find $\delta > 0 $ small enough such that $y(\delta) \leq \alpha(\delta)\, \hat{r}_{\rho}(p_{N+1}, \,) q_{N+1}/\alpha(\delta), \, \infty)$. Thus \eqref{quatre} is valid under the assumptions of the Lemma. 

%where we have to restrict $\delta < p_{N+1}-1-K_{N+1}^{'} \, (s_{N+1}-t)$.

Combining \eqref{evidentes} and \eqref{quatre} and choosing $\epsilon < 1/8$ therein, we obtain that
\begin{align} \label{provi}
	\int_{\Omega} \sum_{i=1}^{N+1} \tilde{\nu}_i\,  |\nabla \tilde{w}_i|^{\tilde{p}_i} \, dx \leq C_{\delta}\, ,
\end{align}
for all $0<\delta < \min\{\delta_0, \, \delta_1\}$ with a number $\delta_1$ determined by the exponents ${\bf r}$, ${\bf s}$, ${\bf k}$, ${\bf K}$, ${\bf q}$ and $t$.
%\begin{align*}
% 0 < \delta < \min\Big\{p_{N+1}-1-k_{N+1}'r_{N+1}, \,p_{N+1}-1-K_{N+1}'s_{N+1},\, \delta_0\Big\} \, .
%\end{align*}
Finally we want to remove the weight affecting $\nabla u$ in \eqref{provi}. Let 
%$1 \leq \sigma < \frac{n}{n-1} \, (p_{N+1}-1)$ 
\begin{align*}
	1 \leq \sigma < p_{N+1} \, \frac{\max\{q_{N+1}\, p_{N+1}^{\prime},p_{N+1}^*\}}{p_{N+1}^{\prime}+\max\{q_{N+1}\, p_{N+1}^{\prime},p_{N+1}^*\}} \, .
\end{align*}
be arbitrary. This guarantees that we can choose a positive $\delta$ such that with $m = \max\{q_{N+1}/\alpha(\delta),p_{N+1}^*\}$
\begin{align}\label{deltamaxrestri}
	0 < \delta < \frac{m \, (p_{N+1}-1) - p_{N+1} \, \sigma}{m \, (p_{N+1}-\sigma) + p_{N+1}\sigma} \, \quad \text{ which implies } \quad \frac{ (1+\delta) \, \sigma}{\alpha \, (p_{N+1}-\sigma)} \leq m \, .
\end{align}
Then, by H\"older's inequality
\begin{align}\label{rsosigma} 
	\int_{\Omega} |\nabla u|^{\sigma} \, dx \leq  \left(\int_{\Omega} \Big(\frac{1}{1+|u|}\Big)^{1+\delta}\, |\nabla u|^{p_{N+1}}  \, dx\right)^{\frac{\sigma}{p_{N+1}}} \, \left(\int_{\Omega} (1+|u|)^{\frac{(1+\delta) \, \sigma}{p_{N+1}-\sigma}} \, dx \right)^{1-\frac{\sigma}{p_{N+1}}}\, \nonumber\\
	\leq (\tilde{\nu}_{N+1})^{-\frac{p_{N+1}}{\sigma}} \, \left(\int_{\Omega}  \tilde{\nu}_{N+1} \, |\nabla \tilde{w}_{N+1}|^{p_{N+1}}  \, dx\right)^{\frac{\sigma}{p_{N+1}}} \, \left(\int_{\Omega} (1+|\tilde{w}_{N+1}|)^{\frac{(1+\delta) \, \sigma}{p_{N+1}-\sigma}} \, dx \right)^{1-\frac{\sigma}{p_{N+1}}}
\end{align}
Under the condition \eqref{deltamaxrestri} it follows from \eqref{Tildeuprops} that
\begin{align*} 
	\int_{\Omega} (1+|\tilde{w}_{N+1}|)^{\frac{(1+\delta)\, \sigma}{p_{N+1}-\sigma}} \, dx \leq & c(\Omega) + \int_{\Omega} |\tilde{w}_{N+1}|^{m} \, dx \\
	\leq & C_1 + C_2 \,  \|\nabla\tilde{w}_{N+1}\|_{L^{p_{N+1}}}^{p_{N+1}^*} + C_3 \,  \|\nabla \tilde{w}_{N+1}\|_{L^{\frac{q_{N+1}}{\alpha}}}^{\frac{q_{N+1}}{\alpha}}  \, .
\end{align*}
Thus, \eqref{provi} allows to estimate $\|\nabla u\|_{L^{\sigma}(\Omega)}$.\\

Next we discuss the reduction of the mixed boundary value problem to a reference configuration using the assumption (C1).

Consider $x^0 \in \partial S$ fixed.
To simplify the notation we write $T = T^{x^0}_S$, and we denote by $\hat{\Omega}$ the image $T(\Omega_{R_0}(x^0))$ with elements $\hat{x}$ such that $\hat{x}_n < 0$ and $\hat{x}_{n-1} < 0$, by $\hat{S}_{1,R_0}$ the image $T(S\cap B_{R_0}(x^0))$ of the Dirichlet part which is planar and contained in the plane $\hat{x}_n = 0$, and by $\hat{S}_{2,R_0}$ the image $T((\partial \Omega \setminus \overline{S})\cap B_{R_0}(x^0))$ of the Neumann part of the boundary, which is a subset of the plane $\hat{x}_{n-1} = 0$.

Consider next $\rho_0 := R_0/L$ (cf.\ (C1), \eqref{Squatre}). The Lipschitz continuity of the inverse ensures that $T^{-1}(\hat{\Omega}_{\rho_0}(0)) \subseteq \Omega_{R_0}(x^0)$, where $\hat{\Omega}_{\rho_0}(0):=\hat{\Omega} \cap \hat{B}_{\rho_0}(0)$. In this way, the properties of the diffeomorphism $T$ ensure that $\partial \hat{\Omega} \cap \hat{B}_{\rho_0}(0)$ is decomposed into two perpendicular planar surfaces: $\hat{S}_{1,\rho_0}$ contained in the plane $\hat{x}_n = 0$ and $\hat{S}_{2,\rho_0}$ in $\hat{x}_{n-1} = 0$.

For ease of writing we let $\Psi := T^{-1}$. If $v = v^0$ solves \eqref{rhoihoelder1}, \eqref{rhoihoelder2} in $\Omega_{R_0}(x^0)$ then we define $\hat{\rho}_i \in W^{1,p_i}(\hat{\Omega}_{\rho_0}(0))$ via $\hat{\rho}_i(\hat{x}) := \rho_i(\Psi(\hat{x}))$, and we let
\begin{align*}
	\widehat{a}_i(\hat{x}, \, z) := \text{det}(D_{\hat{x}}\Psi(\hat{x})) \, (D_{\hat{x}}\Psi(\hat{x}))^{-1} \, \widetilde{a}_i(\Psi(\hat{x}), \, (D_{\hat{x}}\Psi(\hat{x}))^{-\sf T} \, z) \, \quad \text{ for } \hat{x}\in\hat{\Omega}, \, z\in \mathbb{R}^n \, .
\end{align*}
We employ (A1) to verify that
\begin{align*}
	\left(\frac{1}{LM}\right)^2 \, \nu_i \, |z|^{p_i} \leq 	\widehat{a}_i(\hat{x}, \, z) \cdot z \leq 	(L\, M)^2 \, \mu_i \, |z|^{p_i} \quad \text{ for all } \hat{x} \in \hat{\Omega}, \, z \in \mathbb{R}^n\, .
\end{align*}
Moreover, the function $\hat{\rho}_i$ solves in the weak sense
\begin{gather}\begin{split}\label{hatauxiv}
		& -\divv_{\hat{x}} (\widehat{a}_i(\hat{x}, \, \nabla_{\hat{x}} \hat{\rho}_i) + \hat{b}_i) + \hat{B}_i = 0 \text{ in } \hat{\Omega}_{\rho_0}(0) \,, \\
		& \hat{\rho}_i = \hat{\rho}_i^S \text{ on }\hat{S}_{1,\rho_0} \, ,\qquad e^{n-1} \cdot (\widehat{a}_i(\hat{x}, \, \nabla_{\hat{x}} \hat{\rho}_i) + \hat{b}_i) = 0 \text{ on } \hat{S}_{2,\rho_0}\, .
	\end{split}
\end{gather}
%Here $\hat{S} = T^{x_0}(B_{R_0}(x^0) \cap S)$.
%with homogeneous Neumann conditions on $B_{R}(x^0) \cap \partial \Omega \setminus S$.
Let us consider $R_1 := \rho_0/M = R_0/(LM)$, which yields $T(\Omega_{R_1}(x^0)) \subset \hat{B}_{MR_{1}}(0)$ with $MR_1 = \rho_0$. Suppose that we have proved the H\"older continuity in for $\hat{\rho}_i$ in $\hat{\Omega}_{\rho_0}(0)$ with exponent $\lambda_i$, then 
\begin{align*}
	|\rho_i(x_1)-\rho_i(x_2)| \leq &   |\hat{\rho}_i(T(x_1))-\hat{\rho}_i(T(x_2))| \leq [\hat{\rho}_i]_{\hat{\Omega}_{\rho_0},\lambda_i} \, |T(x_1)-T(x_2)|^{\lambda_i} \\
	\leq & M^{\lambda_i} \, [\hat{\rho}_i]_{\hat{\Omega}_{\rho_0},\lambda_i} \, |x_1-x_2|^{\lambda_i}  \, .
\end{align*}
Hence, this preliminary consideration shows that we indeed loose no generality in proving the H\"older bound on a reference configuration.\\

%Thus we now start with the original notations but assume that the boundary is already in the reference configuration near $x^0 = 0$, hence $S_{1,R_0}$ and $S_{2,R_0}$ are planar and orthogonal to each other, while $\Omega_R(x^0)$ is the intersection of $B_R(0)$ with the two half spaces $x_n <0$ and $x_{n-1} < 0$.
%*********************
Next we recall the Sobolev embedding result needed throughout the paper.
\begin{prop}\label{lieberman}
Let $\Omega$ be of class $\mathscr{C}^{0,1}$. Let $1 < p < n$ and $\mu$ be a Radon measure such that with a $0< \theta \leq 1$ we have $|\mu|(\Omega \cap B_\rho(x)) \leq M \, \rho^{n-p+\theta \, p}$ for all $x \in \overline{\Omega}$ and $\rho >0$. Then there is $C = C(n,p,\theta,\Omega,S)$ such that for all $u \in W^{1,p}_S(\Omega)$ (there is $C = C(n,p,\theta,\Omega)$ such that for all $u \in W^{1,p}(\Omega)$ with $\int_{\Omega} u \, dx = 0$)  	
\begin{align*}
\left(\int_{\Omega} |u|^{\kappa \, p} \, d\mu\right)^{\frac{1}{\kappa}} \leq C \, M^{\frac{1}{\kappa}} \, \int_{\Omega} |\nabla u|^p \, dx \quad \text{ with } \quad \kappa = 1 + \frac{\theta \, p}{n-p} \, .
\end{align*}
\end{prop}
\begin{proof}
Let $R_0$ and $y \in \Omega$ be such that $\Omega \subset\!\! B_{R_0}(y)$. Since $\Omega$ is of class $\mathscr{C}^{0,1}$, there is a linear bounded extension operator $E: \, W^{1,p}(\Omega) \rightarrow W^{1,p}_0(B_{R_0}(y))$. We let $\tilde{\mu}$ be the trivial extension $\tilde{\mu}(U) := \mu(U \cap \overline{\Omega})$ for $U \subset B_{R_0(y)}$. Suppose that $B_{\rho}(x) \subset B_{R_0}(y)$. Then, if $B_{\rho}(x)$ intersects $\overline{\Omega}$, we can trivially show that $B_{\rho}(x) \subset B_{2\rho}(x^0)$ for any point $x^0$ of the intersection. This implies that $$|\tilde{\mu}|(B_\rho(x)) \leq |\mu|(\Omega_{2\rho}(x^0)) \leq 2^n \, M \, \rho^{n-p+\theta \, p} \, ,$$ where the assumption was used.

We apply the Lemma 1.1 in \cite{MR1233190}, showing first that
\begin{align*}
	\left(\int_{B_{R_0}(y)} |E \, u|^{\kappa \, p} \, d\mu\right)^{\frac{1}{\kappa}} \leq C_0 \, (2^n \, M)^{\frac{1}{\kappa}} \, \int_{B_{R_0}(y)} |\nabla E\, u|^p \, dx \quad \text{ with } \quad \kappa = 1 + \frac{\theta \, p}{n-p} \, ,
\end{align*}
where $C_0$ depends only on $n$, $p$ and $\theta$. Hence, using the properties of the extension operator, we obtain that
\begin{align*}
	\left(\int_{\Omega} |u|^{\kappa \, p} \, d\mu\right)^{\frac{1}{\kappa}} \leq C_0 \, (2^n \, M)^{\frac{1}{\kappa}} \, \|E\|^p \,  \|u\|_{W^{1,p}(\Omega)}  \, .
\end{align*}
It remains to observe that for $u \in W^{1,p}_S$, we can estimate $\|u\|_{L^p(\Omega)} \leq c_1(p,S,\Omega) \, \|\nabla u\|_{L^p}$, and the claim follows with $C = C_0\, 2^{n/\kappa} \, \|E\|^p \, (1+c_1^{p})$. For the case $u \in W^{1,p}(\Omega)$ such that $\int_{\Omega} u \, dx = 0$ we replace the last inequality by the Poincar\'e inequality for Sobolev functions with average zero.
\end{proof}

\end{document}